\documentclass[trsc,nonblindrev]{informs3} 

\OneAndAHalfSpacedXI 

\usepackage{natbib}
 \bibpunct[, ]{(}{)}{,}{a}{}{,}%
 \def\bibfont{\small}%
\DeclareSymbolFont{matha}{OML}{txmi}{m}{it}
\DeclareMathSymbol{\varv}{\mathord}{matha}{118}

\usepackage{bbm}

\def\R{\mathbb{R}}

\def \-> {\rightarrow}

\def \cotwo {CO$_2$ }

\usepackage{xcolor}
\newenvironment{rev}{\color{black}}{}
\newenvironment{rev2}{\color{black}}{}

\usepackage{enumerate}
\usepackage{enumitem}
\usepackage{tabularx,booktabs,mathrsfs,setspace}%
\usepackage{mathtools}
\usepackage{multirow}
\usepackage{rotating}
\usepackage{float}
\usepackage{pgfplots}
\pgfplotsset{compat=1.9}
\usepackage{csvsimple}
\usepackage{makecell}
\usepackage{multirow}
\usepackage{adjustbox}
\usepackage[caption=false, position=top]{subfig}
\usepackage{alphalph}

\usepackage{mathtools}
\usepackage{comment}
\usepackage{lscape}
\usepackage{longtable}
\usepackage{xltabular}

\TheoremsNumberedThrough     

\EquationsNumberedThrough    

\begin{document}


\RUNAUTHOR{Moradi et al.} 

\RUNTITLE{Load Asymptotics and Dynamic Speed Optimization for the Greenest Path Problem}

\TITLE{Load Asymptotics and Dynamic Speed Optimization for the Greenest Path Problem: \\
A Comprehensive Analysis}

\ARTICLEAUTHORS{%
\AUTHOR{Poulad Moradi, Joachim Arts}
\AFF{Luxembourg Centre for Logistics and Supply Chain Management, University of Luxembourg, Luxembourg City, Luxembourg, 6, rue Richard Coudenhove-Kalergi L-1359, \EMAIL{\{poulad.moradi, joachim.arts\}@uni.lu}}
\AUTHOR{Josué C. Velázquez-Martínez}
\AFF{Center for Transportation and Logistics, Massachusetts Institute of Technology, Cambridge, MA, USA,\\1 Amherst Street, MA 02142,  \EMAIL{josuevm@mit.edu} \URL{}}
} 

\ABSTRACT{%
We study the effect of using high-resolution elevation data on the selection of the most fuel-efficient (greenest) path for different trucks in various urban environments. We adapt a variant of the Comprehensive Modal Emission Model (CMEM) to show that the optimal speed and the greenest path are slope dependent (dynamic). 
When there are no elevation changes in a road network, the most fuel-efficient path is the shortest path with a constant (static) optimal speed throughout. 
However, if the network is not flat, then the shortest path is not necessarily the greenest path, and the optimal driving speed is dynamic. 
We prove that the greenest path converges to an asymptotic greenest path as the payload approaches infinity and that this limiting path is attained for a finite load.
In a set of extensive numerical experiments, we benchmark the \cotwo emissions reduction of our dynamic speed and the greenest path policies against policies that ignore elevation data. We use the geo-spatial data of 25 major cities across 6 continents.
We observe numerically that the greenest path quickly diverges from the shortest path and attains the asymptotic greenest path even for moderate payloads.
Based on an analysis of variance, the main determinants of the \cotwo emissions reduction potential are the variation of the road gradients along the shortest path as well as the relative elevation of the source from the target. 
\begin{rev}
Using speed data estimates for rush hour in New York City, we test \cotwo emissions reduction by comparing the greenest paths with optimized speeds against the fastest paths with traffic speed.
We observe that selecting the greenest paths instead of the fastest paths can significantly reduce \cotwo emissions.
Additionally, our results show that while speed optimization on uphill arcs can significantly help \cotwo reduction, the potential to leverage gravity for acceleration on downhill arcs is limited due to traffic congestion.
\end{rev}
}%



\maketitle

%


\section{Introduction}

The transportation sector is one of the largest sources of anthropogenic \cotwo emissions, as attested by the \cite{IPCCWGI2021}, \cite{GHGInventory}, and the \cite{EEA2021}. 
In 2020, 36.3\% of U.S. \cotwo emissions from fossil fuel combustion came from the transportation sector, of which 45.2\% was generated by heavy-, medium-, and light-duty trucks \citep{GHGInventory}. 
Similarly, the transportation sector accounted for 22\% of the EU's \cotwo emissions in 2020 \citep{EEA2021}. 
Accordingly, there has been considerable attention on reducing \cotwo emissions through ``green routing''; see e.g. \cite{Demir2012, Scora2015, Raeesi2019}. The objective to reduce \cotwo emissions in transportation aligns with efforts to reduce fuel expenditure. 
The reduction of fuel consumption has become imperative as fuel increases in price and volatility due to recent geopolitical events, namely the Russian invasion of Ukraine \citep{WLJ2022}.

Road gradient and vehicle speed are two major factors that influence the carbon footprint of a diesel truck \citep{Demir2014}. 
\cite{Demir2011} demonstrate through numerical analysis that a medium-duty truck may consume an additional six liters of diesel per 100 kilometer while traveling up a hill with a 1\% gradient. The same study also shows that increasing the speed of an empty medium-duty truck from 50 km/h to 100 km/h can raise fuel consumption by more than 3\% on a level path. 
Gravity is an important factor in finding the most efficient path between two points. Johann Bernoulli posed such a problem as early as 1696, in which a path was deemed efficient if the travel time was minimized and only gravity could be used to accelerate. 
The solution to this problem gave rise to so-called brachistochrone curves, which differ from the shortest path between two points. Gravity and speed interact when finding the greenest (or most fuel-efficient) route between two points in a road network. 
The aim of the present paper is to provide a thorough analysis of the difference between shortest paths and greenest paths as a function of speed and vehicle type for a large variety of geographic settings.

In principle, empirical methods are the most precise way of measuring carbon emissions associated with traversing a road with a certain vehicle at a certain speed. 
Unfortunately, it is not practical to empirically find carbon emissions for all roads, speeds, and vehicle types as well as many other parameters (e.g. road surface type) that affect fuel efficiency. Hence, several \cotwo emissions models for trucks have been proposed in literature. 
\cite{Demir2014} offers a summary of these models. The Comprehensive Modal Emission Model (CMEM) is an instantaneous emissions modeling approach that was introduced by \cite{Barth2005, Scora2006, Boriboonsomsin2009}. \cite{Bektas2011} and \cite{Demir2012} present a simplified variant of CMEM that is differentiable with respect to speed. This model is convenient in practical applications. \cite{Rao2016} and \cite{Brunner2021} show that this model can be made more realistic for cases where a vehicle travels downhill. Their modification of the CMEM, unfortunately, renders it no longer differentiable at all speeds. 
Over the past decade, CMEM has been the prevalent emissions model utilized in green/pollution vehicle routing problems \citep[e.g.][]{Bektas2011, Franceschetti2013, Huang2017, Xiao2020}.

We call an optimization problem that seeks a path between an origin and destination a {\em path selection} problem.
In this paper, we focus on the selection of the greenest (most fuel-efficient) path. The greenest path is the path with the least \cotwo emissions. Some authors also call this the eco-friendly path \citep[e.g.][]{Scora2015,Andersen2013,Boriboonsomsin2009,Schroder2019}.  
Path selection is the backbone of a multitude of transport-based supply chain problems, from strategic supply chain network design to operational vehicle routing problems. The complexity of transportation problems forces many solution approaches to use path selection as a
\begin{rev}
    pre-processing
    \end{rev}
activity. 
It is common to use either the shortest or the fastest path in this pre-rocessing step. The implicit assumption is that these paths are also the greenest. In this paper, we show that the actual topology of urban road networks requires that we consider the greenest path selection as a part of the main optimization problem, e.g. vehicle routing problem (VRP) or supply chain network design (SCND).

The development of Geographic Information Systems (GIS) have made high-resolution geospatial data available at low cost. It is not sufficient to only consider the elevation of the origin and destination of a path. Rather, for any path, elevation along different sections of a path determine whether gravity increases or decreases the amount of fuel needed for travel. Thus, detailed elevation data of each segment of a possible path is required to find the greenest path. Furthermore, the slope along different segments of a path also determines the most fuel efficient speed along each segment of a path. 

In this paper, we show that the most fuel efficient speed will change along different segments of any path. Thus, dynamic speed optimization is important to find the greenest path between any origin and destination. The greenest path also depends on the payload of a vehicle. We prove that the greenest path converges to an asymptotic greenest path as the payload approaches infinity and that this limiting path is attained for a finite payload. Our results are illustrated through numerical experiments. 
These experiments consider a setting wherein a logistics service provider seeks to reduce the \cotwo emissions of their transport operations. The company's fleet consists of heavy-, medium-, and light-duty trucks that operate in an urban environment. We use the modified CMEM proposed by \cite{Brunner2021} and focus our analysis on the effects that road gradient, speed, payload, and truck type have on \cotwo emissions. 

We use an extensive numerical study to provide statistical answers to empirical research questions listed at the end of this section. We utilize the real road network and elevation data of 25 cities across six continents. It is worth noting that the closest paper to our work in terms of the \cotwo emissions model is \cite{Brunner2021}. Apart from the differences of our objective functions, the main differences between our study and \cite{Brunner2021} are twofold. Firstly, \cite{Brunner2021} base their analysis on the static speed policy along different segments of a path. We show that a static speed policy can be suboptimal in terms of \cotwo emissions for traversing a path in a city with uneven topography. In addition, we demonstrate that the speed policy influences which path is the greenest. Secondly, \cite{Brunner2021} solve the path selection problem as a
\begin{rev}
    pre-processing
\end{rev}
activity for their main VRP problem. Moving the path selection to a
\begin{rev}
    pre-processing
\end{rev}
step forces them to consider fixed loads and speeds. By contrast, we consider dynamic speed optimization, and study asymptotics greenest paths as payloads increase.
\begin{rev}
We utilize estimated traffic speeds during rush hour for a large subgraph of New York City's road network to study the potential \cotwo reduction by choosing the greenest paths and optimizing speed instead of taking the fastest paths.
We also examine the increased travel duration on the greenest paths, as well as the convergence to the asymptotic greenest path when traffic congestion occurs.
\end{rev}

The main contributions of this paper are listed below:
\begin{enumerate}
\item We show that the greenest path is speed and payload dependent for accurate emission models. We provide a tractable algorithm to optimize the path and the speed jointly, where the speed varies along the path.
\item We show that the greenest path converges to an asymptotic greenest path when the load becomes large and that this path is attained for a finite load in Section \ref{subsect:asymcal}. In Sections \ref{subsec:trajectories} and \ref{subsec:AP_num}, we show that this convergence happens relatively quickly in practice. We also show, in Section \ref{subsec:trajectories}, that the greenest path for a slope-dependent optimal speed policy is quite similar to the one associated with the static speed policy of \cite{Demir2012}, yet significantly different from the shortest path.
\item \label{questions:title} We conduct an extensive numerical study with data from 25 cities over 6 continents and over 3 million origin destination pairs. We use detailed elevation data from \cite{USGS}'s SRTM 1 Arc-Second Global data set. This thorough study allows us to answer the following research questions:
\begin{enumerate}
    \item \label{questions:overalReduction} How much \cotwo emissions can dynamic speed optimization and green path selection reduce jointly? What are the effects of truck type, payload, and city on the carbon reduction potential?
    \item What is the marginal contribution of speed optimization and path optimization in the reduction of \cotwo emissions?
    \item How different is the slope-dependent greenest path from the shortest path (which is the slope-disregarding path)?
    \item \label{questions:speedImpact} What is the impact of the speed policy on the greenest path?
    \item In what settings are the integration of elevation data in path selection most valuable? 
\end{enumerate}
\begin{rev}
\item We conduct an extensive numerical study with road network, elevation, and speed data of New York City for more than 20 thousands origin destination pairs. The traffic speed estimates are collected from Google's Distance Matrix API. With this study we answer to the questions \ref{questions:overalReduction} through \ref{questions:speedImpact} when shortest path is replaced by the fastest path. We also study the increased time in traffic when the greenest path and speed optimization are decided.
\end{rev}

\end{enumerate}

The rest of this paper is organized as follows.
We review the related literature in Section \ref{sec:literature}.
Section \ref{section:model} describes the mathematical model used in this study and the policies that can minimize \cotwo emissions. Section \ref{section:NumericalExperiments} provides the setting and results of the extensive numerical studies
\begin{rev}
under free-flow conditions.
We present the setting and results of our numerical study under traffic congestion in Section \ref{section:NumericalExperimentTraffic}.
\end{rev}
Finally, we offer the conclusions and final remarks in Section \ref{section:conclusion}.

\section{Literature review}
\label{sec:literature}
Green transportation has been studied extensively over various decision-making settings. \cite{Asghari2021,Moghdani2021,Demir2014} give reviews on the most important recent literature on the green VRP. Additionally, \cite{Waltho2019} reviews pivotal studies in the field of green SCND from 2010 to 2017. In most of the main stream green VRP and SCND, the path between every two nodes of interest is computed as a \begin{rev}
    pre-processing
\end{rev}step \citep[e.g.][]{Demir2012}. This has been partially relaxed for the VRP by \cite{Behnke2017}. 
In other words, road networks are reduced to distances between origin and destination pairs to simplify later computations. The implicit assumption is that distances or travel times are the main drivers of costs and/or emissions. This paper extensively studies to which extent this implicit assumption is tenable.

A large body of work in the field of green transportation relies on macroscopic (average aggregate), microscopic (instantaneous) fuel consumption models, or a combination of both. \cite{Demir2014} and \cite{Zhou2016} provide an extensive review of fuel consumption models. 
A number of studies, including \cite{Boriboonsomsin2012,Scora2015}, and \cite{Ericsson2006}, estimate the \cotwo emissions of a specific vehicle based on the measurement of that vehicle.
\cite{Demir2014} explain the main factors that influence fuel consumption in road freight transportation. Among the pertinent determinants for the case of the greenest path are road gradient, speed, truck type, and payload.

The path optimization under environmental consideration (the greenest path) has been explored over the past two decades. This problem can be formulated based on a variant of the shortest path algorithm of \cite{Dijkstra1959} to minimize the total fuel consumption of a vehicle between two nodes. 
\cite{Ericsson2006} studies the \cotwo emissions of light-duty cars by using a navigation system that computes the greenest path based on in-vehicle data and traffic information in Lund, Sweden. They conclude that selecting the greenest path can reduce fuel consumption by 4\% on average in Lund. \cite{Boriboonsomsin2012} presents an Eco-Routing Navigation System (EFNav) as a framework to integrate GIS and traffic data with emissions model estimates to compute eco-friendly paths for light vehicles. 
\cite{Scora2015} extend the EFNav model to heavy-duty trucks (EFNav-HDT) and conduct a numerical study to test the benefits of EFNav-HDT across different vehicle weights in Southern California. \cite{Scora2015} provide excellent insights into the specifications of the greenest path for trucks.
Both \cite{Boriboonsomsin2012} and \cite{Scora2015} base their studies on the CMEM model and estimate the energy/emissions model using linear regression over data from actual measurements. \cite{Boriboonsomsin2012} take advantage of a logarithmic transformation and \cite{Scora2015} use a minimum fuel cutoff point to avoid negative fuel consumption results.
\cite{Andersen2013} take advantage of free road network data, such as OpenStreetMap, and use Controller Area Network (CAN bus) data to compute the greenest path by assigning weights to the different segments of the network. Since this work does not rely on a fuel consumption model, it is very accurate for the paths and vehicles for which fuel consumption data is available, but it does not transfer to other settings without the collection of a large amount of data in that setting.
\cite{Pamucar2016} utilize a similar approach and include other negative externalities associated with transportation, such as noise, land use, and pollutants other than \cotwo.
\cite{Schroder2019} consider a Digital Elevation Model and Copert III emissions model to compute the greenest path.
\cite{Dundar2022} propose an approach to increase the resolution of the road network and compute the fuel consumption over along a path more accurately.  

Speed optimization as a means to reduce the emissions and driving costs was first introduced by \cite{Demir2012}. \cite{Franceschetti2013} present a speed optimization technique that can also be used for traffic congestion. Both of these works, as well as many other well-cited papers, such as \cite{Lai2021}, rely on the CMEM model of \cite{Demir2011}, which results in negative fuel consumption over many downhill paths \citep{Brunner2021}. \cite{Brunner2021} modify the fuel consumption model, yet only consider a constant travel speed. Some papers consider the \begin{rev}fastest path or the emissions minimizing path\end{rev} under dynamic speeds induced by congestion \citep[e.g.][]{Ehmke2016a, Ehmke2016b, Huang2017, Ehmke2018}. 

In our paper, we consider the modified CMEM \citep{Brunner2021}. We explore the individual and combined effects of elevation, speed optimization, truck type, payload, and characteristic city topography on \cotwo emissions reduction and the greenest path policies. This paper, is the first paper to provide asymptotic results for a path selection problem and the the greenest path problem in particular.

\section{Model Description}
\label{section:model}
In this section, we introduce the notations (Section \ref{subsec:net-notation}) and mathematical foundations of our research, including the \cotwo emissions models (Section \ref{subsec:emissionsmodels}) together with the optimal speed policies (Section \ref{subsec:speedanalysis}). We formally introduce the greenest paths between two locations in a city road network and discuss how optimal speed policies complicate the computation of the greenest path (Section \ref{subsect:gpp}). We study the asymptotic behavior of the greenest path when the payload increases (Section \ref{subsect:asymcal}).
In Section \ref{section:model}, we only consider speed, payload, and/or path (or a single arc) as the explicit arguments of functions, since these three factors are the focus of our analysis in Sections \ref{subsec:speedanalysis}, \ref{subsect:gpp}, and \ref{subsect:asymcal}.

\subsection{City Network and Notations}
\label{subsec:net-notation}
Let a directed graph $\mathcal{G}=(V,A)$ represent the road network of a city, where $V = \{1, \dots, m\}$ is the set of $m$ vertices, the points of interest along the roads (e.g. road intersections), and $A \subseteq V \times V$ is the set of arcs (road segments) that connects the vertices. Any arc $a \in A$ has the following features: the length $\delta:A\to  \mathbb{R}_{++}$, the angle $\theta: A \to \mathbb{R}$, the maximum allowable speed by $\varv^{\max}: A \to  \mathbb{R}_{++}$, and the minimum allowable speed by $\varv^{\min}: A \to  \mathbb{R}_{++}$, where $\mathbb{R}_{++}=\{x \in \mathbb{R} : x>0\}$.
We consider an internal combustion engine truck that traverses an arc $a\in A$ with speed $v \in [\varv^{\min} (a) , \varv^{\max} (a)]$. $v$ is constant along arc $a$, but the speed of the truck can vary on other arcs. 
The truck consumes $f_a$ liters of diesel fuel and produces $e_a$ kilograms of \cotwo to traverse arc $a \in A$ ($v$ will be selected to minimize $f_a$ and $e_a$ according to different emission models).
Notation, including those of truck properties, are listed in Table \ref{tab:notation}.
\begin{table*}[t!]
  \caption{\textsf{Overview of notation.}}
	\fontsize{8pt}{9pt}\selectfont
\scriptsize
  \label{tab:notation}
  \begin{tabularx}{\textwidth}{l X}
    \toprule
    Notation& Description\\
		\midrule
	Sets & \\
        $\mathcal{G}$ 	& Directed multigraph representing the urban road network, $\mathcal{G}=(V,A)$.\\
        $V$			    & Set of vertices of $\mathcal{G}$, $V=\{1, \dots, m\}$, where $m$ is the number of vertices.\\
        $A$             & Set of arcs of $\mathcal{G}$, $A \subseteq V \times V$.\\
        $\Pi$   & Set of all paths between a pair of nodes $n_s, n_t \in A$.\\
        $\Pi_d$   & Subset of $\Pi$ such that for all $a \in \Pi_d$, $tan \theta (a) < -C_r$.\\ \\
    Network Features & \\
        $\delta (a)$      & $\delta: A \to \mathbb{R}_{++}$, length of arc $a \in A$.\\
        $\theta(a)$      & $\theta: A \to \mathbb{R}$, angle of arc $a \in A$.\\
        $h^\prime (a)$    & $h^\prime: A \to \mathbb{R}$, augmented ascent of arc $a \in A$, i.e. $h^{\prime} (a) = \delta(a) \sin \left(\theta(a) + \arctan C_r \right)^+$.\\
        $\varv^{\max} (a)$ & $\varv^{\max}: A \to \mathbb{R}_{++}$, maximum allowable speed for traversing arcs $a \in A$.\\
        $\varv^{\min} (a)$ & $\varv^{\min}: A \to \mathbb{R}_{++}$, minimum allowable speed for traversing arcs $a \in A$.\\ \\
	Parameters &  \\
	    $\xi$           & Fuel-to-air mass ratio. \\
        g             & Gravitational constant (m/s$^2$).\\
        $\rho$          & Air density (kg/m$^3$).\\
        $C_r$           & Coefficient of rolling resistance.\\
        $\eta$          & Efficiency parameter for diesel engines\\
        $\eta_{tf}$     & Vehicle drivetrain efficiency.\\
        $\kappa$        & Heating value of a typical diesel fuel (kJ/g).\\
        $\psi$          & Conversion factor (g/s to L/s).\\
        $w$             & Curb weight (kg).\\
        $L$             & Maximum payload (kg).\\
        $k$             & Engine friction factor (kJ/rev/L).\\
        $N$             & Engine speed (rps).\\
        $D$             & Engine displacement (L).\\
        $C_d$           & Coefficient of aerodynamic drag.\\
        $S$             & Frontal surface area (m$^2$).\\
        $c_e$           & Fuel's \cotwo Emissions Coefficient.\\ \\
	Variables & \\ 
        $l$             & $l \in \mathbb{R}_+$, payload (kg).\\
        $v$             & $v \in \mathbb{R}_{++}$, speed of a truck to traverse arc $a \in A$.\\
        $c_{\varv}$     & Constant speed (Equation \eqref{eq:constant_opt_speed}) that minimizes the standard emissions model (Equation \eqref{eq:cmem}).\\
	    $\varv (a)$        & $\varv: A \to \mathbb{R}_{++}$, speed policy for a truck to traverse arc $a \in V$.\\
        $\varv^d (a,l)$        & $\varv^d: A \times \mathbb{R}_+ \to \mathbb{R}_{++}$, dynamic speed policy on arc $a \in A$ with payload $l$ as per Proposition \ref{prop:opt_speed}.\\
        $\varv^s (a)$        & $\varv^s: A \to \mathbb{R}_{++}$, static speed policy on arc $a \in A$ as per Equation \eqref{eq:constant_opt_speed}.\\
        $\varv^t (a,l)$        & $\varv^t: A \times \mathbb{R}_+ \to \mathbb{R}_{++}$, terminal velocity on arc $a \in A$ with payload $l$ as per Proposition \ref{prop:opt_speed}.\\
        $t_a (v)$       & $t_a: \mathbb{R_{++}} \to \mathbb{R_{++}}$, traveling time on arc $a \in A$ with speed $v \in \mathbb{R_{++}}$.\\
        $\pi$       & $\pi \in \Pi$, path between a pair of nodes $n_s, n_t \in A$.\\
        $\pi^{sp}$  & $\pi^{sp} \in \Pi$, the shortest path between a pair of nodes (Equation \eqref{eq:shortestPath}).\\
        $\pi^g (\varv, l)$  & $\pi^g (\varv, l) \in \Pi$, the most fuel-efficient (greenest) path between a pair of nodes under the speed policy $\varv$ and payload $l$ (Equation \eqref{eq:greenestPath}).\\
        $\pi^{\infty} (\varv)$  & $\pi^{\infty} (\varv) \in \Pi$, the asymptotic greenest path between a pair of nodes under the speed policy $\varv$, i.e. the greenest  path when the payload is arbitrarily large (Proposition \ref{prop:AP}).\\
        $\Tilde{f}_a (v,l)$           & $\Tilde{f}_a: \mathbb{R}_{++} \times \mathbb{R}_+ \to {R}_{++}$, amount of fuel (liter) that a truck consumes for traversing arc $a \in V$ with speed $v$ and payload $l$ under the standard emissions model (Section \ref{subsubsec:standardEmissionsModel}).\\
        $f_a (v,l)$           & $f_a: \mathbb{R}_{++} \times \mathbb{R}_+ \to {R}_{++}$, amount of fuel (liter) that a truck consumes for traversing arc $a \in V$ with speed $v$ and payload $l$ under the improved emissions model (Section \ref{subsubsec:improvedEmissionsModel}).\\
        $\Tilde{e}_a (v,l)$           & $\Tilde{e}_a: \mathbb{R}_{++} \times \mathbb{R}_+ \to {R}_{++}$, amount of \cotwo (kg) that a truck emits for traversing arc $a \in V$ with speed $v$ and payload $l$ under the standard emissions model (Section \ref{subsubsec:standardEmissionsModel}).\\
	    $e_a (v,l)$           & $e_a: \mathbb{R}_{++} \times \mathbb{R}_+ \to {R}_{++}$, amount of \cotwo (kg) that a truck emits for traversing arc $a \in V$ with speed $v$ and payload $l$ under the improved emissions model (Section \ref{subsubsec:improvedEmissionsModel}).\\
        $\mathcal{E}(\pi,\varv,l)$  & Total amount of \cotwo emitted by a truck along a path $\pi$ under the speed policy $\varv$ and payload $l$, i.e.  $\mathcal{E}(\pi,\varv,l) = \sum_{a \in \pi} e_a(\varv(a),l)$.\\
    \bottomrule
  \end{tabularx}
\end{table*}

\subsection{Emission Models}
\label{subsec:emissionsmodels}
We discuss two emission models. The first of these models is most commonly used in recent papers on the green/pollution routing problem \cite[e.g.][]{Bektas2011, Demir2012, Franceschetti2013, Dabia2017}. We will call this the standard model. The second model is a small improvement on the standard model to disallow negative fuel consumption on downward sloping road segments.

\subsubsection{Standard Emissions Model}
\label{subsubsec:standardEmissionsModel}
The CMEM \citep{Barth2005,Scora2006,Boriboonsomsin2009} is a microscopic truck fuel consumption model that has been widely used in literature for pollution/green vehicle routing problems. The Standard model is an instantiation of the CMEM approach. Suppose a truck with the parameters given in Table \ref{tab:notation} and payload $l$ travels along arc $a \in V$ with speed $v$. In the standard emission model introduced by \cite{Bektas2011} and \cite{Demir2012}, the truck's fuel consumption is given by:

\begin{alignat}{2}
    & \Tilde{f}_a(v,l) = \frac{P \delta(a)}{v} + Q \delta(a) \left( g \sin \theta(a) + C_r g \cos \theta(a) \right) (w+l) + R \delta(a) v^2 \label{eq:cmem}\\
    \mbox{with }& P = \frac{\xi kND}{\kappa \psi}, \label{eq:P_coff}\\
    & Q = \frac{\xi}{1000 \eta \eta _{tf}\kappa \psi},\label{eq:Q_coeff}\\
    & R = \frac{ \xi C_d \rho S}{2000 \eta \eta _{tf}\kappa \psi}. \label{eq:R_coeff}
\end{alignat}
The main assumption behind Equation \eqref{eq:cmem} is that the truck parameters remain constant along each arc $a \in V$. This model sets aside a number of minor sources of fuel consumption, such as air conditioning and compressed air systems. Burning one liter of diesel in a combustion engine produces $c_e=2.67$ kg/L of \cotwo \citep{EmissionFacts2005}. Thus we find that the \cotwo emissions associated with traversing an arc $a$ with load $l$ at speed $v$ is given by,
\[
\tilde{e}_a(v,l) = c_e \Tilde{f}_a(v,l) = 2.67 \Tilde{f}_a(v,l)
\]
under the standard model.

\subsubsection{Improved Emissions Model}
\label{subsubsec:improvedEmissionsModel}
\cite{Rao2016} and \cite{Brunner2021} establish that the standard emission model gives rise to negative fuel consumption on some negative road angles that are not realistic for internal combustion engine vehicles. Thus, \cite{Brunner2021} propose the following adjustment to Equation \eqref{eq:cmem}:
\begin{equation}
    f_a(v,l) = \frac{P \delta(a)}{v} + \Big(Q \delta(a) \left(g \sin \theta(a) + C_r g \cos \theta(a) \right) (w+l) + R \delta(a) v^2 \Big)^+, \label{eq:cmem_modified}
\end{equation}
where $(\cdot)^+=\max \{\cdot,0\}$. Equation \eqref{eq:cmem_modified} shows that gravity works in favor of the vehicle over downhills arcs to compensate the energy loss caused by drag and rolling resistance force. This equation assumes that any engine-powered brakes consume a negligible amount of fuel. 
The standard and improved emission models \eqref{eq:cmem} and \eqref{eq:cmem_modified} are identical on a flat network ($\theta(a)=0$ for all $a\in A$). We note that a slight modification of the above models can also allow for electric vehicles with regenerative braking; see \cite{Larminie2012}. As before we now find that
the \cotwo emissions associated with traversing arc $a$ with load $l$ at speed $v$ is given by,
\begin{equation}
    \label{eq:FE_relation}
    e_a(v,l) = c_e f_a(v,l) = 2.67 f_a(v,l).
\end{equation}

\subsection{Optimal Speed}
\label{subsec:speedanalysis}
The most fuel efficient speed to traverse an arc depends on the emission model that is used. We will show below that there is one optimal speed for all arcs in a network under the standard emission model, but that the optimal speed may differ by arc for the improved emission model.

\subsubsection{Static Speed Optimization}

The \textit{speed optimization problem (SO)} is to compute the speed policy which minimizes the carbon emissions when a vehicle travels across an arc $a \in A$. Under the standard emissions model, SO can be formulated as,

\begin{equation*}
    \Tilde{e}_a^*= \min_{v \in [\varv^{\min}(a), \varv^{\max}(a)]} {\Tilde{e}_a(v,l)} \quad \mbox{ and } \quad \varv^s (a) = \arg\min_{v \in [\varv^{\min}(a), \varv^{\max}(a)]}{\Tilde{e}_a(v,l)}.
\end{equation*}
This implies that the most fuel efficient speed is the same along any arc $a \in A$ and is given by $\varv^s: A \to \mathbb{R}_{++}$ that is defined by, 
\begin{equation}
    \label{eq:constant_opt_speed}
    \varv^s (a) :=
    \begin{cases}
        \varv^{\min}(a)& \mbox{if }  c_{\varv} \leq \varv^{\min}(a)\\
        c_{\varv}& \mbox{if } \varv^{\min}(a) < c_{\varv} \leq \varv^{\max}(a)\\
        \varv^{\max}(a)& \mbox{if }   \varv^{\max}(a) < c_{\varv},
    \end{cases}
\end{equation}
where $c_{\varv}$,
\begin{equation}
\label{eq:freeflowstatic}
    c_{\varv}=\sqrt[3]{\frac{P}{2R}}.
\end{equation}
is the optimal speed without any speedlimits. Expressions for $P$ and $R$ are given by Equations \eqref{eq:P_coff} and \eqref{eq:R_coeff}. 
Equation \eqref{eq:freeflowstatic} is obtained by solving the first order conditions to minimize \eqref{eq:cmem} with respect to $v$. Since $c_{\varv}$ is constant along all arcs, we use the term {\em static speed} policy to denote a policy that will have a vehicle traverse every arc at the speed $\varv^s$. We note that for practically meaningful values of $\varv^{\min}(a)$ and $\varv^{\max}(a)$ the optimal speed is usually given by \eqref{eq:freeflowstatic}, or the second case in \eqref{eq:constant_opt_speed}.

\subsubsection{Dynamic Speed Optimization}
In the improved emissions model, the most fuel-efficient speed to traverse an arc depends on its slope. Under the improved emissions model \eqref{eq:cmem_modified}, the speed optimization problem is formulated as,
\begin{equation}
    \label{eq:dynamicSpeedOpt}
    e_a^*= \min_{v \in [\varv^{\min}(a), \varv^{\max}(a)]} {e_a(v,l)} \quad \mbox{ and } \quad \varv^d (a,l) = \arg\min_{v \in [\varv^{\min}(a), \varv^{\max}(a)]}{e_a(v,l)}.
\end{equation}
Note that the derivative of $e_a(v,l)=\frac{c_e P \delta(a)}{v} + c_e \left(Q \delta(a) (g \sin \theta(a) + C_r g \cos \theta(a)) (w+l) + R \delta(a) v^2 \right)^+$ with respect to $v$ is given by
\[
    \frac{\partial e_a(v,l)}{\partial v}=
    \begin{dcases}
        -\frac{c_e P \delta(a)}{v^2} & \mbox{if } 0 \leq v < \varv^t(a,l)\\
        -\frac{c_e P \delta(a)}{v^2}+ 2 c_e R \delta(a) v & \mbox{if } \varv^t(a,l) < \varv (a),
    \end{dcases}
\]
where $\varv^t: A \times \mathbb{R}_+ \to \mathbb{R}_{+}$ is defined by,
\begin{equation}
    \varv^t(a,l) :=
        \begin{dcases}
            \sqrt{\frac{-Q (g \sin \theta(a) + C_r g \cos \theta(a)) (w+l)}{R}}, & \text{if } \tan \theta(a) < -C_r\\
            0, & \text{if } \tan \theta(a) \geq -C_r.
        \end{dcases}
\end{equation}
This derivation shows that the \cotwo emissions of an arc $e_a(v,l)$ is not differentiable with respect to $v$ at the point $\varv^t(a,l)$. The speed $\varv^t$ has a physical interpretation as the terminal velocity of a vehicle on a slope with inclination $\theta$. It is the speed at which the gravitational force along the slope equals the sum of the drag and rolling resistance forces (see e.g. \cite{fox2020}).
A vehicle reaches a non-zero terminal velocity on an arc if the angle falls below $-\arctan{C_r}$. The optimal solution to the speed optimization problem in \eqref{eq:dynamicSpeedOpt} is slightly more involved as it needs to account for the terminal velocity. The solution is given in Proposition \ref{prop:opt_speed}.
 \begin{proposition}
    \label{prop:opt_speed}
    The optimal solution to the speeds optimization problem in \eqref{eq:dynamicSpeedOpt}
    is given by $\varv^d: A \times \mathbb{R}_+ \to \mathbb{R}_{++}$ that is defined by,
    \begin{equation}
        \varv^d (a,l):=
        \arg\min_{v \in [\varv^{\min}(a), \varv^{\max}(a)]}{e_a(v,l)}=
            \begin{dcases}
                \varv^{\min}(a), & \text{if } \max \{c_{\varv}, \varv^t(a,l)\} \leq \varv^{\min}(a)\\
                \max \{c_{\varv}, \varv^t(a,l)\}, & \text{if } \varv^{\min}(a) < \max \{c_{\varv}, \varv^t(a,l)\} \leq \varv^{\max}(a)\\
                \varv^{\max}(a), & \text{if } \varv^{\max}(a) < \max \{c_{\varv}, \varv^t(a,l)\}.
            \end{dcases}
    \end{equation}
\end{proposition}

\proof{Proof of Proposition \ref{prop:opt_speed}.}
We consider the case where the terminal velocity is zero, and where it is strictly positive separately.

    Case 1 ($\tan \theta(a) \geq -C_r$; $\varv^t(a,l)=0$): Equation \eqref{eq:cmem_modified} reduces to Equation \eqref{eq:cmem} so that one may verify that
    \[
        \varv^d (a,l) = c_{\varv} > \varv^t(a,l) = 0.
    \]
    
    Case 2 ($\tan \theta(a) < -C_r$; $\varv^t(a,l) >0$): In this case,
    \begin{equation*}
        e_a(v,l) = 
        \begin{dcases}
            e^1_a = \frac{c_e P \delta(a)}{v}, & \text{if } 0 \leq v < \varv^t(a,l)\\
            e^2_a = c_e \left(\frac{P \delta(a)}{v} + Q \delta(a) (g \sin \theta(a) + C_r g \cos \theta(a))(w+l) + R \delta(a) v^2 \right), & \text{if } \varv^t(a,l) \leq v.
        \end{dcases}
    \end{equation*}
    It is straightforward to verify that $e_a(v,l)$ is continuous, $e_a^1$ is convex and non-increasing in $v$, and $e_a^2$ is convex in $v$ with a minimum at $c_{\varv}$. Consequently, the optimal speed exceeds the terminal velocity, i.e. $\varv^d (a,l)\geq \varv^t(a,l)$.
    
   As $e_a(v,l)$ is convex in $v$ on $[\varv^t(a,l),\infty)$, it has an extremum at $c_{\varv}$ if $\varv^t(a,l) \leq c_{\varv}$, or at $\varv^t(a,l)$ if $\varv^t(a,l) > c_{\varv}$. It follows that the optimal speed is $\max \{c_{\varv}, \varv^t(a,l)\}$ if it lies within the allowable speed interval, $[\varv^{\min}(a),\varv^{\max}(a)]$. In case $\max \{c_{\varv}, \varv^t(a,l)\} < \varv^{\min}(a)$, then $e_a(v,l)$ is non-decreasing in $v \in [\varv^{\min}(a),\varv^{\max}(a)]$ and the optimal speed is $\varv^{\min}(a)$. If $\max \{c_{\varv}, \varv^t(a,l)\} \geq \varv^{\max}(a)$, then $e_a(v,l)$ is non-increasing in $v \in [\varv^{\min}(a),\varv^{\max}(a)]$ and the optimal speed is $\varv^{\max}(a)$.
    \Halmos
\endproof
The main insight from Proposition \ref{prop:opt_speed} is that it is efficient to use gravity to reduce the required engine power and emission.
Proposition \ref{prop:opt_speed} indicates that a static speed policy is not optimal on a path that contains downhill arcs. Thus, a {\em dynamic speed} policy ($\varv^d$), as per Proposition \ref{prop:opt_speed}, reduces a truck's fuel consumption, \cotwo emissions, and travel time since it requires higher speeds on downhills.

\subsection{The Greenest Path}
\label{subsect:gpp}
Let $n_s$ and $n_t$ be two different
vertices of $\mathcal{G}$ such that $n_t$ is reachable from $n_s$. Let $\Pi$ be the set of all possible paths between $n_s$ and $n_t$. Under a given speed policy $\varv: A \to \mathbb{R}_{++}$ and a constant payload $l$, the total \cotwo emissions of a truck to travel between $n_s$ and $n_t$ along a path $\pi \in \Pi$, $\mathcal{E}(\pi, \varv, l)$, is defined as,
\begin{equation}
    \label{eq:totalEmissions}
    \mathcal{E}(\pi, \varv, l) = \sum_{a \in \pi} e_a(\varv(a),l).
\end{equation}
Based on this definition, the \textit{greenest path problem (GPP)} is to compute the path with the least \cotwo emissions, $\pi^g$, between $n_s$ and $n_t$, i.e.
\begin{equation}
    \label{eq:greenestPath}
    \mathcal{E}^*(\varv, l) = \min_{\pi \in \Pi} \mathcal{E} (\pi, \varv, l) = \min_{\pi \in \Pi} \sum_{a \in \pi} e_a (\varv(a),l) \quad \mbox{and} \quad
    \pi^g (\varv , l) = \arg \min_{\pi \in \Pi} \mathcal{E} (\pi, \varv, l) = \arg \min_{\pi \in \Pi} \sum_{a \in \pi} e_a (\varv(a),l).
\end{equation}
We define the \textit{shortest path problem (SPP)} as the computation of the minimum-distance path ($\pi^{sp}$) between $n_s$ and $n_t$, i.e.
\begin{equation}
    \label{eq:shortestPath}
    \delta^{sp} = \min_{\pi \in \Pi} \sum_{a \in \pi} \delta(a)
    \quad \mbox{and} \quad
    \pi^{sp} = \arg \min_{\pi \in \Pi} \sum_{a \in \pi} \delta(a).
\end{equation}
The following proposition
\begin{rev}
    shows that if the elevation data is ignored and the speeds are identical along all arcs then the shortest path, $\pi^g$, is an optimal solution to GPP.
\end{rev}

\begin{rev}
    \begin{proposition}
        \label{prop:gppmodel1}
        If the road gradient $\theta(a)=0$ and the speeds $\varv(a)$ are identical for all arcs $a \in A$, then the Greenest Path ($\pi^g (\varv , l)$) is the Shortest Path ($\pi^{sp}$) .
    \end{proposition}
\end{rev}

\proof{Proof of Proposition \ref{prop:gppmodel1}}
    Let angle $\theta(a)=0$ for all $a \in A$, and the payload $l$ and speed policy $\varv(a)$ be identical, i.e. $\varv(a)= \varv^*$, where $\varv^* \in \mathbb{R}$ is constant. Taking into account that $\sin \theta(a) = 0$ and $\cos \theta(a) =1$ for all $a \in A$, the GPP implies that,
    \begin{alignat*}{2}
         \mathcal{E}^* & = c_e \min_{\pi \in \Pi} \sum_{a \in \pi}  \frac{P \delta(a)}{\varv(a)} + Q \delta(a) (g \sin \theta(a) + C_r g \cos \theta(a)) (w+l) + R \delta(a) \varv(a)^2\\
        & = c_e \min_{\pi \in \Pi} \sum_{a \in \pi}  \frac{P \delta(a)}{\varv^*} + Q \delta(a) C_r g (w+l) + R \delta(a) {\varv^*}^2\\
        & = c_e \left(\frac{P}{\varv^*} + Q C_r g (w+l) + R {\varv^*}^2 \right) \min_{\pi \in \Pi} \sum_{a \in \pi} \delta(a)\\
        & = c_e \left(\frac{P}{\varv^*} + Q C_r g (w+l) + R {\varv^*}^2\right) \delta_{n_s,n_t}^*.
    \end{alignat*}
    Thus, the $\pi^{sp}$ satisfies this problem that proves the proposition.
    \Halmos
\endproof
\begin{rev}
    When the speeds are bounded by traffic or variable speed limits, then the analogous result holds for the fastest path.
\end{rev}
It is straightforward to verify that the greenest path is the fastest path when all road gradients are zero \begin{rev} and the speeds are constant \end{rev}; see Proposition \ref{prop:gppmodel1}.
Further notice that by Proposition \ref{prop:opt_speed}, the speed $c_{\varv}$ in \eqref{eq:constant_opt_speed} is optimal for all arcs when $\theta(a)=0$ for all $a\in A$. This implies that a decision maker will believe the shortest path is the greenest path when she ignores elevation data. 

Nonetheless, the improved emissions model and Proposition \ref{prop:opt_speed} show that if the elevation data is considered, the speed along each segment of a path can change. Even under the static speed policy, the greenest path is not necessarily the shortest due to the non-linearity of emission along an arc when in the gradient. 
We note that the emission model does not explicitly account for acceleration and deceleration of a vehicle and so the estimates emissions $\mathcal{E}^*(\varv^d,l)$ are a lower-bound for the \cotwo emissions of a truck traveling from $n_s$ to $n_t$.

\subsection{The Asymptotic Greenest Path}
\label{subsect:asymcal}
In this section, we explore the greenest path as the payload becomes arbitrarily large. Let $e_a^\prime (v)$ be the \cotwo emissions per unit payload when a truck traverses arc $a \in A$ with speed $v$, that is to say,
\begin{equation}
    \label{eq:emissionsperload}
    e_a^\prime (v) = \frac{e_a (v,l)}{l} =
    \frac{c_e P \delta(a)}{v l} + c_e\left(Q \delta(a) (g \sin \theta(a) + C_r g \cos \theta(a)) \left(1+\frac{w}{l}\right) + \frac{R}{l} \delta(a) v^2\right)^+.
\end{equation}
Consider two distinct connected vertices $n_s$ and $n_t$.
Observe that under a speed policy $\varv:A\to\R^+$ and a constant load $l \in \mathbb{R}_+$, the greenest path, i.e. $\pi^g (\varv,l)$, minimizes the total \cotwo emissions and the total \cotwo emissions per unit payload between $n_s$ and $n_t$. Thus, we can interchangeably use the total \cotwo emissions and the total \cotwo emissions per unit payload to compute the greenest path.

Let $\Pi$ be the set of all paths from $n_s$ to $n_t$. Let $\Pi_d \subseteq \Pi$ be the subset of paths $\Pi$ that are entirely downhill with a slope below $\arctan(-C_r)$, i.e. $\tan \theta(a) < -C_r$ for all $a\in \pi$ with $\pi\in\Pi_d$. We can now state the definition of the asymptotic greenest path:

\begin{definition}
The asymptotic greenest path satisfies
\begin{rev}
\begin{equation}
        \label{eq:piInftydef}
        \pi^{\infty} (\varv) \in
        \begin{dcases}
            \arg \min_{\pi \in \Pi} \lim_{l \to \infty} \sum_{a \in \pi} e_a(\varv(a)) & \mbox{if } \Pi_d \not = \emptyset\\
            \arg \min_{\pi \in \Pi} \lim_{l \to \infty} \sum_{a \in \pi} e_a^\prime(\varv(a),l) & \mbox{if } \Pi_d = \emptyset.
        \end{dcases}
\end{equation}
\end{rev}
\end{definition}
Note that the set $\Pi_d$ plays an important role in this definition. The emission per load vanishes for any sufficiently steep down downhill path ($\Pi_d\neq\emptyset$) because gravity will get the vehicle to its destination. Among all sufficiently steep downhill paths ($\pi\in\Pi_d$), the one with the least absolute emission is given by the second case in \eqref{eq:piInftydef}. When gravity does not suffice to move a vehicle from its origin to its destination ($\Pi_d=\emptyset$) then the asymptotic greenest path is the one that minimizes emissions per load; see case 1 in \eqref{eq:piInftydef}.
The following proposition demonstrates that $\pi^{\infty}(\varv)$ exists and provides an explicit form to compute it.
\begin{proposition}
    \label{prop:AP}
     $\pi^{\infty} (\varv)$ exists and can be computed as follows.
    \begin{equation}
    \label{eq:piInf}
    \pi^{\infty} (\varv)\in
    \begin{cases}
        \arg \min_{\pi \in \Pi_d} \sum_{a \in \pi} t_a (\varv (a)) & \text{if } \Pi_d \not = \emptyset\\
        \arg \min_{\pi \in \Pi} \sum_{a \in \pi} h^{\prime} (a) & \text{if } \Pi_d = \emptyset,
    \end{cases}
\end{equation}
    where $t_a:\mathbb{R}_{++} \to \mathbb{R}_{++}$, is defined by
    \[t_a (v) = \frac{\delta(a)}{v},\]
    and $h^{\prime}: A \to \mathbb{R}_+$, is defined by  
    \[h^{\prime} (a) = \delta(a) \sin \left(\theta(a) + \arctan C_r \right)^+,\]
    if $-90^\circ < \theta(a)+\arctan C_r < 90^\circ$ for all $a \in A$.
\end{proposition}
We call $\pi^{\infty} (\varv)$ the asymptotic greenest path. Proposition \ref{prop:AP} explains that $\pi^{\infty} (\varv)$ is the fastest downward path $\pi \in \Pi_d$,  if $\Pi_d$ is non-empty. Otherwise, it is the path with the minimum total augmented ascents, $h^{\prime}$. Evidently, $\pi^{\infty} (\varv)$ can be computed using the algorithms offered to solve the shortest path problem (e.g. \cite{Dijkstra1959}). The requirement that $-90^\circ < \theta(a)+\arctan C_r < 90^\circ$ for all $a \in A$ is completely benign.
\proof{Proof of Proposition \ref{prop:AP}}
    For all payloads $l \in \mathbb{R}_+$, and any speed policy $\varv$, $\pi^g(\varv, l)$ exists from $n_s$ to $n_t$, since by Equations \eqref{eq:FE_relation} and \eqref{eq:emissionsperload} there are no negative emissions cycles between the vertices. Suppose that the payload $l$ satisfies,
    \begin{equation}
        \label{eq:LargeLCond}
        l \geq \max_{a \in A} \left\{  \frac{R(\varv^{\max}(a))^2}{-Q(g \sin \theta(a) + C_r g \cos \theta(a))}-w \right\}.
    \end{equation}
    Then for arc $a \in A$,
    \begin{subequations}
    \begin{alignat}{2}
        \label{eq:e_aLargeL}
        & e_a(\varv(a),l) =
        \begin{dcases}
            \frac{c_e P \delta(a)}{\varv(a)} & \text{if } \tan \theta(a) < -C_r\\
            \frac{c_e P \delta(a)}{\varv(a)} + c_e \left(Q \delta(a) (g \sin \theta(a) + C_r g \cos \theta(a)) (w + l) + R \delta(a) \varv(a)^2 \right) & \text{if } \tan \theta(a) \geq -C_r,
        \end{dcases}\\
        \label{eq:el_aLargeL}
        & e_a^\prime(\varv(a)) =
        \begin{dcases}
            \frac{c_e P \delta(a)}{l\varv(a)} & \text{if } \tan \theta(a) < -C_r\\
            \frac{c_e P \delta(a)}{l\varv(a)} + c_e \left(Q \delta(a) (g \sin \theta(a) + C_r g \cos \theta(a)) (1 + \frac{w}{l}) + \frac{R \delta(a) \varv(a)^2}{l}  \right) & \text{if } \tan \theta(a) \geq -C_r,
        \end{dcases}
    \end{alignat}
    \end{subequations}
    by Equations \eqref{eq:FE_relation} and \eqref{eq:emissionsperload}.
    
    Suppose that $\Pi_d$ is a non-empty set. For this case, we use the total \cotwo emissions to compute the $\pi^{\infty}(\varv)$. Thus, by Equations \eqref{eq:e_aLargeL},
    \begin{equation*}
        \label{eq:inftye}
        \lim_{l \to \infty} \sum_{a \in \pi} e_a(\varv(a),l)=
        \begin{cases}
            \sum_{a \in \pi} \frac{c_e P}{\varv(a)} \delta(a) & \text{if } \pi \in \Pi_d\\
            \infty & \text{if } \pi \in \Pi \setminus \Pi_d.
        \end{cases}
    \end{equation*}
    Then it follows that from Equation \eqref{eq:piInftydef} that
    \begin{equation*}
        \pi^{\infty} (\varv) = \arg \min_{\pi \in \Pi_d} \sum_{a \in \pi} \frac{c_e P \delta(a)}{\varv(a)} = \arg \min_{\pi \in \Pi_d} \sum_{a \in \pi} \frac{\delta(a)}{\varv(a)}=\arg \min_{\pi \in \Pi_d} \sum_{a \in \pi} t_a (\varv (a)),
    \end{equation*}
    since $P$ and $c_e$ are constant across all arcs $a \in A$.
    
    Now, suppose that $\Pi_d$ is an empty set. For this case, we use the total \cotwo emissions per unit load to compute $\pi^{\infty}(\varv)$.
    Thus, by Equation \eqref{eq:el_aLargeL},
    \begin{alignat*}{2}
        \label{eq:inftyel}
        & \lim_{l \to \infty} \sum_{a \in \pi} e_a^\prime(\varv(a))= \sum_{a \in \pi} c_e Q (g \sin \theta(a) + C_r g \cos \theta(a))^+ \delta(a) \\
        = \ \ & c_e Q g \sqrt{1+{C_r}^2} \sum_{a \in \pi} \delta(a) \left(\sin \left(\theta(a) + \arctan C_r \right)\right)^+ = c_e Q g \sqrt{1+{C_r}^2} \sum_{a \in \pi} \delta(a) \sin \left(\theta(a) + \arctan C_r \right)^+,
    \end{alignat*}
    as  $-90^\circ < \theta(a)+\arctan C_r < 90^\circ$ for all $a \in A$ by supposition. Again, since $c_e$, $Q$, and $C_r$ are constant for all $a \in A$, by Equation \eqref{eq:piInftydef},
    \begin{equation*}
        \pi^{\infty} (\varv) = \arg \min_{\pi \in \Pi} \sum_{a \in \pi} \delta(a) \sin \left(\theta(a) + \arctan C_r \right)^+ = \arg \min_{\pi \in \Pi} \sum_{a \in \pi} h^{\prime} (a).
        \Halmos
    \end{equation*}
\endproof
Proposition \ref{prop:AP} demonstrates the convergence of the $\pi^g(\varv , l)$ to the $\pi^{\infty}(\varv)$ for a very large payload. On the other hand, by Proposition \ref{prop:gppmodel1} the shortest path, $\pi^{sp}$, is the greenest path under the dynamic speed policy, i.e. $\pi^g(\varv^d,l)$, if $w+l=0$ and $\varv^{\min}(a) \leq c_{\varv} \leq\varv^{\max}(a)$ for all $a \in A$. The reason is that if $w+l=0$, the dynamic speed policy equals the static speed policy ($\varv^d=\varv^s$) by Proposition \ref{prop:opt_speed}. Therefore, one can argue that $\pi^g(\varv , l)$ diverges from $\pi^{sp}$ and converges to the $\pi^{\infty}(\varv)$ as the load increases. We will explore this idea in Section \ref{subsec:AP_num} through numerical experiments.

Finally, if the payload $l$ satisfies Inequality \eqref{eq:LargeLCond}, by Proposition \ref{prop:opt_speed}, $\varv^d (a,l)$, for arc $a \in A$ can be computed as follows.
\begin{equation*}
    \varv^d(a) =
    \begin{cases}
        \varv^{\min}(a), & \text{if } \tan \theta(a) \geq -C_r \wedge c_{\varv} \leq \varv^{\min}(a)\\
        c_{\varv}, & \text{if } \tan \theta(a) \geq -C_r \wedge \varv^{\min}(a) < c_{\varv} \leq \varv^{\max}(a)\\
        \varv^{\max}(a), & \text{if } \tan \theta(a) < -C_r \vee \tan \theta(a) \geq -C_r \wedge \varv^{\max}(a) < c_{\varv}.\\
    \end{cases}
\end{equation*}
Consequently, if $\varv^{\min}(a)$ and $\varv^{\max}(a)$ are constant for all arcs $a \in A$ and if $\Pi_d$ is non-empty, then $\pi^{\infty}(\varv^d) = \pi^{\infty}(\varv^s)$, by Proposition \ref{prop:AP}. Evidently, if $\Pi_d$ is empty then Proposition \ref{prop:AP} requires $\pi^{\infty}(\varv,l)$ to be independent of the speed policy $\varv$. As a result, $\pi^{\infty}(\varv^d) = \pi^{\infty}(\varv^s)$ if $\varv^{\min}(a)$ and $\varv^{\max}(a)$ are constant for all arcs $a \in A$.

\section{Numerical Experiments}
\label{section:NumericalExperiments}
In this section, we explore the value of using elevation data to inform routing and speed decisions to reduce emissions over a comprehensive data set. 
Additionally, we explore the major drivers of \cotwo emissions reduction. We benchmark the greenest path ($\pi^g$) and dynamic speed policy ($\varv^d$) against the shortest path ($\pi^{sp}$) and the static speed policy ($\varv^s$). Note that the shortest path is also the greenest path under a dynamic speed policy (i.e. $\pi^g(\varv^d,l) = \pi^{sp}$) if the effect of road gradients is ignored, as shown in Proposition \ref{prop:gppmodel1}. We also study how the greenest path changes, $\pi^g (\varv , l)$, as the payload $l$ increases and how the asymptotic greenest path $\pi^{\infty}(\varv)$ performs in terms of \cotwo emissions reduction and similarity to $\pi^g (\varv , l)$. 
In our numerical experiments the asymptotic greenest path under the dynamic speed policy, i.e. $\pi^{\infty}(\varv^d)$, and the one under the static speed policy, i.e. $\pi^{\infty}(\varv^s)$, are identical since $\varv^{\min}(a)$ and $\varv^{\max}(a)$ are constant for all $a \in A$ (see Section \ref{subsect:asymcal}), i.e. $\pi^{\infty} =\pi^{\infty}(\varv^d)=\pi^{\infty}(\varv^s)$. 

Given a pair of source and target vertices and a constant payload $l$, we compute the relative \cotwo emissions reduction of one policy in comparison with another. In particular we study the \cotwo reduction of using path-speed policy 2, $d_2= (\pi_2, \varv_2, l)$, relative to path-speed policy 1, $d_1= (\pi_1, \varv_1, l)$, ($\%\mathcal{E}_{d_1}^{d_2}$) to quantify the benefit of using the elevation data in \cotwo reduction. That is to say,
\[
    \%\mathcal{E}_{d_1}^{d_2}=100 \cdot \frac{\mathcal{E} (\pi_1, \varv_1, l) - \mathcal{E} (\pi_2, \varv_2, l)}{\mathcal{E} (\pi_1, \varv_1, l)},
\]
where $\mathcal{E} (\pi_i, \varv_i, l), i=1,2$ is the total \cotwo emissions as per Equation \eqref{eq:totalEmissions}.
If $\pi_i$, $i=1,2$, is a greenest path then $\pi_i = \pi^g(\varv_i,l)$.
Similarly, we compute the relative distinction between the paths of policies $\pi_1$ and $\pi_2$ ($\%\delta_{\pi_1}^{\pi_2}$) weighted by distance, as follows. 
\[
    \%\delta_{\pi_1}^{\pi_2} =100 \cdot \sum_{a \in \pi_1 \setminus \pi_2} \delta(a) \mathbin{/} \sum_{a \in \pi_1} \delta(a).
\]
Table \ref{tab:compstud} briefly summarizes the ratios that we use in our comparative studies.

\begin{table*}[!htbp]
  \caption{\textsf{List of ratios used in the comparative studies.}}
	\fontsize{8pt}{9pt}\selectfont
\scriptsize
  \label{tab:compstud}
  \begin{tabularx}{\textwidth}{l X}
    \toprule
    Ratio & Description\\
		\midrule
    $\%\mathcal{E}_{(\pi^{sp}, \varv^s, l)}^{(\pi^g, \varv^d, l)}$ \vspace{.2cm}  &
    Relative \cotwo emissions reduction by selecting the greenest path with the dynamic speed policy relative to the shortest path with the static speed policy given the load $l$.\\
    $\%\mathcal{E}_{(\pi^{sp}, \varv^s, l)}^{(\pi^g, \varv^s, l)}$ \vspace{.2cm} &
    Relative \cotwo emissions reduction by selecting the greenest path with the static speed policy relative to the shortest path with the static speed policy given the load $l$.\\
    $\%\mathcal{E}_{(\pi^{sp}, \varv^d, l)}^{(\pi^g, \varv^d, l)}$ \vspace{.2cm} &
    Relative \cotwo emissions reduction by selecting the greenest path with the dynamic speed policy relative to the shortest path with the dynamic speed policy given the load $l$.\\
    $\%\mathcal{E}_{(\pi^g, \varv^s, l)}^{(\pi^g, \varv^d, l)}$ \vspace{.2cm} &
    Relative \cotwo emissions reduction by selecting the greenest path with the dynamic speed policy relative to the greenest path with the static speed policy given the load $l$.\\
    $\%\mathcal{E}_{(\pi^{sp}, \varv^d, l)}^{(\pi^{\infty}, \varv^d, l)}$ \vspace{.2cm} &
    Relative \cotwo emissions reduction by selecting the asymptotic greenest path with the dynamic speed policy relative to the shortest path with the dynamic speed policy given the load $l$.\\
    $\%\mathcal{E}_{(\pi^{sp}, \varv^s, l)}^{(\pi^{\infty}, \varv^s, l)}$ \vspace{.2cm} &
    Relative \cotwo emissions reduction by selecting the asymptotic greenest path with the static speed policy relative to the shortest path with the static speed policy given the load $l$.\\
    $\%\mathcal{E}_{(\pi^g, \varv^d, l)}^{(\pi^{\infty}, \varv^d, l)}$ \vspace{.2cm} &
    Relative \cotwo emissions reduction by selecting the asymptotic greenest path with the dynamic speed policy relative to the greenest path with the dynamic speed policy given the load $l$.\\
    $\%\mathcal{E}_{(\pi^g, \varv^s, l)}^{(\pi^{\infty}, \varv^s, l)}$ \vspace{.2cm} &
    Relative \cotwo emissions reduction by selecting the asymptotic greenest path with the static speed policy relative to the greenest path with the static speed policy given the load $l$.\\
    $\%\delta_{\pi^{sp}}^{\pi^g (\varv^d, l)}$ \vspace{.2cm} &
    Ratio of the length of the shortest path that is not shared with the greenest path under the dynamic speed policy given the load $l$.\\
    $\%\delta_{\pi^{sp}}^{\pi^g (\varv^s, l)}$ \vspace{.2cm} &
    Ratio of the length of the shortest path that is not shared with the greenest path under the static speed policy given the load $l$.\\
    $\%\delta_{\pi^g (\varv^d, l)}^{\pi^g (\varv^s, l)}$ \vspace{.2cm} &
    Ratio of the length of the greenest path under the dynamic speed policy that is not shared with the greenest path under the static speed policy given the load $l$.\\
    $\%\delta_{\pi^g (\varv^d, l)}^{\pi^{\infty}}$ \vspace{.2cm} &
    Ratio of the length of the greenest path under the dynamic speed policy that is not shared with the asymptotic greenest path given the load $l$.\\
    $\%\delta_{\pi^g (\varv^s, l)}^{\pi^{\infty}}$ \vspace{.2cm} &
    Ratio of the length of the greenest path under the static speed policy that is not shared with the asymptotic greenest path given the load $l$.\\
    \bottomrule
  \end{tabularx}
\end{table*}

In Section \ref{subsec:testbed}, we outline the test-bed that we consider. This test-bed comprises 25 cities and all the ratios in Table\ref{tab:compstud} are computed for instances in this test-bed. 
We present the results of our computations in Sections \ref{subsec:emissionreduction} through \ref{subsec:drivers}. Section \ref{subsec:emissionreduction} focuses on the results for the \cotwo emissions reduced by $\pi^g$ and $\varv^d$ relative to $\pi^{sp}$ and $\varv^s$. In Section \ref{subsec:trajectories}, we address the distinctions between $\pi^g$ and $\pi^{sp}$ and the effect of the speed policies $\varv^s$ and $\varv^d$ on the greenest path. In Section \ref{subsec:AP_num}, we study the asymptotic greenest path $\pi^{\infty}$ and explore the performance of  $\pi^{\infty}$ relative to the shortest path $\pi^{sp}$ and the greenest path $\pi^g$ in terms of \cotwo emissions reduction. In Sections \ref{subsec:emissionreduction} through \ref{subsec:AP_num}, we elaborate on how payload affects our results.
Finally, Section \ref{subsec:drivers} concentrates on the major determinants of \cotwo emissions reduction and path alteration.

\subsection{Data and test-bed}
\label{subsec:testbed}
We consider the 25 cities shown in Table \ref{tab:cities} and three truck types, namely heavy-duty diesel (HDD), medium-duty diesel (MDD), and light-duty diesel (LDD) for which we utilise the typical parameters as found in Table \ref{tab:trucks} of \cite{koc2014}.

\begin{table*}[!htbp]
  \fontsize{8pt}{9pt}\selectfont
  \scriptsize
  \parbox[t]{.45\linewidth}{
  \centering
  \caption{\textsf{List of cities and sample sizes.}}
  \label{tab:cities}
  \begin{tabular}{l l r}
    \toprule
        City & Country & Number of S-T pairs \\ \midrule
        Amsterdam & Netherlands & 146842 \\ 
        Ankara & Turkey & 114675 \\ 
        Athens & Greece & 114687 \\ 
        Barcelona & Spain & 100649 \\ 
        Canberra & Australia & 131565 \\ 
        Geneva & Switzerland & 125605 \\ 
        Guadalajara & Mexico & 126600 \\ 
        Guangzhou & China & 124940 \\ 
        Huston & US & 148548 \\ 
        Istanbul & Turkey & 94999 \\ 
        Johannesburg & South Africa & 130054 \\ 
        Lima & Peru & 119174 \\ 
        Los Angeles & US & 142452 \\ 
        Luxembourg & Luxembourg & 136201 \\ 
        Madrid & Spain & 117686 \\ 
        Mexico City & Mexico & 93076 \\ 
        Monterrey & Mexico & 128619 \\ 
        Mumbai & India & 143116 \\ 
        New York & US & 146875 \\ 
        Rome & Italy & 124452 \\ 
        San Francisco & US & 93504 \\ 
        Santiago & Chile & 142212 \\ 
        Shiraz & Iran & 90935 \\ 
        Tehran & Iran & 110104 \\ 
        Tel Aviv & Israel & 131066 \\ \bottomrule
    \end{tabular}
  }
  \hfill
  \parbox[t]{.45\linewidth}{
  \centering
  \caption{\textsf{Truck parameters.}}
  \label{tab:trucks}
    \begin{tabular}{l r r r}
        \toprule 
        Parameter & HDD & MDD & LDD \\ \midrule
        $w$ & 14000 & 5500 & 3500 \\ 
        $L$ & 26000 & 12500 & 4000 \\ 
        $k$ & 0.15 & 0.2 & 0.25 \\ 
        $N$ & 30 & 36.67 & 38.34 \\ 
        $D$ & 10.5 & 6.9 & 4.5 \\ 
        $\eta$ & 0.45 & 0.45 & 0.45 \\ 
        $\eta_{tf}$ & 0.45 & 0.45 & 0.45 \\ 
        $\xi$ & 1 & 1 & 1 \\ 
        $\kappa$ & 44 & 44 & 44 \\ 
        $\phi$ & 737 & 737 & 737 \\ 
        $C_d$ & 0.9 & 0.7 & 0.6 \\ 
        $\rho$ & 1.2041 & 1.2041 & 1.2041 \\ 
        $A$ & 10 & 8 & 7 \\ 
        $g$ & 9.81 & 9.81 & 9.81 \\ 
        $C_r$ & 0.01 & 0.01 & 0.01 \\ \bottomrule
    \end{tabular}
  }
\end{table*}

We use OpenStreetMap's database \citep{OpenStreetMap} to obtain the information of a 2D road network including all vertices within a 20 km radius around a manually selected point for each city. We only use roads that the database designates as public and driveable \citep{OMNXDoc}. We only consider arcs with a gradient ranging from $-10\%$ to $10\%$ (i.e. $[-5.71^{\circ},5.71^{\circ}]$)
so that gradients are in line with the implicit assumptions of the modified emissions model \eqref{eq:cmem_modified}.
We retrieve the elevation (height above sea level) of the vertices from the \cite{USGS}'s SRTM 1 Arc-Second Global data sets. 
We consider payloads of $30\%$, $40\%$, $50\%$, $60\%$, $70\%$, and $80\%$ of the maximum capacity for each truck type.
For all arcs the $\varv^{\max}=90$ km/h and $\varv^{\min}=20$ km/h.

We select several unique pairs of source and target vertices uniformly at random for each city. We make sure that the vertices in each pair are non-identical and connected. The number of selected pairs of vertices (sample size) for each city is presented in Table \ref{tab:cities}.

The Dijkstra algorithm \citep{Dijkstra1959} is used to solve the shortest path and the greenest path problems. We use the arcs' distance $\delta (a)$, $a \in A$, to compute the shortest path $\pi^{sp}$.
We consider two speed policies, namely dynamic speed policy, $\varv^d$, and static speed policy, $\varv^s$ to calculate the the \cotwo emissions, $e_a(\varv,l)$, for all arcs. Then we use the calculated $e_a(\varv,l)$ to compute the greenest paths ($\pi^g(\varv^d,l)$ and $\pi^g(\varv^s,l)$).
We use the Dijkstra algorithm to compute $\pi^{\infty}$ as per Proposition \ref{prop:AP}.

We consider two sets of ratios, as shown in 
Table \ref{tab:compstud}, to compare the different path ($\pi^{sp}$, $\pi^g$, and $\pi^{\infty}$) and speed ($\varv^s$ and $\varv^d$) policies. The first group of ratios measure the relative \cotwo emissions reduction and the second group measures the geometrical distinctions between the paths.  
We compute the ratios for a full factorial combination of trucks and payloads traversing all samples, resulting in a total of more than $58.5$ million path selection instances with a total shortest distance of more than $1.27$ billion km. Evidently, it is hardly possible to determine \cotwo emissions experimentally by letting trucks drive $1.27$ billion km as the approaches of \cite{Boriboonsomsin2012} and \cite{Scora2015}. The confidence intervals of any estimate reported later are negligibly small due to the large sample size.
Considering the large test-bed, we notice that the distribution and the sample mean of the ratios varies between different cities. For a given city, we use the overbar to denote the average of a ratio across all instances within a city. For instance, $\overline{\%\mathcal{E}}_{(\pi^{sp}, \varv^s, l)}^{(\pi^g, \varv^d, l)}$ for a city represents the sample mean of $\%\mathcal{E}_{(\pi^{sp}, \varv^s, l)}^{(\pi^g, \varv^d, l)}$ for that city.

\subsection{Results: \cotwo Emissions Reduction by Greenest Path and Dynamic Speed Policy}
\label{subsec:emissionreduction}
In this section, we consider the payload as a percentage of the truck's maximum carrying capacity rather than the payload in kilograms, to make the notations simpler. For instance, $l=60\%$ indicates that the payload equals 60\% of the maximum capacity of the truck. Since the payload varies the results, we use $l=60\%$ as our base case to maintain consistency.

\begin{figure}[t!]
	\centering
  \includegraphics[width=1\textwidth]{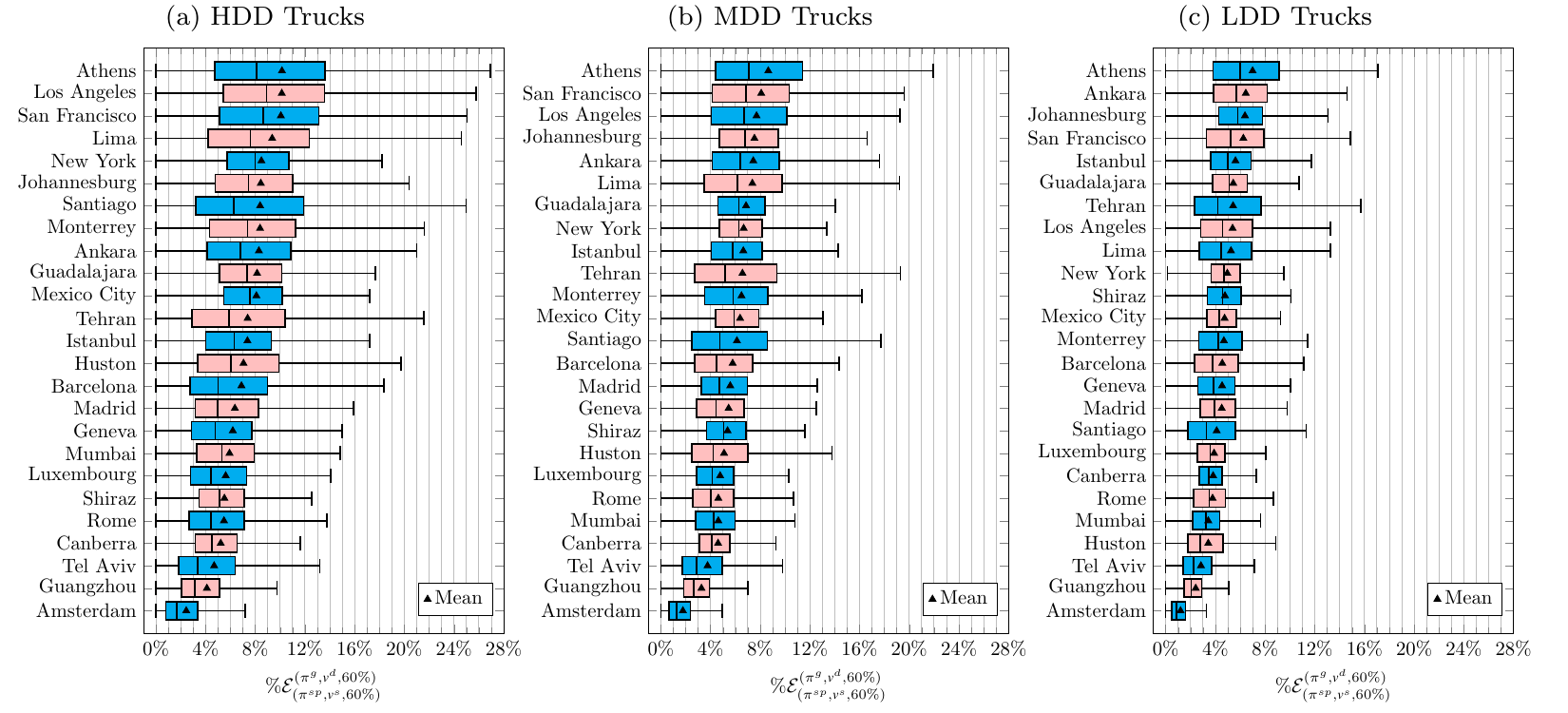}\\
	\caption{\textsf{Relative \cotwo emissions reduction by selecting $(\pi^g, \varv^d, 60\%)$ rather than $(\pi^{sp}, \varv^s, 60\%)$.}}
  \label{fig:baseSaving_new_old}
\end{figure}

\begin{figure}[t!]
	\centering
  \includegraphics[width=1\textwidth]{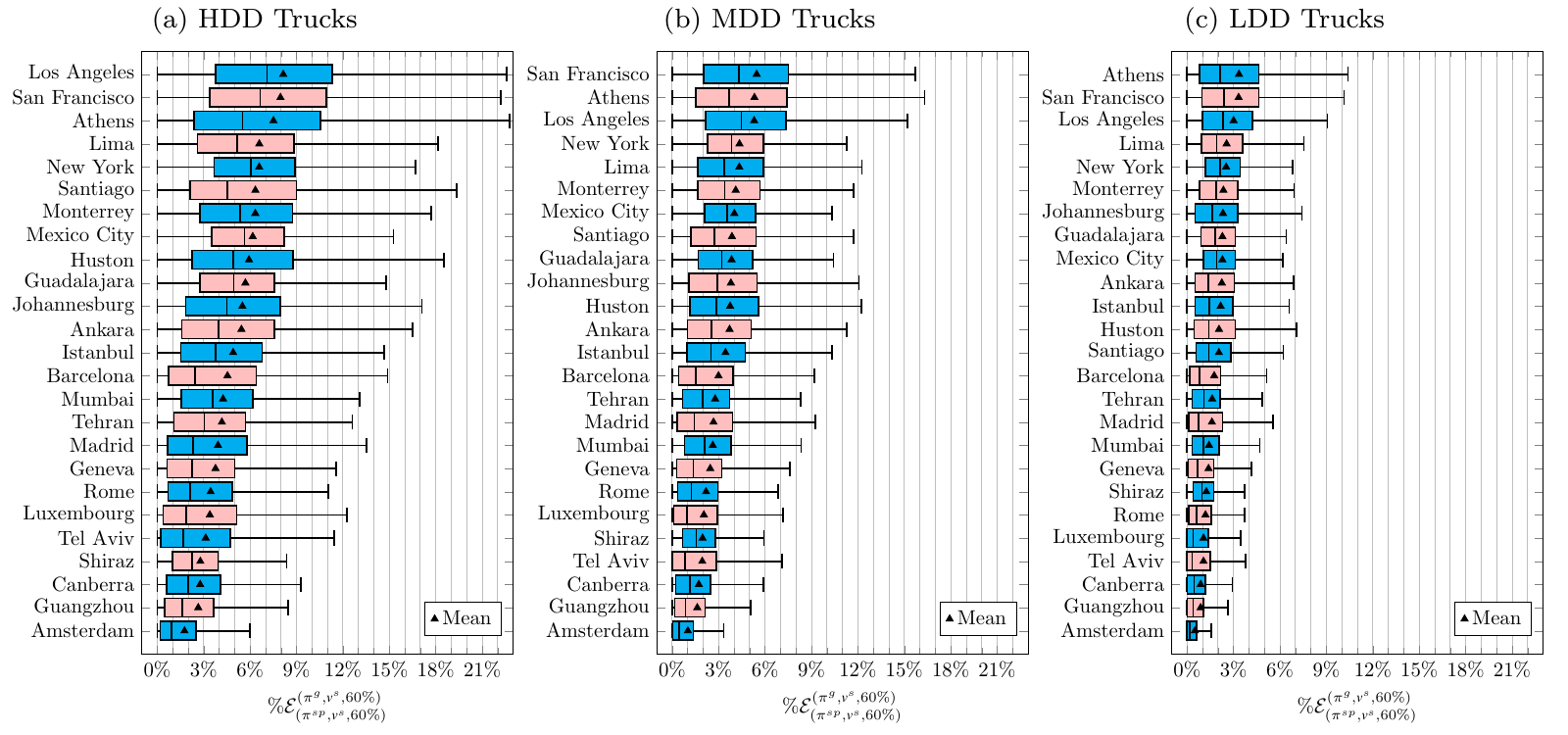}\\
	\caption{\textsf{Relative \cotwo emissions reduction by selecting $(\pi^g, \varv^s, 60\%)$ rather than $(\pi^{sp}, \varv^s, 60\%)$.}}
  \label{fig:baseSaving_old_old}
\end{figure}

\begin{figure}[b!]
	\centering
  \includegraphics[width=1\textwidth]{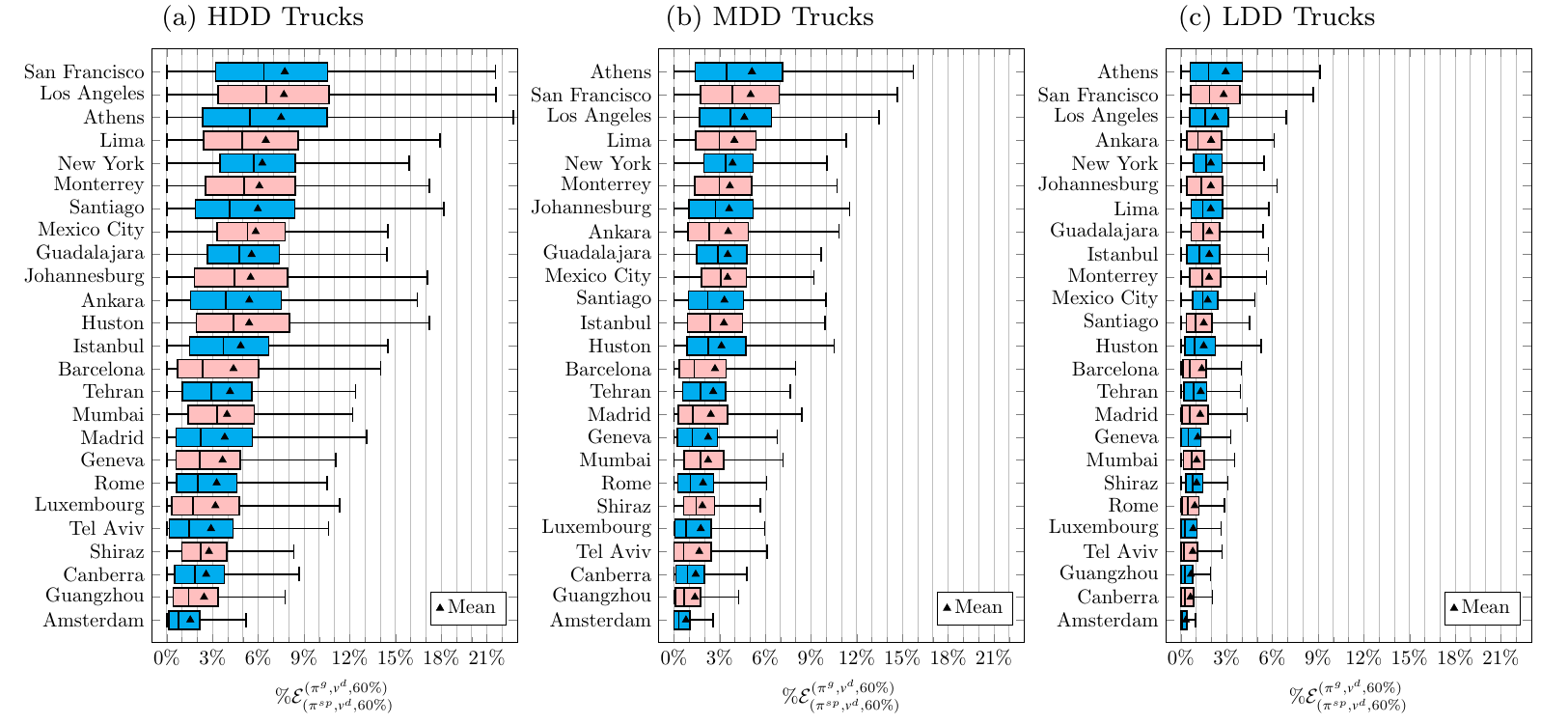}\\
	\caption{\textsf{Relative \cotwo emissions reduction by selecting $(\pi^g, \varv^d, 60\%)$ rather than $(\pi^{sp}, \varv^d, 60\%)$.}}
  \label{fig:baseSaving_new_new}
\end{figure}

\begin{figure}[t!]
	\centering
  \includegraphics[width=1\textwidth]{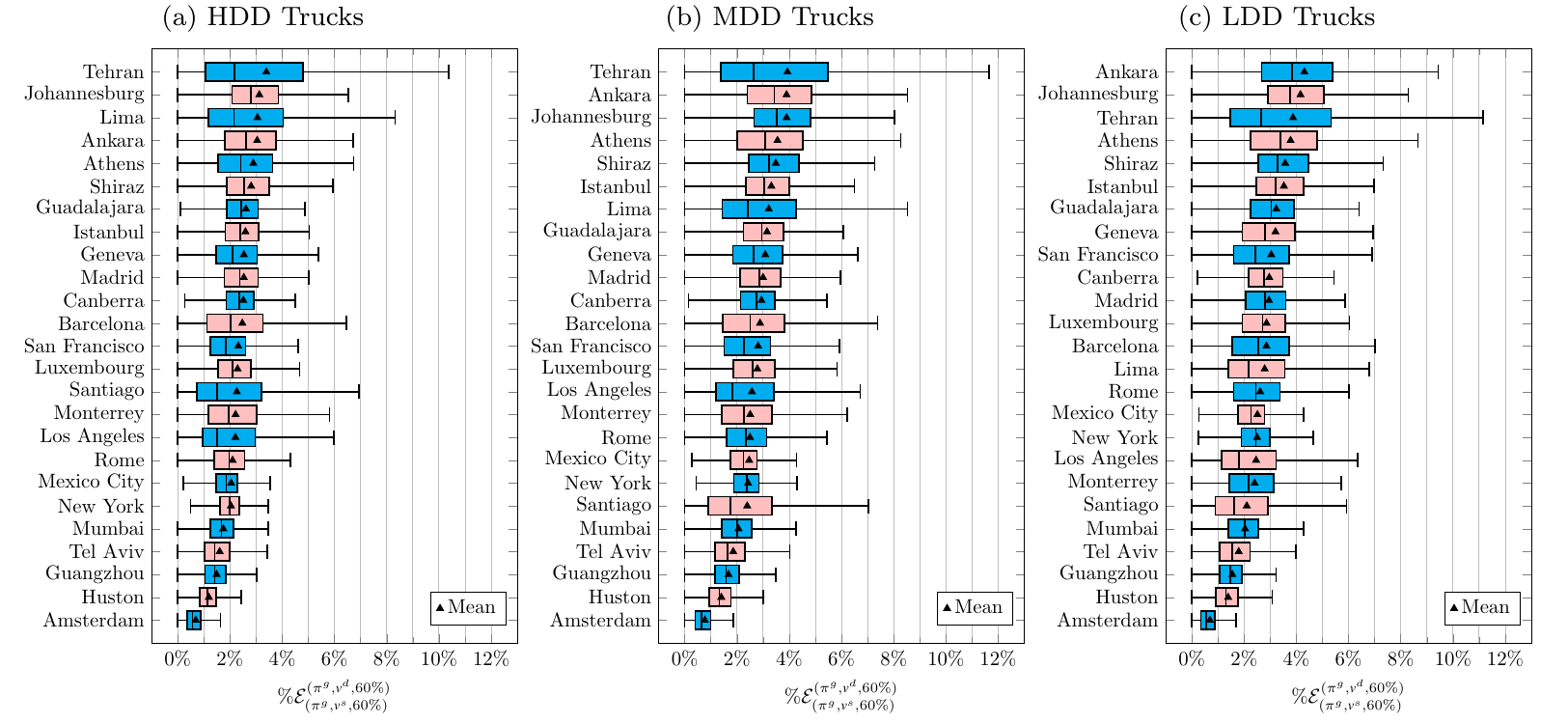}\\
	\caption{\textsf{Relative \cotwo emissions reduction by selecting $(\pi^g, \varv^d, 60\%)$ rather than $(\pi^g, \varv^s, 60\%)$.}}
  \label{fig:baseSaving_GPP}
\end{figure}

Figures~\ref{fig:baseSaving_new_old} through \ref{fig:baseSaving_GPP} visualize the empirical distribution of \cotwo emissions reduction ratios for the base case instances.
We present the distributions separately for each truck type and each city. The sample size of each empirical distribution is listed in Table \ref{tab:cities}.

Figure~\ref{fig:baseSaving_new_old} shows that $\overline{\%\mathcal{E}}_{(\pi^{sp}, \varv^s, 60\%)}^{(\pi^g, \varv^d, 60\%)}$ lies between $4.11\%$ and $10.15\%$ for HDD trucks across all cities except Amsterdam. 
Figure~\ref{fig:baseSaving_new_old} also shows that $\%\mathcal{E}_{(\pi^{sp}, \varv^s, 60\%)}^{(\pi^g, \varv^d, 60\%)}$ decreases in truck class such that $\overline{\%\mathcal{E}}_{(\pi^{sp}, \varv^s, 60\%)}^{(\pi^g, \varv^d, 60\%)}$ ranges from $3.27\%$ to $8.65\%$ for MDD, and from $2.41\%$ to $7.00\%$ for LDD trucks, in the same cities. 
Amsterdam, a known flat city, is the lone exception, but even here $\overline{\%\mathcal{E}}_{(\pi^{sp}, \varv^s, 60\%)}^{(\pi^g, \varv^d, 60\%)}$ are $2.44\%$, $1.78\%$, and $1.19\%$, respectively, showing that it is possible to use significantly more fuel-efficient paths.
The distribution of $\%\mathcal{E}_{(\pi^{sp}, \varv^s, 60\%)}^{(\pi^g, \varv^d, 60\%)}$, on the other hand, sheds more light on the potential \cotwo emissions reduction by using the greenest path with a dynamic speed policy, $\pi^g (\varv^d,60\%)$. 
In Los Angeles, for instance, $25\%$ of cases have a $\%\mathcal{E}_{(\pi^{sp}, \varv^s, 60\%)}^{(\pi^g, \varv^d, 60\%)}$ of at least $13.57\%$ for HDD, $10.14\%$ for MDD, and  $7.00\%$ for LDD trucks. 
It may appear that these effects are larger than the numerical results of earlier studies, for instance \cite{Scora2015, Schroder2019} and \cite{Brunner2021}. We submit that this is due to the long tail of the distribution of fuel savings which is found only with a sufficiently large sample.

To discern the individual effect of road gradient on \cotwo emissions reduction, we fix a speed policy $\varv \in \{\varv^s,\varv^d\}$ and then take into account the \cotwo emissions reduction by traveling along the greenest path $\pi^g(\varv , l)$ rather than the shortest path $\pi^{sp}$.
We consider two ratios $\%\mathcal{E}_{(\pi^{sp}, \varv^s, l)}^{(\pi^g, \varv^s, l)}$ and $\%\mathcal{E}_{(\pi^{sp}, \varv^d, l)}^{(\pi^g, \varv^d, l)}$ to assess this effect.
Figures~\ref{fig:baseSaving_old_old} and \ref{fig:baseSaving_new_new} indicate the distribution and mean of these ratios for base case instances in different cities. Figure~\ref{fig:baseSaving_old_old} demonstrates that the selection of $\pi^g(\varv^s,60\%)$ rather than $\pi^{sp}$ can reduce, on average, $1.76\%$ to $8.15\%$ of the \cotwo emissions, if $l=60\%$ and $\varv^s$ is decided.
As explained before, this \cotwo emissions reduction capacity is lower for the MDD and LDD trucks, yet remains substantive. In the case of a dynamic speed policy $\varv^d$, the statistics, i.e. $\%\mathcal{E}_{(\pi^{sp}, \varv^d, l)}^{(\pi^g, \varv^d, l)}$, remain close to that of $\varv^s$, i.e. $\%\mathcal{E}_{(\pi^{sp}, \varv^s, l)}^{(\pi^g, \varv^s, l)}$, but they are slightly smaller. This implies that regardless of speed, taking into account the road gradient results in significant reductions in CO2 emissions.

Next, to investigate the effect of speed policies on fuel-efficient paths and CO2 emissions reduction, we appraise the carbon reduction by modifying the policy from $(\pi^g, \varv^s, l)$ to $(\pi^g, \varv^d, l)$ for the same truck, i.e. $\%\mathcal{E}_{(\pi^g, \varv^s, l)}^{(\pi^g, \varv^d, l)}$. Figure~\ref{fig:baseSaving_GPP} shows that for most cities, $\overline{\%\mathcal{E}}_{(\pi^g, \varv^s, 60\%)}^{(\pi^g, \varv^d, 60\%)}$ is between $2\%$ and $4\%$, and in all cases the estimates do not depend on the truck type. 

\begin{figure}[b!]
	\centering
  \includegraphics[width=1\textwidth]{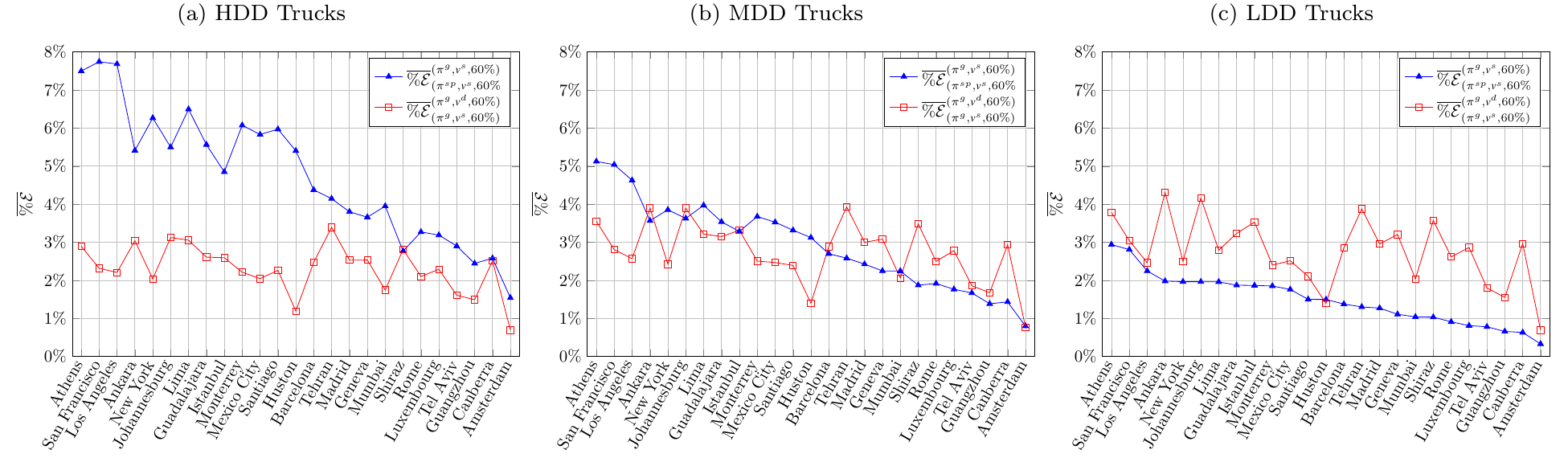}\\
	\caption{\textsf{Comparison of the \cotwo emissions reduction potential: sole $\pi^g$ ($\overline{\%\mathcal{E}}^{(\pi^g, \varv^s , 60\%)}_{(\pi^{sp}, \varv^s , 60\%)}$) vs. sole $\varv^d$ ($\overline{\%\mathcal{E}}^{(\pi^g, \varv^d , 60\%)}_{(\pi^g, \varv^s , 60\%)}$) across different cities.}}
  \label{fig:SlopeSpeedComparison}
\end{figure}

We contrast $\%\mathcal{E}_{(\pi^g, \varv^s, 60\%)}^{(\pi^g, \varv^d, 60\%)}$ and $\%\mathcal{E}_{(\pi^{sp}, \varv^d, 60\%)}^{(\pi^g, \varv^d, 60\%)}$, as shown in Figure~\ref{fig:SlopeSpeedComparison}, in order to experimentally evaluate the relative efficacy of the greenest path and speed optimization in reducing CO2 emissions for each type of vehicle (i.e. HDD, MDD, and LDD). For HDD trucks, the road gradient is more crucial than the dynamic speed policy, whereas the dynamic speed policy has more impact for LDD trucks. The greenest path and dynamic speed policy can bring down the \cotwo emissions of MDD trucks to the same extent.

\begin{figure}[t!]
	\centering
  \includegraphics[width=1\textwidth]{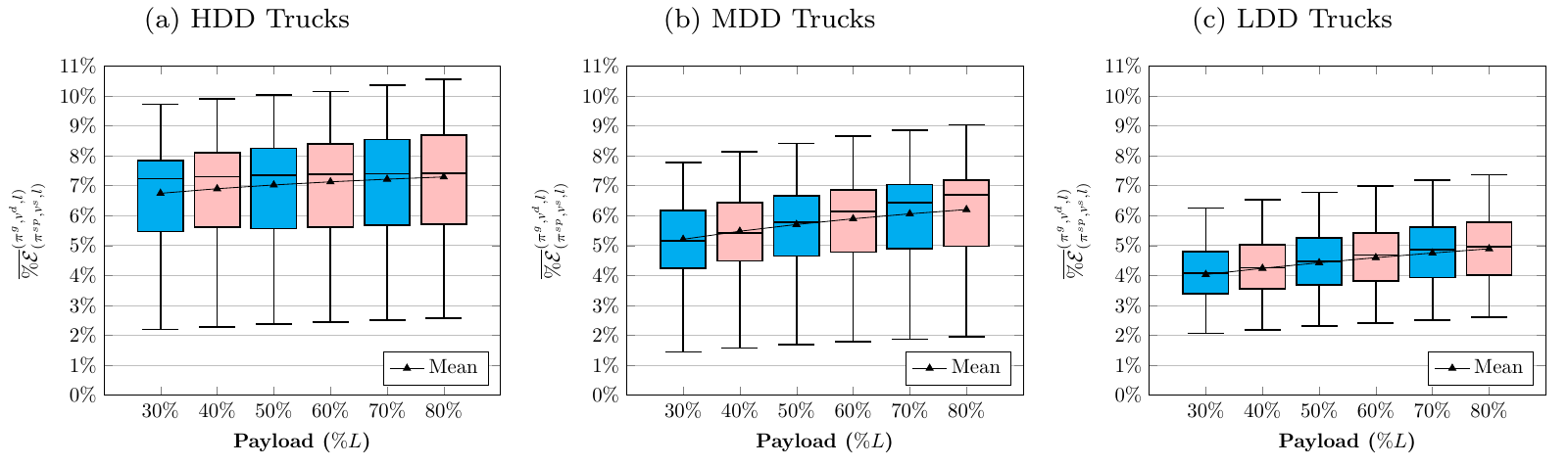}\\
	\caption{\textsf{Effect of payload on $\overline{\%\mathcal{E}}_{(\pi^{sp}, \varv^s, l)}^{(\pi^g, \varv^d, l)}$ across 25 cities.}}
  \label{fig:Saving_new_old_varyingLoad}
\end{figure}

\begin{figure}[t!]
	\centering
  \includegraphics[width=1\textwidth]{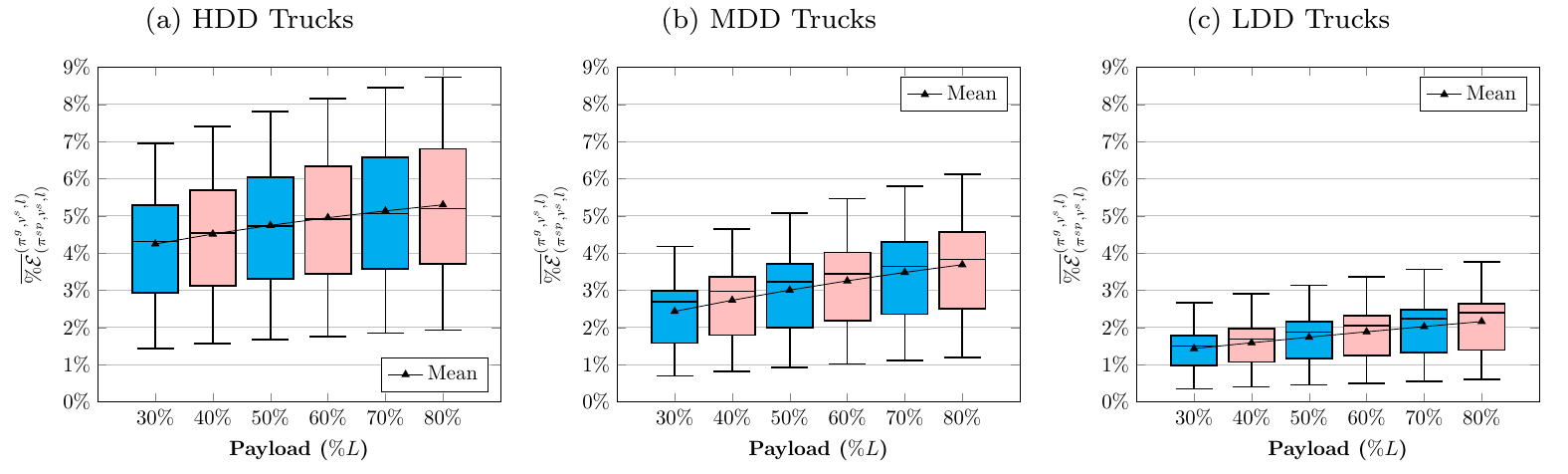}\\
	\caption{\textsf{Effect of payload on $\overline{\%\mathcal{E}}_{(\pi^{sp}, \varv^s, l)}^{(\pi^g, \varv^s, l)}$ across 25 cities.}}
  \label{fig:Saving_old_old_varyingLoad}
\end{figure}

\begin{figure}[t!]
	\centering
  \includegraphics[width=1\textwidth]{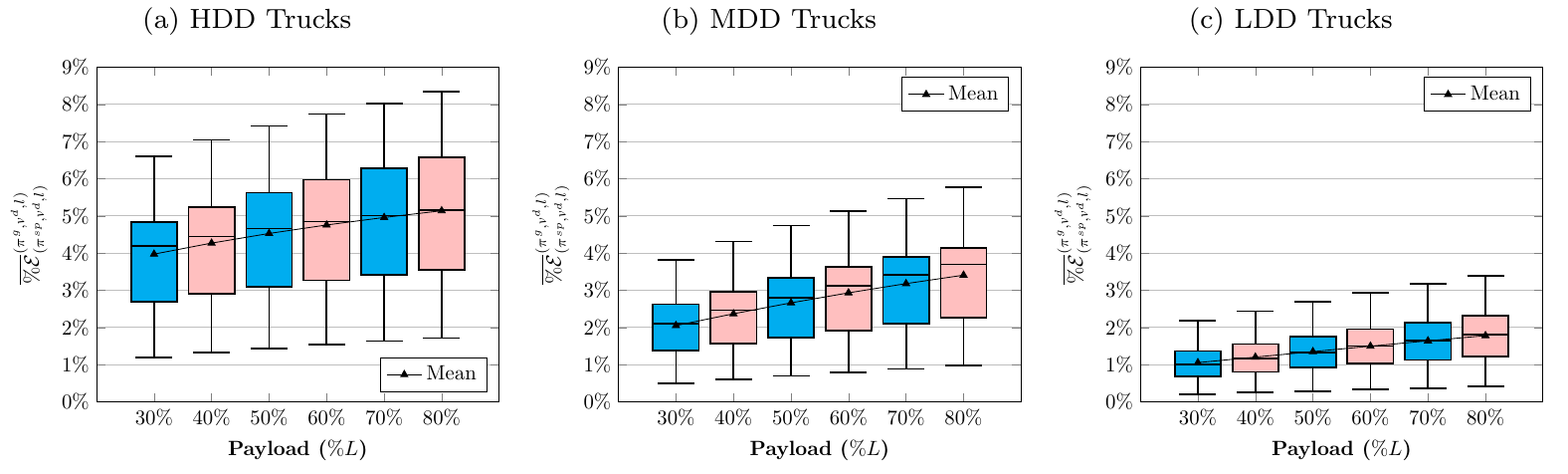}\\
	\caption{\textsf{Effect of payload on $\overline{\%\mathcal{E}}_{(\pi^{sp}, \varv^d, l)}^{(\pi^g, \varv^d, l)}$ across 25 cities.}}
  \label{fig:Saving_new_new_varyingLoad}
\end{figure}

\begin{figure}[htb!]
	\centering
  \includegraphics[width=1\textwidth]{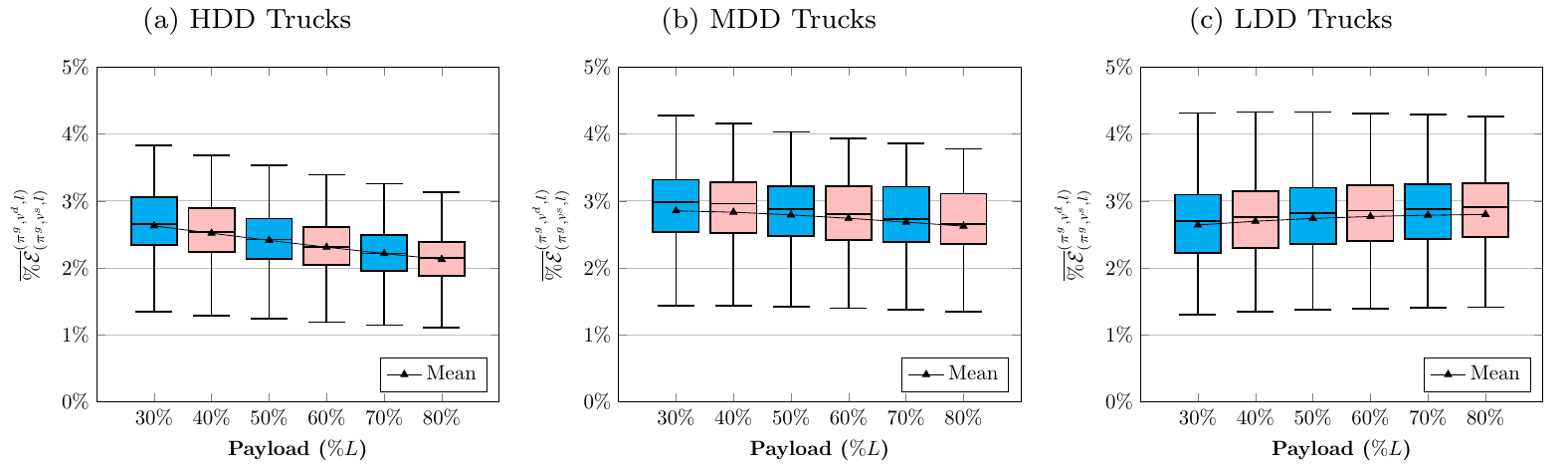}\\
	\caption{\textsf{Effect of payload on $\overline{\%\mathcal{E}}_{(\pi^g, \varv^s, l)}^{(\pi^g, \varv^d, l)}$ across 25 cities.}}
  \label{fig:Saving_GPP_new_old_varyingLoad}
\end{figure}

To analyze the effect of payload on \cotwo emissions reduction, we vary payload ratio for the base case instances ($30\%$, $40\%$, $50\%$, $70\%$, and $80\%$) and repeat the same experiments. Figures~\ref{fig:Saving_new_old_varyingLoad} through \ref{fig:Saving_GPP_new_old_varyingLoad} present the distributions of the sample mean of the relative \cotwo emissions reduction ratios over the 25 cities, where the sample size of each box plot is 25. The figures also present the alteration of the distributions as the payload increases. These results show that, on average,  $\overline{\%\mathcal{E}}_{(\pi^{sp}, \varv^s, l)}^{(\pi^g, \varv^d, l)}$ (Figure~\ref{fig:Saving_new_old_varyingLoad}), $\overline{\%\mathcal{E}}_{(\pi^{sp}, \varv^s, l)}^{(\pi^g, \varv^s, l)}$ (Figure~\ref{fig:Saving_old_old_varyingLoad}), and $\overline{\%\mathcal{E}}_{(\pi^{sp}, \varv^d, l)}^{(\pi^g, \varv^d, l)}$ (Figure~\ref{fig:Saving_new_new_varyingLoad}) are non-decreasing in payload. However, all graphs are concave, so that the growth rate of $\overline{\%\mathcal{E}}$ decreases in payload. In many cities, this phenomenon results in a slow increase, and in one case (Shiraz) slight decrease of $\%\mathcal{E}_{(\pi^{sp}, \varv^s, l)}^{(\pi^g, \varv^d, l)}$ for HDD trucks. The same concave pattern takes place for $\overline{\%\mathcal{E}}_{(\pi^g, \varv^s, l)}^{(\pi^g, \varv^d, l)}$ (Figure~\ref{fig:Saving_GPP_new_old_varyingLoad}) with the exception that the maxima of the concave functions are typically in the MDD or LDD regions. 
This result can explain the close range of $\%\mathcal{E}_{(\pi^g, \varv^s, l)}^{(\pi^g, \varv^d, l)}$ across different truck types as shown in Figure~\ref{fig:baseSaving_GPP}.

\subsection{Results: Paths of the $\pi^g(\varv^d,l)$, $\pi^g(\varv^s,l)$, and $\pi^{sp}$}
\label{subsec:trajectories}
The differences between the greenest path and the shortest path have been covered in earlier sections, along with an analysis of the impact of speed and road gradient. Although our findings indicate significant differences in fuel consumption and CO2 emissions, it is important to consider whether the shortest path's trajectory differs significantly from the trajectory produced by the greenest path.

\begin{figure}[tb!]
	\centering
  \includegraphics[width=.4\textwidth, angle=-90]
  {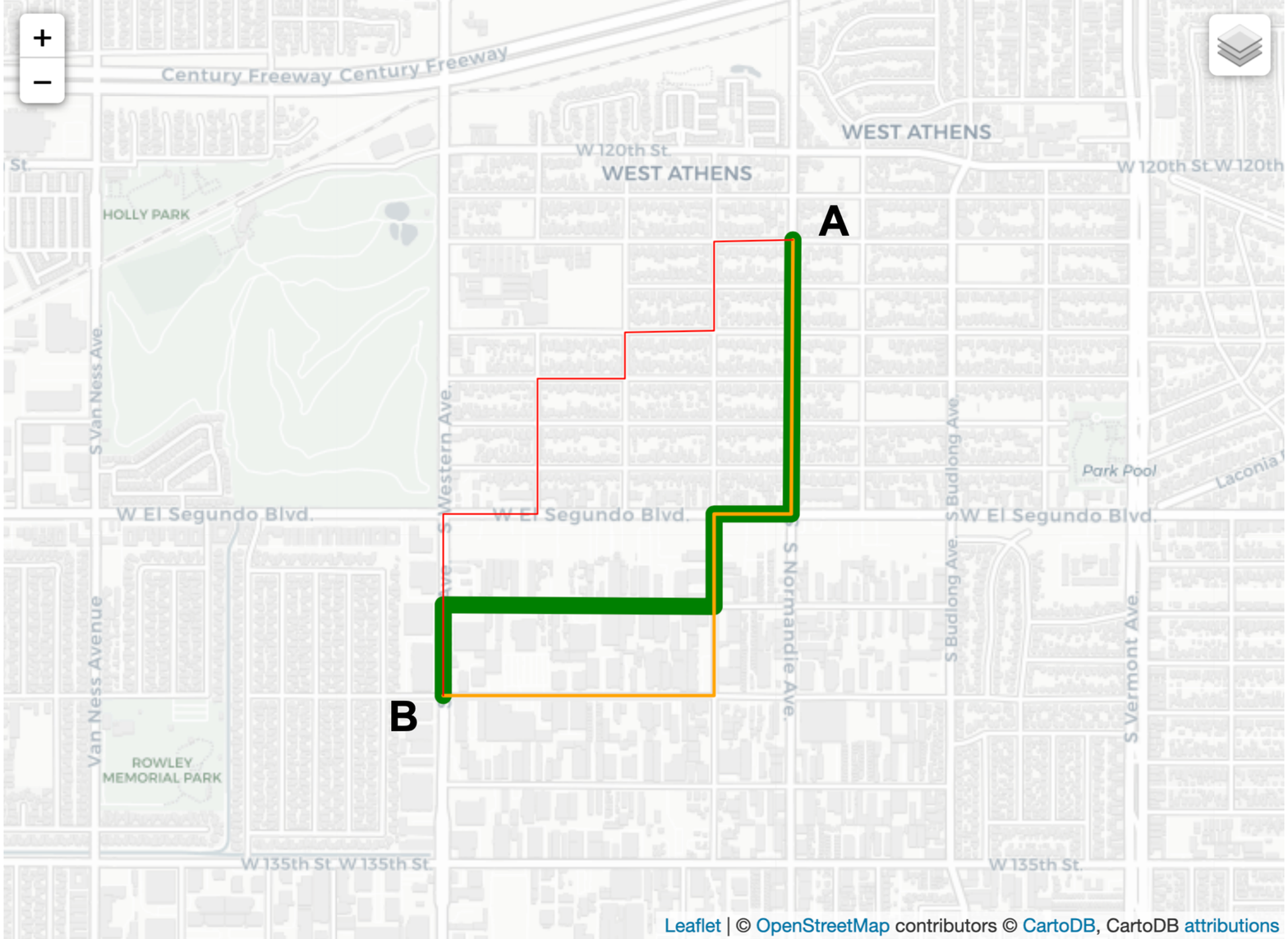}\\
	\caption{\textsf{Example of $\pi^g(\varv^d , 60\%)$, $\pi^g(\varv^s , 60\%)$, $\pi^{sp}$ (Green/Bold: $\pi^g(\varv^d , 60\%)$, Orange/Medium: $\pi^g(\varv^s , 60\%)$, and Red/Thin: $\pi^{sp}$).}}
  \label{fig:GPExample_path}
\end{figure}

\begin{figure}[tb!]
	\centering
  \includegraphics[width=.95\textwidth]{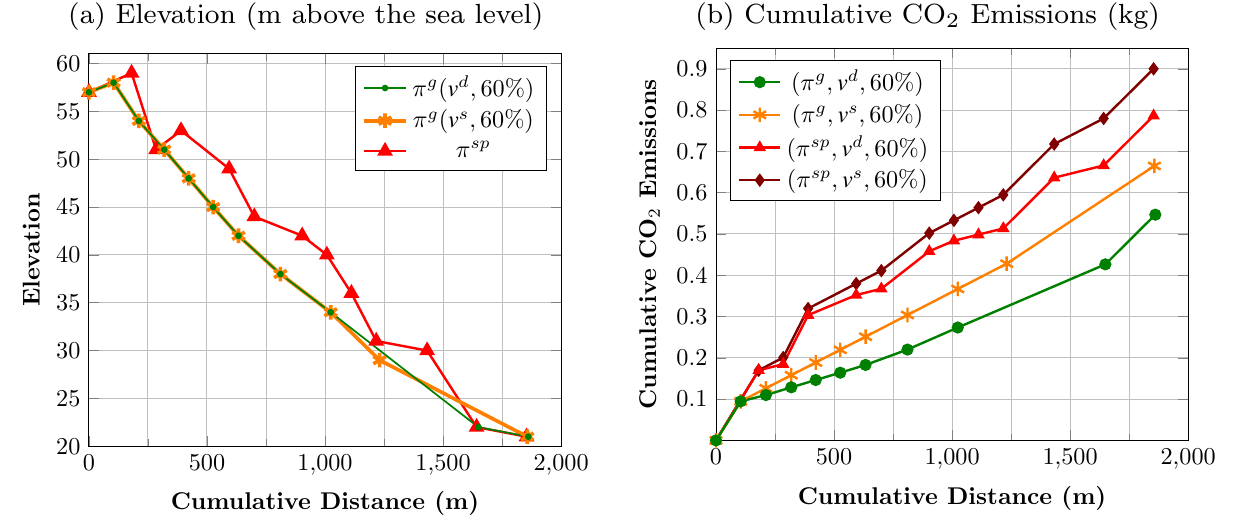}\\
	\caption{\textsf{Elevation of the vertices and \cotwo emissions along the paths as per Figure~\ref{fig:GPExample_path}.}}
  \label{fig:GPExample_emissions}
\end{figure}

To illustrate this difference, we consider a LDD truck that delivers cargo weighing 60\% of its maximum capacity from point A to B within a district of Los Angeles, see Figure~\ref{fig:GPExample_path}.
Figure~\ref{fig:GPExample_path} displays the greenest paths ($\pi^g (\varv^d,l)$ and $\pi^g (\varv^s,l)$) and the shortest path path ($\pi^{sp}$) on the map, and Figure~\ref{fig:GPExample_emissions} shows the elevation of the vertices and total \cotwo emissions for different path and speed choices. In this instance, $\pi^{sp}$ differs significantly from $\pi^g (\varv^d,l)$ and $\pi^g (\varv^s,l)$, whereas the two greenest paths share a number of arcs. In this section, we examine whether such an observation is common throughout our test-bed.

\begin{figure}[t!]
	\centering
  \includegraphics[width=1\textwidth]{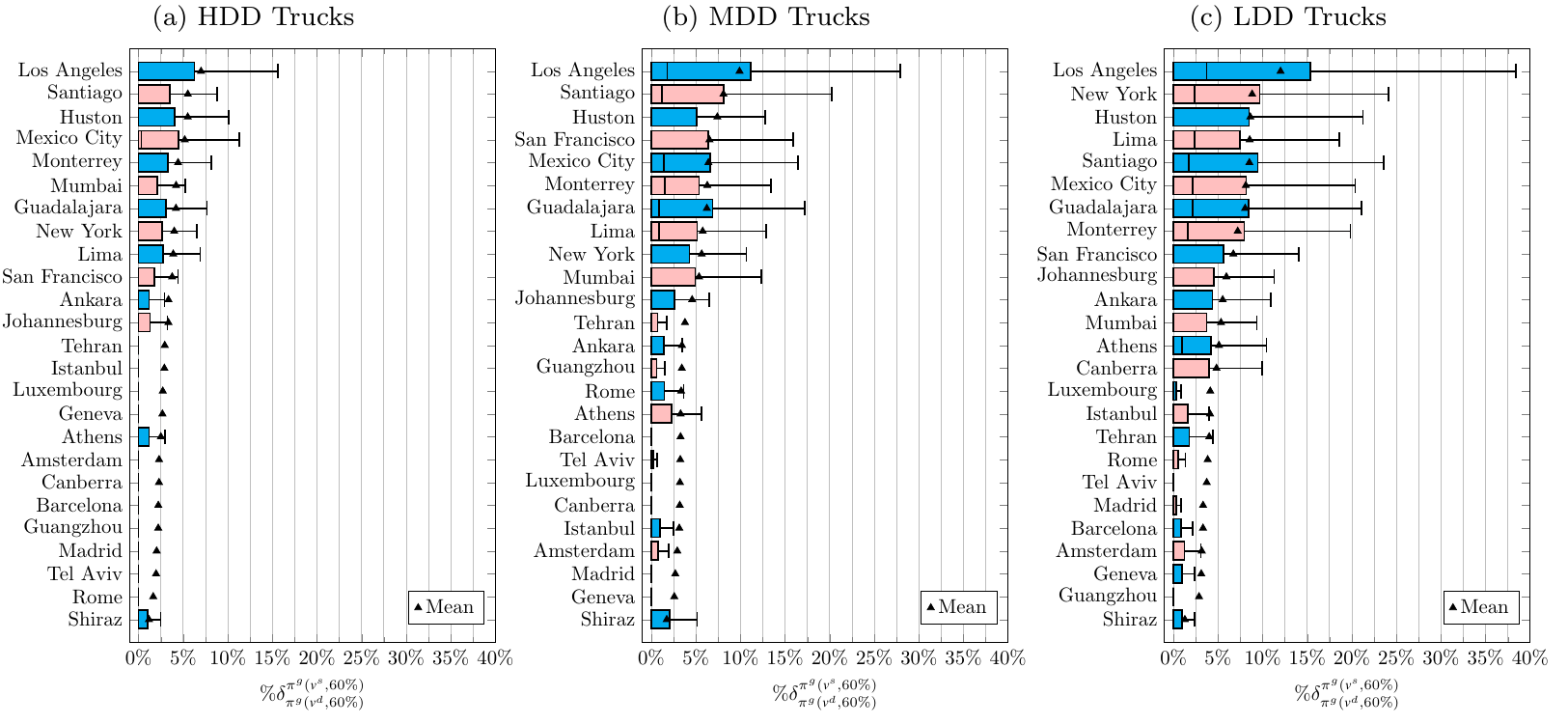}\\
	\caption{\textsf{Ratio of the length of $\pi^g (\varv^d, 60\%)$ that is not shared with $\pi^g (\varv^s, 60\%)$.}}
  \label{fig:baseLength_GPP}
\end{figure}

\begin{figure}[htb!]
	\centering
  \includegraphics[width=1\textwidth]{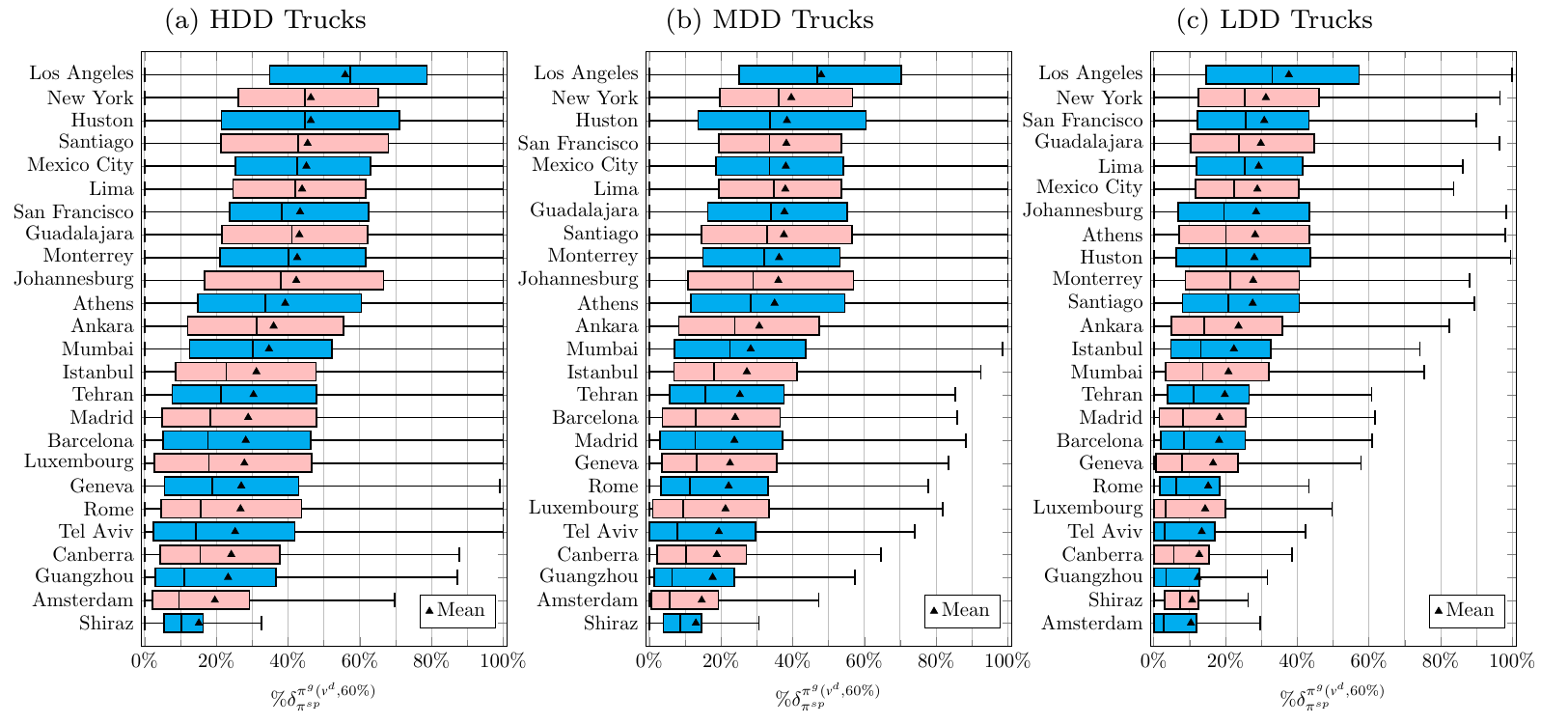}\\
	\caption{\textsf{Ratio of the length of $\pi^{sp}$ that is not shared with $\pi^g (\varv^d, 60\%)$.}}
  \label{fig:baseLength_new}
\end{figure}

\begin{figure}[htb!]
	\centering
  \includegraphics[width=1\textwidth]{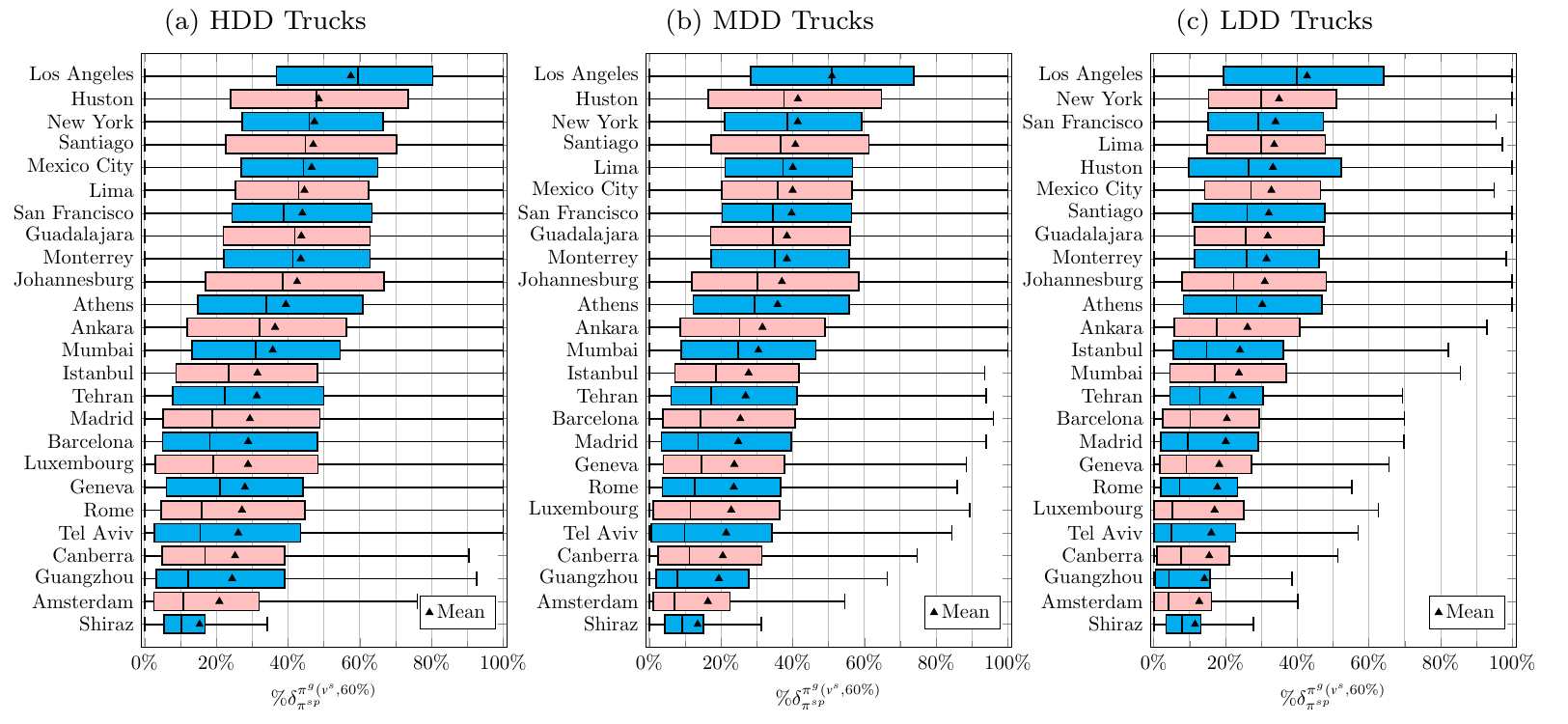}\\
	\caption{\textsf{Ratio of the length of $\pi^{sp}$ that is not shared with $\pi^g (\varv^s, 60\%)$.}}
  \label{fig:baseLength_old}
\end{figure}

Figures~\ref{fig:baseLength_GPP} through \ref{fig:baseLength_old} encapsulate the distribution and sample mean of $\%\delta_{(\pi^g, \varv^d, 60\%)}^{(\pi^g, \varv^s, 60\%)}$, $\%\delta_{\pi^{sp}}^{\pi^g (\varv^d, 60\%)}$, and $\%\delta_{\pi^{sp}}^{\pi^g (\varv^s, 60\%)}$ for the base cases. Figure~\ref{fig:baseLength_GPP} shows that the average difference of $\pi^g (\varv^s, 60\%)$ and $\pi^g (\varv^d, 60\%)$ is between $1.16\%$ and $12.01\%$ across the cities. In fact, the quartiles of $\%\delta_{(\pi^g, \varv^d, 60\%)}^{(\pi^g, \varv^s, 60\%)}$ show that for the most part $\pi^g (\varv^d, 60\%)$ are quite similar to $\pi^g (\varv^s, 60\%)$. In other words, in a majority of instances, the greenest path is independent of the speed policy. Additionally, for heavier trucks the greenest path is less likely to vary as a result of speed optimization. Figures~\ref{fig:baseLength_new} and \ref{fig:baseLength_old} show that the distinction between the shortest and the greenest paths, i.e. $\%\delta_{\pi^{sp}}^{\pi^g(\varv^d,60\%)}$ and $\%\delta_{\pi^{sp}}^{\pi^g(\varv^s,60\%)}$, are conspicuously larger than the distinction between the greenest paths, i.e. $\%\delta_{\pi^g(\varv^d,60\%)}^{\pi^g(\varv^s,60\%)}$. This gap intensifies with heavier truck classes.

\begin{figure}[htb!]
	\centering
  \includegraphics[width=1\textwidth]{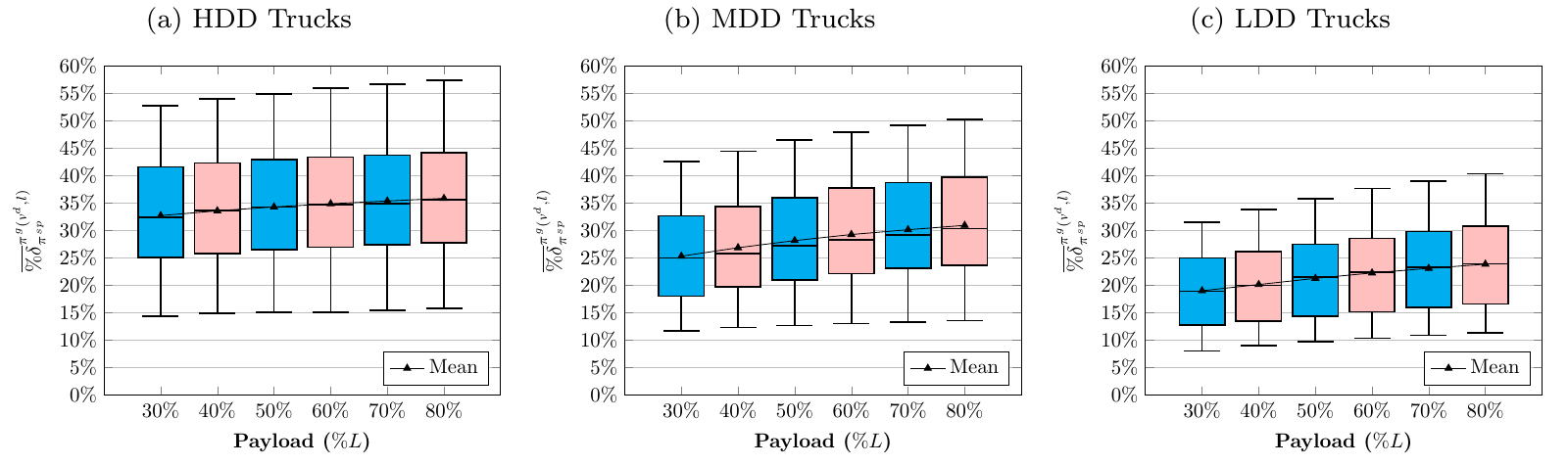}\\
	\caption{\textsf{Effect of payload on $\overline{\%\delta}_{\pi^{sp}}^{\pi^g (\varv^d, l)}$ across 25 cities.}}
  \label{fig:Length_new_varyingLoad}
\end{figure}

\begin{figure}[htb!]
	\centering
  \includegraphics[width=1\textwidth]{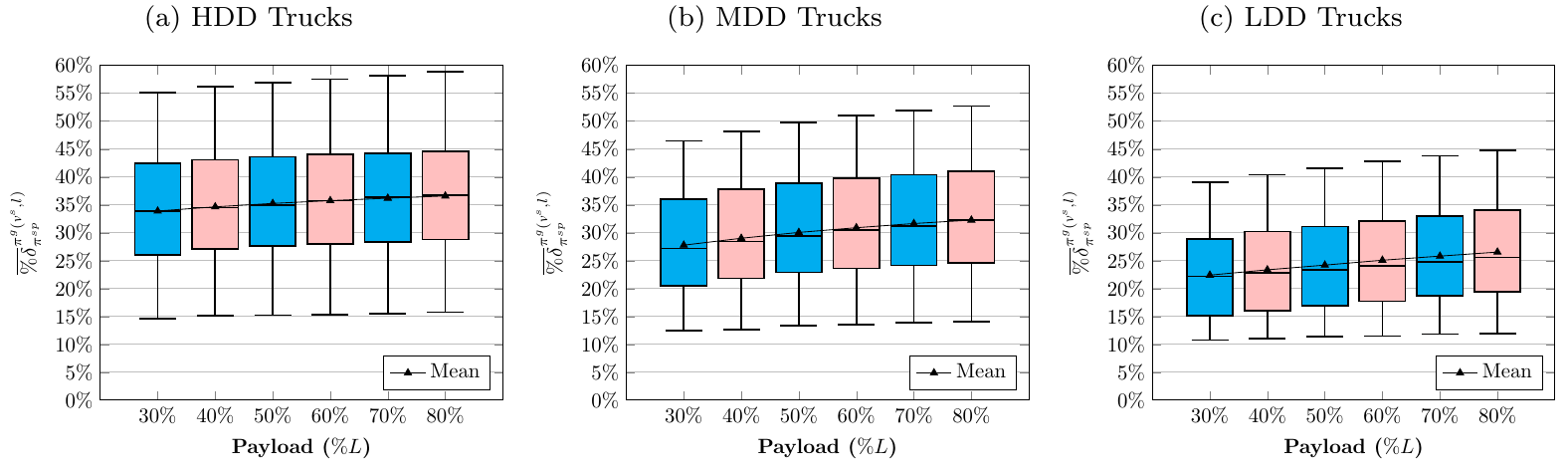}\\
	\caption{\textsf{Effect of payload on $\overline{\%\delta}_{\pi^{sp}}^{\pi^g (\varv^s, l)}$ across 25 cities.}}
  \label{fig:Length_old_varyingLoad}
\end{figure}

\begin{figure}[htb!]
	\centering
  \includegraphics[width=1\textwidth]{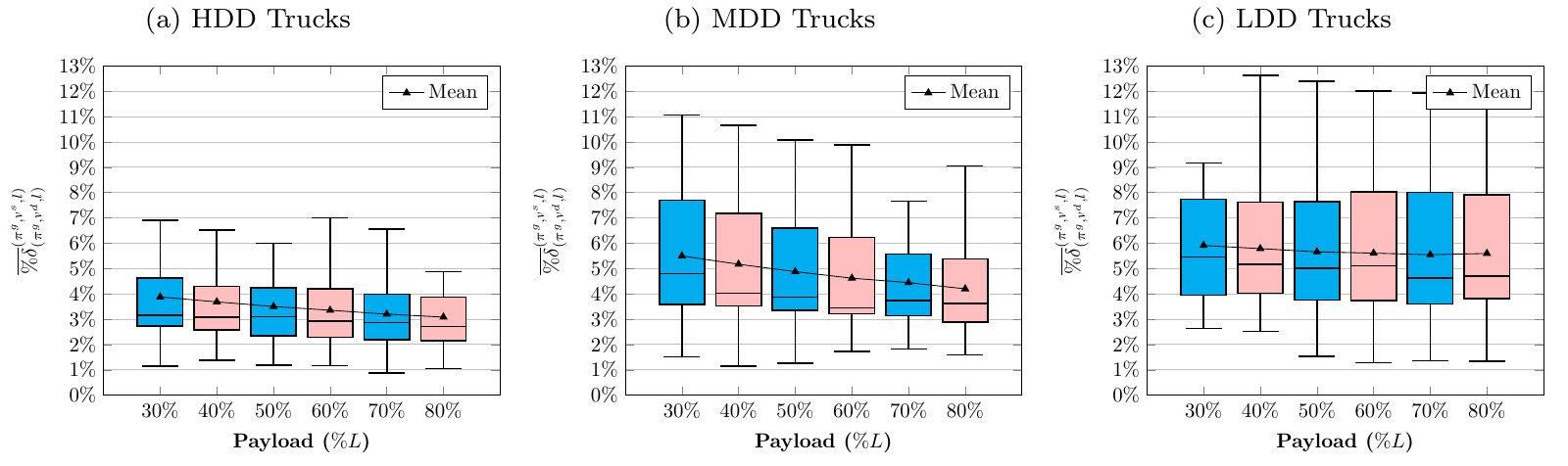}\\
	\caption{\textsf{Effect of payload on $\overline{\%\delta}_{(\pi^g, \varv^d, l)}^{(\pi^g, \varv^s, l)}$ across 25 cities.}}
  \label{fig:Length_GPP_varyingLoad}
\end{figure}

To expand our understanding of the payload's influence on paths, we study whether $\pi^g (\varv^s, l)$ and $\pi^g (\varv^d, l)$ converge to each other and diverge from $\pi^{sp}$ as the payload increases.
Figures~\ref{fig:Length_new_varyingLoad} and \ref{fig:Length_old_varyingLoad} demonstrate that both $\overline{\%\delta}_{\pi^{sp}}^{\pi^g (\varv^d, l)}$ and $\overline{\%\delta}_{\pi^{sp}}^{\pi^g (\varv^s, l)}$ are non-decreasing in payload in contrast to $\overline{\%\delta}_{(\pi^g, \varv^d, l)}^{(\pi^g, \varv^s, l)}$ which is mostly decreasing, as indicated by Figure~\ref{fig:Length_GPP_varyingLoad}. 
Note that, $\overline{\%\delta}_{\pi^{sp}}^{\pi^g (\varv^s, l)}$ is always higher than $\overline{\%\delta}_{\pi^{sp}}^{\pi^g (\varv^d, l)}$ since the more efficient dynamic speed policy of $\pi^g (\varv^d , l)$ usually allows for a shorter (and faster) path relative to $\pi^g (\varv^s , l)$. However, increase in payload erodes the impact of speed policy.

\subsection{Results: Performance of the Asymptotic Greenest Paths}
\label{subsec:AP_num}
\begin{figure}[htb!]
	\centering
  \includegraphics[width=1\textwidth]{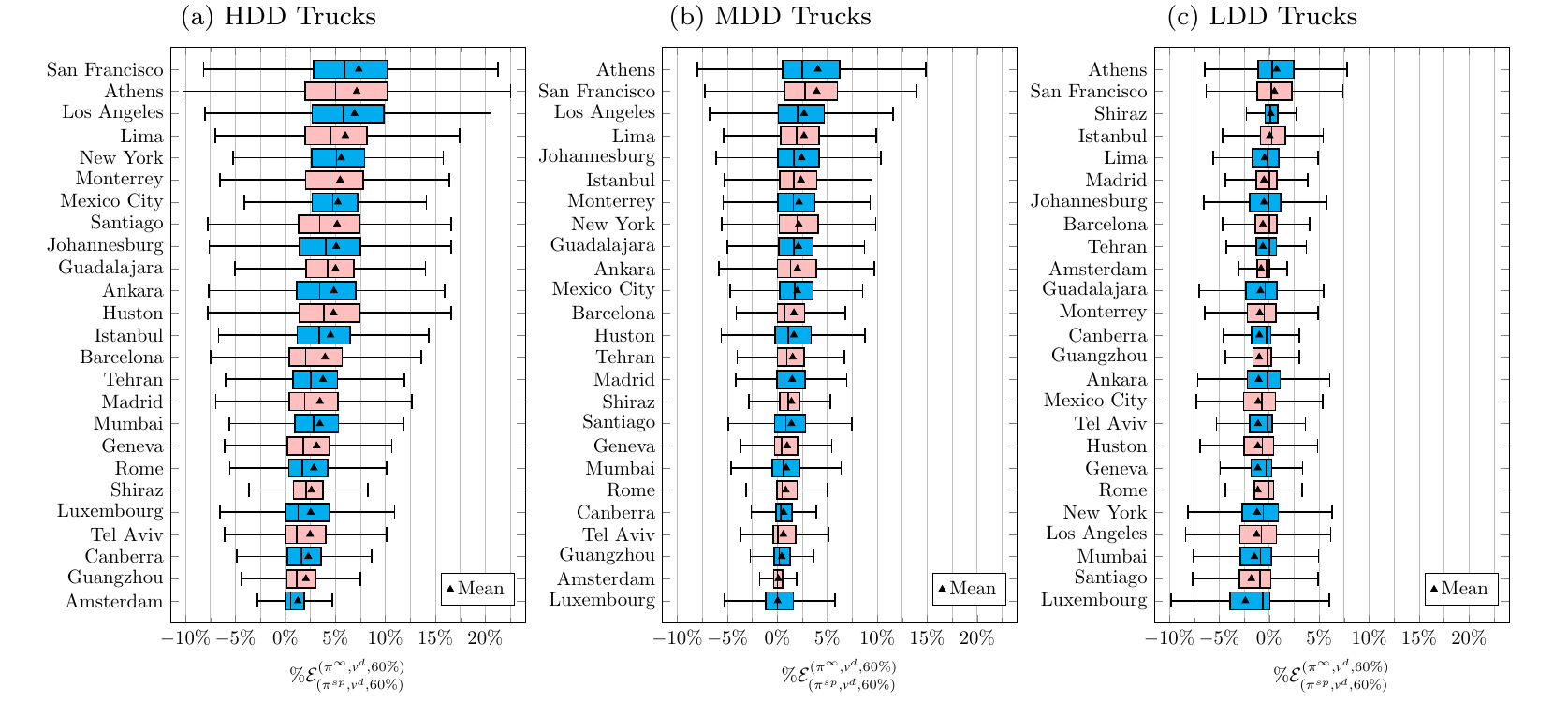}\\
	\caption{\textsf{Relative \cotwo emissions reduction by selecting $(\pi^{\infty}, \varv^d, 60\%)$ rather than $(\pi^{sp}, \varv^d, 60\%)$.}}
  \label{fig:relFCSPS}
\end{figure}

\begin{figure}[htb!]
	\centering
  \includegraphics[width=1\textwidth]{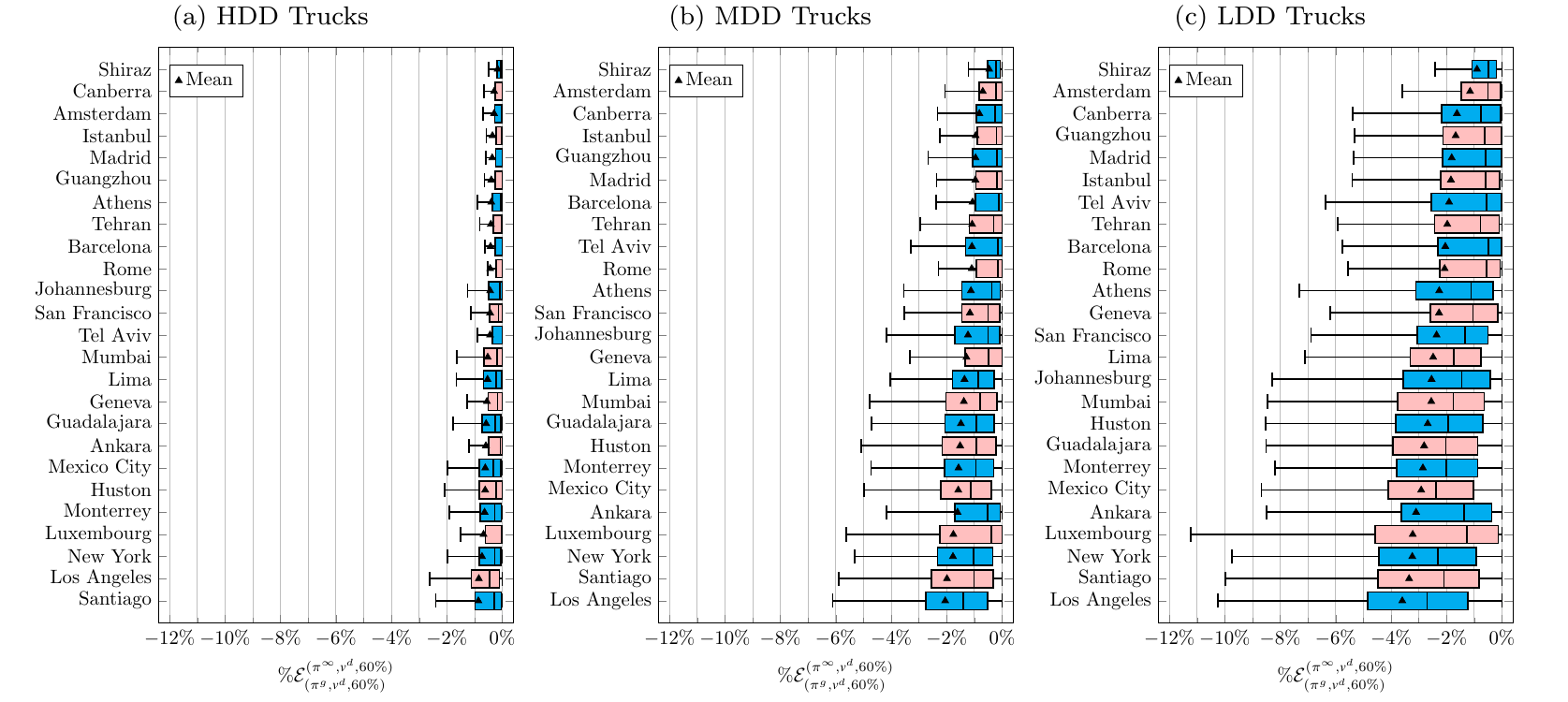}\\
	\caption{\textsf{Relative \cotwo emissions reduction by selecting $(\pi^{\infty}, \varv^d, 60\%)$ rather than $(\pi^g, \varv^d, 60\%)$.}}
  \label{fig:relFCGPS}
\end{figure}

\begin{figure}[htb!]
	\centering
  \includegraphics[width=1\textwidth]{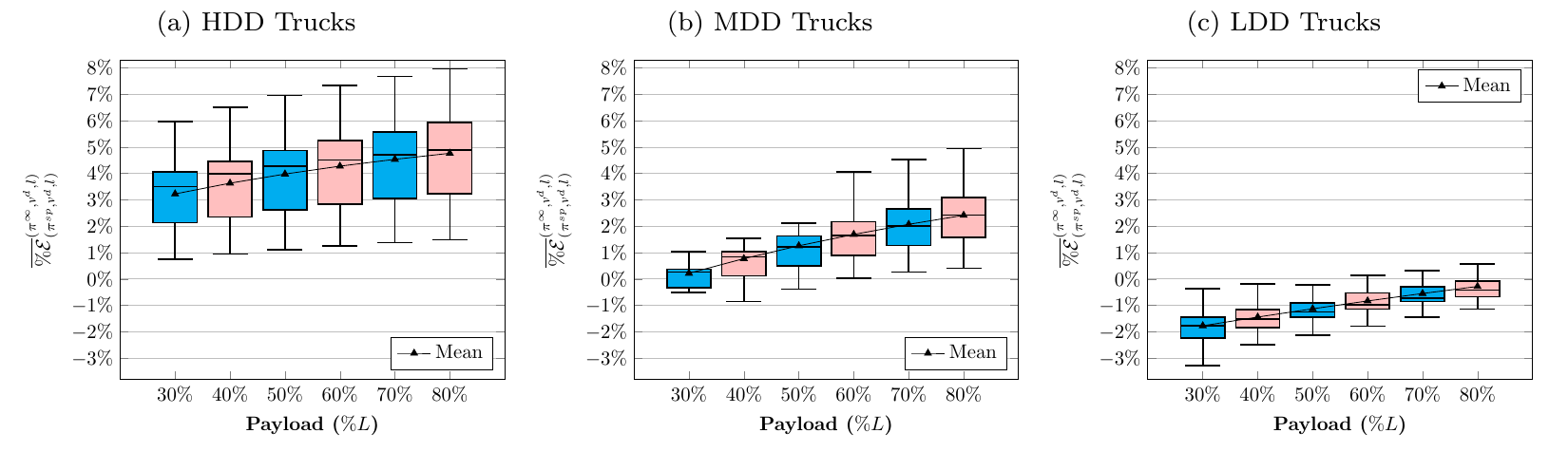}\\
	\caption{\textsf{Effect of payload on $\overline{\%\mathcal{E}}_{(\pi^{sp}, \varv^d, l)}^{(\pi^{\infty}, \varv^d, l)}$ across 25 cities.}}
  \label{fig:varLrelFCSPS}
\end{figure}

\begin{figure}[htb!]
	\centering
  \includegraphics[width=1\textwidth]{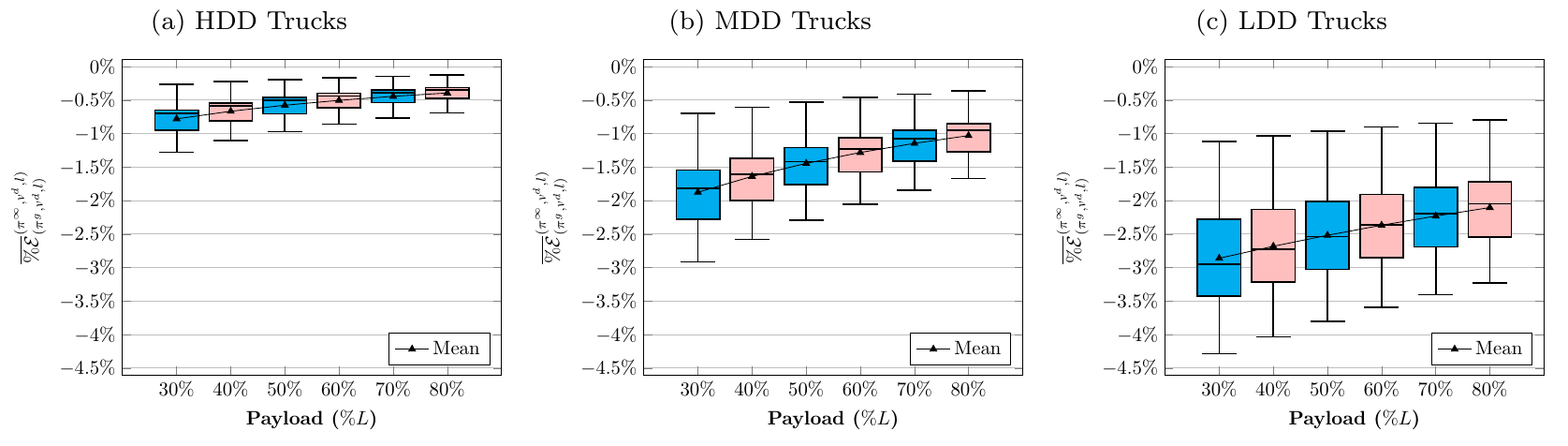}\\
	\caption{\textsf{Effect of payload on $\overline{\%\mathcal{E}}_{(\pi^g, \varv^d, l)}^{(\pi^{\infty}, \varv^d, l)}$ across 25 cities.}}
  \label{fig:varLrelFCGPS}
\end{figure}

The greenest path converges to the asymptotic greenest path for arbitrarily large payloads as shown in Section \ref{subsect:asymcal}. 
In this section, we study the performance of the asymptotic greenest path relative to the shortest path and the greenest path. Then we study the rate of convergence of the greenest path to the asymptotic greenest path for the dynamic speed policy.
In Appendix \ref{section:APC_Performance} we study these things under the static speed policy.
Figures~\ref{fig:relFCSPS} and \ref{fig:relFCGPS} show that the distribution of the \cotwo emissions reduction of $\pi^{\infty}$ relative to $\pi^{sp}$ and $\pi^g(\varv^d, 60\%)$ for different cities. Similar to Section \ref{subsec:emissionreduction}, we present the results for the base cases (60\% payload ratio). 
Figure~\ref{fig:relFCSPS} shows that for the most part an LDD truck emits slightly more \cotwo if it traverses $\pi^{\infty}$ instead of $\pi^{sp}$ in 18 cities. Whereas, the $\pi^{\infty}$ is greener than the $\pi^{sp}$ for MDD and HDD trucks in more than 50\% of the instances in all cities.
The \cotwo emissions reduction of $(\pi^{\infty},\varv^d , 60\%)$ relative to $(\pi^g,\varv^d , 60\%)$, i.e. $\%\mathcal{E}_{(\pi^g, \varv^d, 60\%)}^{(\pi^{\infty}, \varv^d, 60\%)}$, is consistent with this observation. 
Figure~\ref{fig:relFCGPS} shows that the median of extra \cotwo emissions along $\pi^{\infty}$ compared to the $\pi^g(\varv^d,60\%)$ ranges from 0.48\% to 2.70\% for LDD trucks. This range decreases to between 0.11\% and 1.40\% for MDD trucks, and 0\% and 0.45\% for HDD trucks.
Figure~\ref{fig:varLrelFCSPS} shows the distribution of the sample mean of $\%\mathcal{E}_{(\pi^{sp}, \varv^d, l)}^{(\pi^{\infty}, \varv^d, l)}$ across the 25 cities for various payload ratios, i.e. $\overline{\%\mathcal{E}}_{(\pi^{sp}, \varv^d, l)}^{(\pi^{\infty}, \varv^d, l)}$. Correspondingly, Figure~\ref{fig:varLrelFCGPS} presents $\overline{\%\mathcal{E}}_{(\pi^g, \varv^d, l)}^{(\pi^{\infty}, \varv^d, l)}$.
The two figures show that the average \cotwo emissions along $\pi^{\infty}$ relative to $\pi^{sp}$ and $\pi^g (\varv^d,l)$ is non-increasing in load, $l$. Evidently, $\pi^{sp}$ outperforms $\pi^{\infty}$ in terms of average \cotwo emissions for LDD trucks with any payload ratio. Whereas,  $\pi^{\infty}$ is on average greener than $\pi^{sp}$ for MDD and HDD truck types for all payload ratios. The mean excess \cotwo emissions of $\pi^{\infty}$ relative to $\pi^g (\varv^d,l)$ is less than 1\%, 2\%, and 3\% for HDD, MDD, and LDD trucks, respectively.

\begin{figure}[htb!]
	\centering
  \includegraphics[width=1\textwidth]{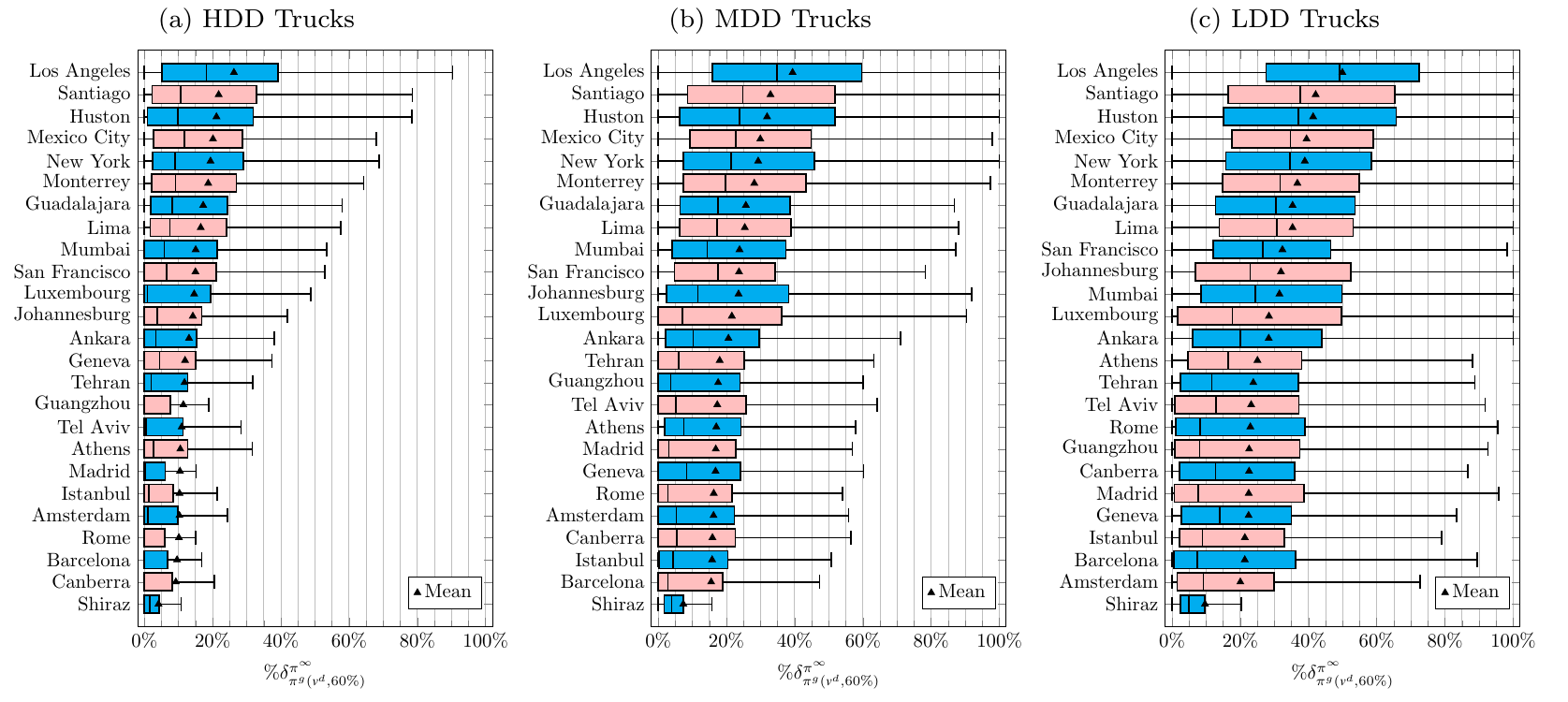}\\
	\caption{\textsf{Ratio of the length of $\pi^g (\varv^d, 60\%)$ that is not shared with $\pi^{\infty}$.}}
  \label{fig:relGPS}
\end{figure}

The median of the difference between $\pi^g (\varv^d, l)$ and $\pi^{\infty}$, i.e. $\%\delta^{\pi^{\infty}}_{\pi^g (\varv^d, 60\%)}$, varies between 4.97\% and 49.14\% for the LDD trucks in base cases as Figure~\ref{fig:relGPS} shows. However, the similarity increases in MDD and HDD truck types as the median $\%\delta^{\pi^{\infty}}_{\pi^g (\varv^d, 60\%)}$ ranges from 2.84\% to 34.81\% for MDD trucks and 0\% to 18.14\% for HDD trucks.
Moreover, the difference between the $\pi_{(\pi^g, \varv^d, l)}$ and $\pi^{\infty}$ reduces in the payload ratio in all truck types.
\begin{figure}[htb!]
	\centering
  \includegraphics[width=1\textwidth]{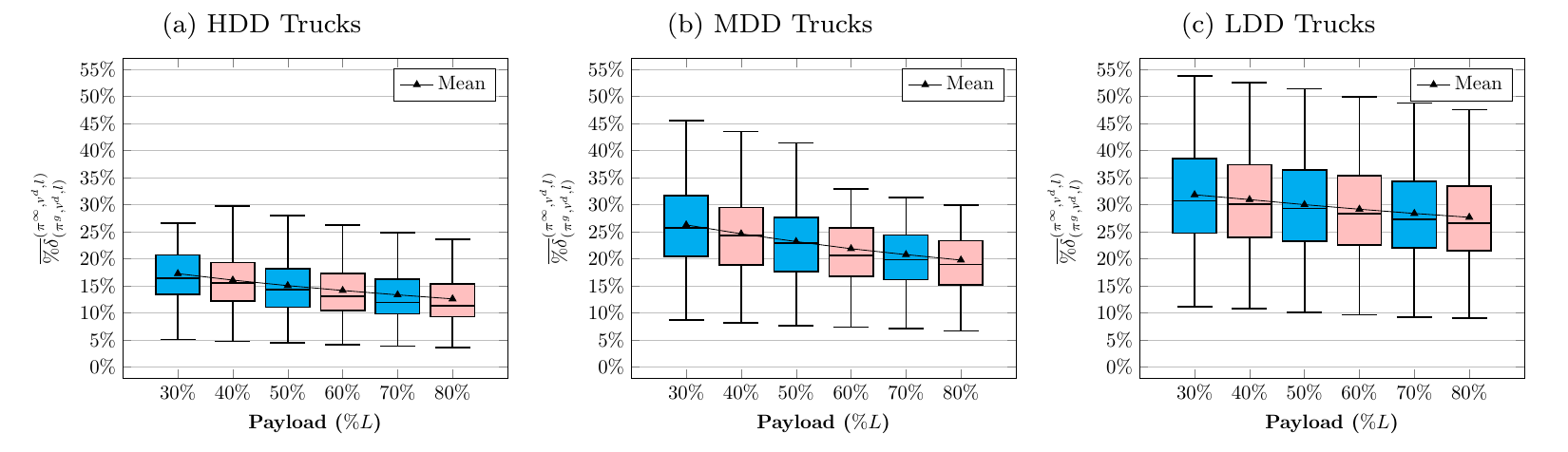}\\
	\caption{\textsf{Effect of payload on $\overline{\%\delta}_{(\pi^g, \varv^d, l)}^{(\pi^{\infty}, \varv^d, l)}$ across 25 cities.}}
  \label{fig:varLrelGPS}
\end{figure}
Figures~\ref{fig:relGPS} and \ref{fig:varLrelGPS} show that $\pi^g ( \varv^d, l)$ converges to $\pi^{\infty}$ in the payload ratio as established in Proposition \ref{prop:AP}.
These results confirm that $\pi^g ( \varv^d, l)$ diverges from $\pi^{sp}$ (Figures~\ref{fig:baseLength_new} and \ref{fig:Length_new_varyingLoad}) and converges to  $\pi^{\infty}$ (Figures~\ref{fig:relGPS} and \ref{fig:varLrelGPS}) as the payload (and curb weight) increases.

\subsection{Results: Main Determinants}
\label{subsec:drivers}
In this section, we address the major determinants of the \cotwo emissions reduction and path alteration. 
We consider the following input features: city, truck type, payload, the elevation difference of source and target ($\Delta h$), the distance of the shortest path ($\delta^{sp}$) and the standard deviation of the gradients along the shortest path ($\sigma^{sp}({\theta})$). 
The latter characterizes the hilliness of the shortest path. All of these features can be efficiently computed.
We use linear regression accompanied by the analysis of variance (ANOVA) to regress these features against seven responses, namely $\%\mathcal{E}_{(\pi^{sp}, \varv^s, l)}^{(\pi^g, \varv^d, l)}$, $\%\mathcal{E}_{(\pi^{sp}, \varv^d, l)}^{(\pi^g, \varv^d, l)}$, $\%\mathcal{E}_{(\pi^{sp}, \varv^s, l)}^{(\pi^g, \varv^s, l)}$, $\%\mathcal{E}_{(\pi^g, \varv^s, l)}^{(\pi^g, \varv^d, l)}$, $\%\delta_{\pi^{sp}}^{\pi^g (\varv^d, l)}$, $\%\delta_{\pi^{sp}}^{\pi^g (\varv^s, l)}$, and $\%\delta_{(\pi^g, \varv^d, l)}^{(\pi^g, \varv^s, l)}$.
We apply min-max normalization for the continuous features and dummy encode the categorical features. The encoding removes the redundant dummy features including Canberra and HDD among cities and trucks, respectively. We use the type III sum of squares in the ANOVA. The full report is available in Appendix \ref{section:ANOVA}.

\begin{table*}[!htbp]
    \centering
    \caption{\textsf{Summary of the linear regression and ANOVA for seven different responses.}}
    \fontsize{8pt}{9pt}\selectfont
    \scriptsize
    \label{tab:feature_rank}
    \begin{tabularx}{.83\textwidth}{l l r|l r|l r|l r|l r|l r|l r|l r|l r|l r|l r|l r|l}
        \toprule
         && \multicolumn{2}{c}{$\%\mathcal{E}_{(\pi^{sp}, \varv^s, l)}^{(\pi^g, \varv^d, l)}$} & \multicolumn{2}{c}{$\%\mathcal{E}_{(\pi^{sp}, \varv^s, l)}^{(\pi^g, \varv^s, l)}$} & \multicolumn{2}{c}{$\%\mathcal{E}_{(\pi^{sp}, \varv^d, l)}^{(\pi^g, \varv^d, l)}$} & \multicolumn{2}{c}{$\%\mathcal{E}_{(\pi^g, \varv^s, l)}^{(\pi^g, \varv^d, l)}$} & \multicolumn{2}{c}{$\%\delta_{\pi^{sp}}^{\pi^g (\varv^s, l)}$} & \multicolumn{2}{c}{$\%\delta_{\pi^{sp}}^{\pi^g (\varv^d, l)}$} & \multicolumn{2}{c}{$\%\delta_{(\pi^g, \varv^s, l)}^{(\pi^g, \varv^s, l)}$}\\
        \cmidrule{3-16}
        $Features$  & $df$  & R & S & R & S & R & S & R & S & R & S  & R & S & R & S\\
	    \midrule
	    $\sigma^{sp}(\theta)$  & 1     & 2 & $+$ & 1 & $+$ & 1 & $+$ & 3 & $+$ & 1 & $+$ & 1 & $+$ & 1 & $+$\\
        $\Delta h$              & 1     & 1 & $-$ & 2 & $-$ & 2 & $-$ & 1 & $-$ & 4 & $-$ & 4 & $-$ & 4 & $-$\\
        $\delta^{sp}$           & 1     & 7 & $-$ & 7 & $+$ & 7 & $+$ & 4 & $-$ & 2 & $+$ & 2 & $+$ & 2 & $+$\\
        $l$                     & 1     & 6 & $+$ & 3 & $+$ & 3 & $+$ & 5 & $-$ & 6 & $+$ & 6 & $+$ & 5 & $-$\\
        City                    & 24    & 4 & $\pm$   & 4 & $\pm$   & 5 & $\pm$   & 6 & $\pm$   & 3 & $\pm$    & 3 & $\pm$  & 3 & $\pm$\\
        Truck                   & 2     & 5 & $-$ & 5 & $-$ & 4 & $-$ & 7 & $\pm$   & 5 & $-$ & 5 & $-$ & 7 & $+$\\
        (Intercept)             & 1     & 3 & $-$ & 6 & $-$ & 6 & $-$ & 2 & $-$   & 7 & $-$ & 7 & $-$ & 6 & $+$\\
        \bottomrule
        \multicolumn{16}{l}{\tiny R: Feature's ranking in ANOVA}\\
        \multicolumn{16}{l}{\tiny S: Sign of the feature's weight in linear regression}\\
  \end{tabularx}
\end{table*}

Table \ref{tab:feature_rank} summarizes the ranking and sign of different features in the ANOVA as per Appendix \ref{section:ANOVA}. By the results, $\sigma^{sp}(\theta)$, i.e. the standard deviation of road gradient along the $\pi^{sp}$ has the most explanatory power for \cotwo reduction capacity.
In addition, $\sigma^{sp}(\theta)$ is has the strongest association with the dissimilarity of the greenest and shortest paths. That is to say, a higher $\sigma^{sp}(\theta)$ indicates a higher potential of \cotwo emissions reduction by selecting the greenest path instead of the shortest path. 
Next comes difference in elevation between the target and the source, $\Delta h$, which is negatively associated with the \cotwo emissions reduction capacity. This relation is strongest when comparing the dynamic speed policy with the static speed policy as in $\%\mathcal{E}_{(\pi^g,\varv^s,l)}^{(\pi^g,\varv^d,l)}$. 
This implies that using elevation data in routing policies is more pivotal for downward trips. Table \ref{tab:feature_rank} also reveals the positive association of $\%\mathcal{E}_{(\pi^{sp}, \varv^d,l)}^{(\pi^g, \varv^d,l)}$ and $\%\mathcal{E}_{(\pi^{sp}, \varv^s,l)}^{(\pi^g, \varv^s,l)}$ with payload.
Our analysis shows that relative dissimilarity of the shortest and greenest paths increases in distance of the shortest path, i.e. $\delta^{sp}$. 
However, $\delta^{sp}$ is less important for \cotwo emissions reduction.
A city's individual characteristics have a fair impact on the \cotwo emissions reduction capacity, albeit this effect is not comparable with that of $\sigma^{sp}(\theta)$ and $\Delta h$. 
Finally, the truck type has an effect that is similar to the payload.
It follows that curb weight and payload of truck are more significant than other parameters for the \cotwo emissions reduction.

\begin{rev}

\section{Numerical Experiments with Traffic Information}
\label{section:NumericalExperimentTraffic}
In this section, we examine how the simultaneous utilization of elevation and traffic data can guide routing decisions with lower \cotwo emissions over a large dataset.
We compare the \cotwo emissions of different types of trucks traveling along three types of routes: the greenest path, $\pi^g(\varv,l)$, the asymptotic greenest path, $\pi^{\infty}(\varv)$, and the path with minimum possible travel duration, the fastest path $\pi^{fp}$.
These comparisons are made under three different speed decisions: traffic speed $\varv^f$, dynamic speed $\varv^d$, and static speed $v^s$.
For all arcs $a \in A$, the maximum speed $\varv^{\max}$ is set to the traffic speed $\varv^f$ and the minimum speed $\varv^{\min}$ is set to zero.

We consider a strongly connected subgraph of New York city's road network for our experiments.
The subgraph comprises 39,143 arcs and 23,091 vertices.
Similar to Section \ref{section:NumericalExperiments}, we obtain road network data from OpenStreetMap \citep{OpenStreetMap} and elevation data from the \cite{USGS}'s SRTM 1 Arc-Second Global datasets.
The specifications of the trucks used in the study are provided in Table \ref{tab:trucks}.
Since traffic speed information is not publicly available for all arcs, we calculate traffic speeds using travel distance and duration inquiries from Google's Distance Matrix API.
We select a time point with anticipated heavy traffic, particularly Wednesday, October 9, 2024, at 7:00 a.m., and set the traffic model to ``best-guess".
Given that $t_a(\varv^f(a))$ is the time to traverse arc $a \in A$ with traffic speed of arc $a$, i.e. $\varv^f(a)$, one can compute $\varv^f(a)$ by, $\varv^f(a) = \delta(a)/t_a(\varv^f(a))$.

We randomly select 20,098 unique pairs of non-identical source and target vertices.
For each pair of source and target and each pair of path-speed policies $d_i = (\pi_i, \varv_i ,l), i =1,2$, we compute three metrics including the relative \cotwo reduction $\% \mathcal{E}^{d_2}_{d_1}$, the relative path distinction $\% \delta^{\pi_2}_{\pi_1}$ and the relative time increase of selecting $d_2$ instead of $d_1$, $\%t^{d_2}_{d_1}$, defined by,
\[
    \%t^{d_2}_{d_1} = 100 \cdot \frac{\sum_{a \in \pi_2} t_a(\varv_2) - \sum_{a \in \pi_1} t_a(\varv_1)}{\sum_{a \in \pi_1} t_a(\varv_1)}.
\]

Table \ref{tab:compstudtraffic} briefly summarizes the additional ratios that we use in our studies under traffic conditions.
\begin{table*}[!htbp]
  \caption{\textsf{List of ratios used in the comparative studies in addition to Table \ref{tab:compstud}.}}
	\fontsize{8pt}{9pt}\selectfont
\scriptsize
  \label{tab:compstudtraffic}
  \begin{tabularx}{\textwidth}{l X}
    \toprule
    Ratio & Description\\
		\midrule
    $\%\mathcal{E}_{(\pi^{fp}, \varv^f, l)}^{(\pi^g, \varv^d, l)}$ \vspace{.2cm} &
    Relative \cotwo emissions reduction by selecting the greenest path with the dynamic speed policy relative to the fastest path with the traffic speed given the load $l$.\\
    $\%\mathcal{E}_{(\pi^{fp}, \varv^f, l)}^{(\pi^g, \varv^s, l)}$ \vspace{.2cm} &
    Relative \cotwo emissions reduction by selecting the greenest path with the static speed policy relative to the fastest path with the traffic speed given the load $l$.\\
    $\%\mathcal{E}_{(\pi^{fp}, \varv^d, l)}^{(\pi^g, \varv^d, l)}$ \vspace{.2cm} &
    Relative \cotwo emissions reduction by selecting the greenest path with the dynamic speed policy relative to the fastest path with the dynamic speed policy given the load $l$.\\
    $\%\mathcal{E}_{(\pi^{fp}, \varv^s, l)}^{(\pi^g, \varv^s, l)}$ \vspace{.2cm} &
    Relative \cotwo emissions reduction by selecting the greenest path with the static speed policy relative to the fastest path with the static speed policy given the load $l$.\\
    $\%\mathcal{E}_{(\pi^{fp}, \varv^d, l)}^{(\pi^g(\varv^f,l), \varv^d, l)}$ \vspace{.2cm} &
    Relative \cotwo emissions reduction by selecting the greenest path under the assumption of driving at traffic speed but using the dynamic speed policy upon path traversal relative to the fastest path with the dynamic speed policy given the load $l$.\\
    $\%\mathcal{E}_{(\pi^g(\varv^f,l), \varv^d, l)}^{(\pi^g, \varv^d, l)}$ \vspace{.2cm} &
    Relative \cotwo emissions reduction by selecting the greenest path with the dynamic speed policy relative to the greenest path under the assumption of driving at traffic speed but using the dynamic speed policy upon path traversal given the load $l$.\\
    $\%\mathcal{E}_{(\pi^{fp}, \varv^d, l)}^{(\pi^{\infty}, \varv^d, l)}$ \vspace{.2cm} &
    Relative \cotwo emissions reduction by selecting the asymptotic greenest path with the dynamic speed policy relative to the fastest path with the static speed policy given the load $l$.\\
    $\%\delta_{\pi^{fp}}^{\pi^g (\varv^d, l)}$ \vspace{.2cm} &
    Ratio of the length of the fastest path that is not shared with the greenest path under the dynamic speed policy given the load $l$.\\
    $\%\delta_{\pi^{fp}}^{\pi^g (\varv^s, l)}$ \vspace{.2cm} &
    Ratio of the length of the fastest path that is not shared with the greenest path under the static speed policy given the load $l$.\\
    $\%t_{(\pi^{fp}, \varv^f, l)}^{(\pi^g, \varv^d, l)}$ \vspace{.2cm} &
    Relative time increase by selecting the greenest path with the dynamic speed policy relative to the fastest path with the traffic speed given the load $l$.\\
    $\%t_{(\pi^{fp}, \varv^d, l)}^{(\pi^g, \varv^d, l)}$ \vspace{.2cm} &
    Relative time increase by selecting the greenest path with the dynamic speed policy relative to the fastest path with the dynamic speed policy given the load $l$.\\
    $\%t_{(\pi^{fp}, \varv^f, l)}^{(\pi^g, \varv^s, l)}$ \vspace{.2cm} &
    Relative time increase by selecting the greenest path with the static speed policy relative to the fastest path with the traffic speed given the load $l$.\\
    $\%t_{(\pi^{fp}, \varv^s, l)}^{(\pi^g, \varv^s, l)}$ \vspace{.2cm} &
    Relative time increase by selecting the greenest path with the static speed policy relative to the fastest path with the static speed policy given the load $l$.\\
    \bottomrule
  \end{tabularx}
\end{table*}

\subsection{Results: Impact of Path and Speed Decisions on \cotwo Emissions Reduction}
\label{subsec:emissionreduction_traffic}

Figure \ref{fig:traffic_fpvf_gvd} demonstrates that $\overline{\%\mathcal{E}}^{(\pi^g, \varv^d, l)}_{(\pi^{fp}, \varv^f ,l)}$ ranges from 19.40\% to 26.02\% for different truck types and payloads, whereas Figure \ref{fig:traffic_fpvd_gvd} shows that $\overline{\%\mathcal{E}}^{(\pi^g, \varv^d, l)}_{(\pi^{fp}, \varv^d ,l)}$ varies between 5.80\% and 8.91\%.
Similar statistics are observed for $\overline{\%\mathcal{E}}^{(\pi^g, \varv^s, l)}_{(\pi^{fp}, \varv^s ,l)}$ as illustrated in Figure \ref{fig:traffic_fpvs_gvs}.
\begin{figure}[htb!]
	\centering
  \includegraphics[width=1\textwidth]{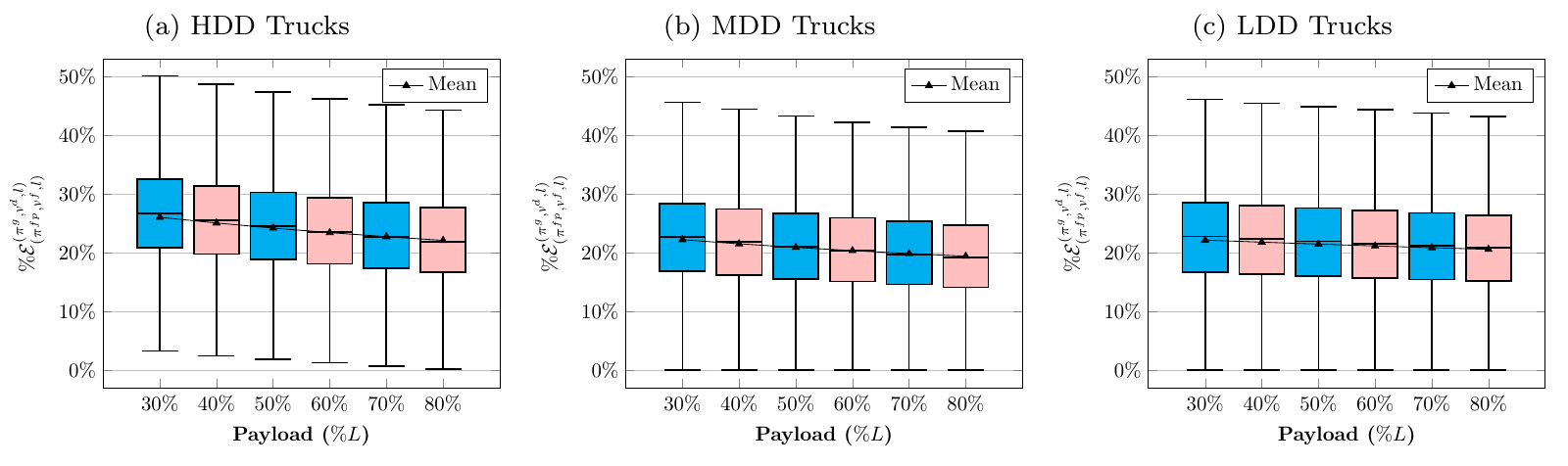}\\
	\caption{\textsf{$\%\mathcal{E}^{(\pi^g, \varv^d, l)}_{(\pi^{fp}, \varv^f ,l)}$ across truck types and payloads in traffic condition.}}
  \label{fig:traffic_fpvf_gvd}
\end{figure}
\begin{figure}[htb!]
	\centering
  \includegraphics[width=1\textwidth]{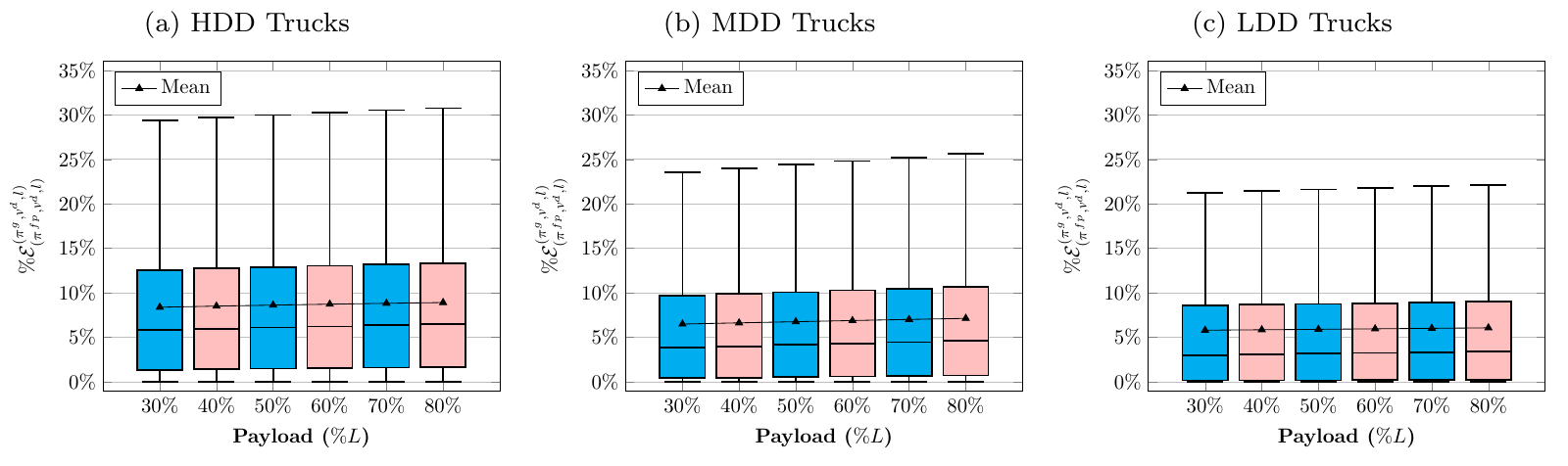}\\
	\caption{\textsf{$\%\mathcal{E}^{(\pi^g, \varv^d, l)}_{(\pi^{fp}, \varv^d ,l)}$ across truck types and payloads in traffic condition.}}
  \label{fig:traffic_fpvd_gvd}
\end{figure}
\begin{figure}[htb!]
	\centering
  \includegraphics[width=1\textwidth]{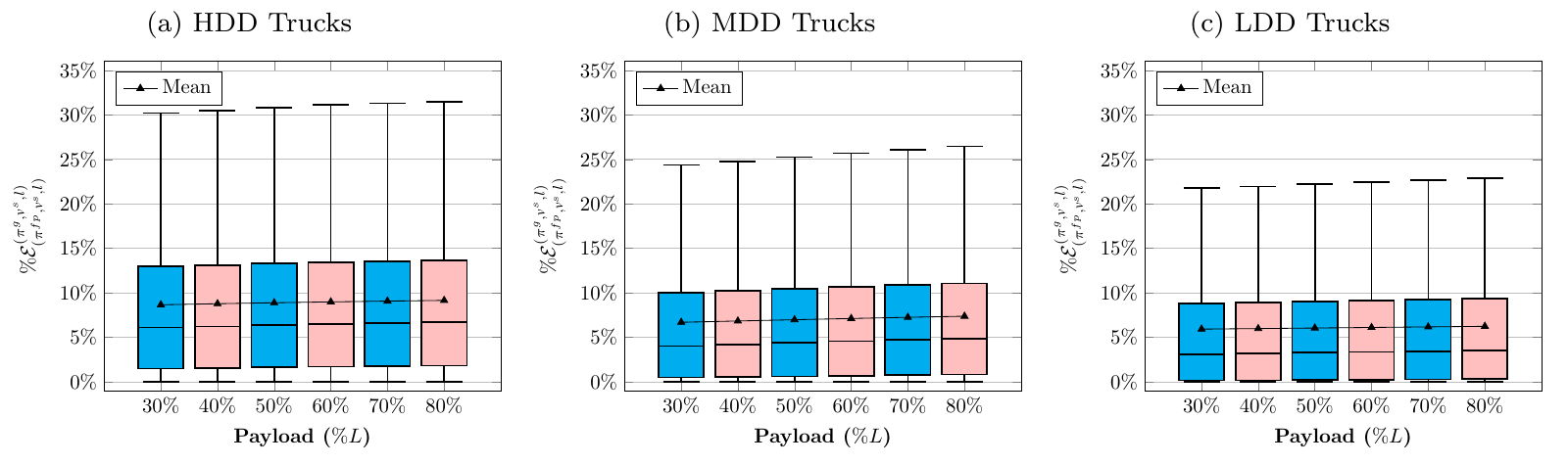}\\
	\caption{\textsf{$\%\mathcal{E}^{(\pi^g, \varv^s, l)}_{(\pi^{fp}, \varv^s ,l)}$ across truck types and payloads in traffic condition.}}
  \label{fig:traffic_fpvs_gvs}
\end{figure}
These results provide significant evidence that both path selection and speed optimization can contribute to reducing \cotwo emissions in intra-city truck transportation.
Additionally, the reduction potential in \cotwo emissions is greater on the greenest path during traffic conditions compared to free-flow situations (cf. Figures \ref{fig:baseSaving_old_old} and \ref{fig:baseSaving_new_new}).
However, the potential reduction in \cotwo emissions through a dynamic speed policy versus a static speed policy is negligible in most instances, as illustrated in Figure \ref{fig:traffic_gvs_gvd}.
\begin{figure}[htb!]
	\centering
  \includegraphics[width=1\textwidth]{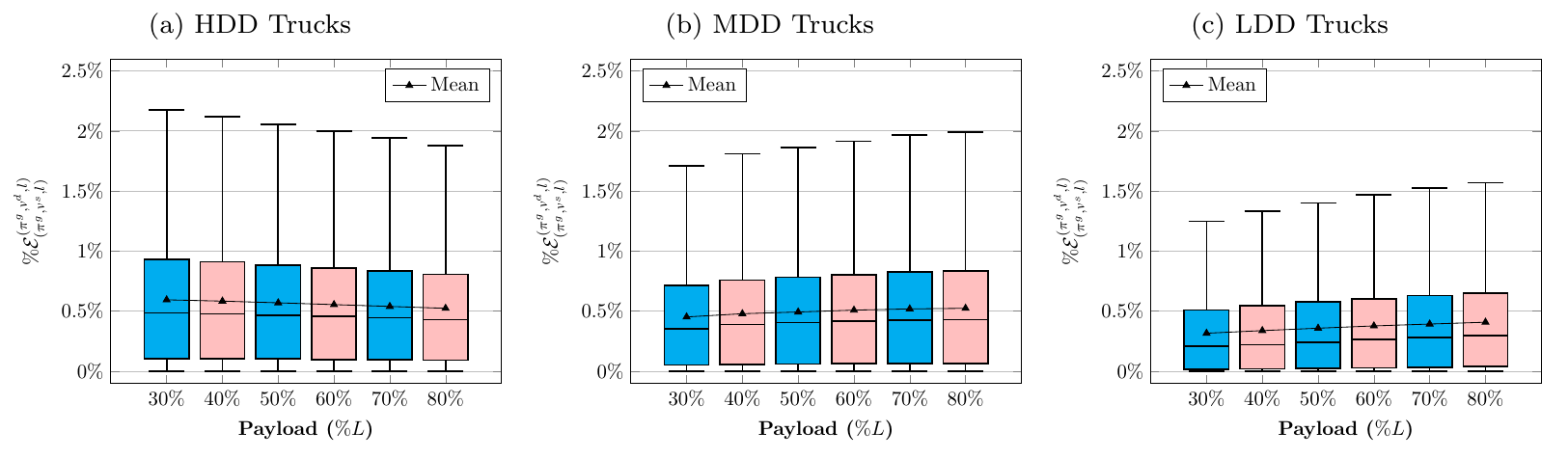}\\
	\caption{\textsf{$\%\mathcal{E}^{(\pi^g, \varv^s, l)}_{(\pi^g, \varv^d, l)}$ across truck types and payloads in traffic condition.}}
  \label{fig:traffic_gvs_gvd}
\end{figure}
This figure demonstrates that the sample mean and third quartiles of $\%\mathcal{E}^{(\pi^g, \varv^d, l)}_{(\pi^g, \varv^s, l)}$ are below 1\% across all truck types and payloads.
This result is primarily because the traffic conditions hinder trucks from utilizing gravity for acceleration on downhill segments, in most instances.
Nevertheless, optimizing the speed on uphill segments can substantially reduce \cotwo emissions.
If the route planner selects the greenest path for traffic speed, $\pi^g(\varv^f,l)$, rather than the fastest path, and the traveling speed is $\varv^d$, the average \cotwo emissions reduction, $\overline{\%\mathcal{E}}^{(\pi^g(\varv^f,l), v^d , l)}_{(\pi^{fp}, v^d , l)}$, ranges from 3.01\% to 7.03\% (see Figure \ref{fig:traffic_fpvd_gvfvd}).
\begin{figure}[htb!]
	\centering
  \includegraphics[width=1\textwidth]{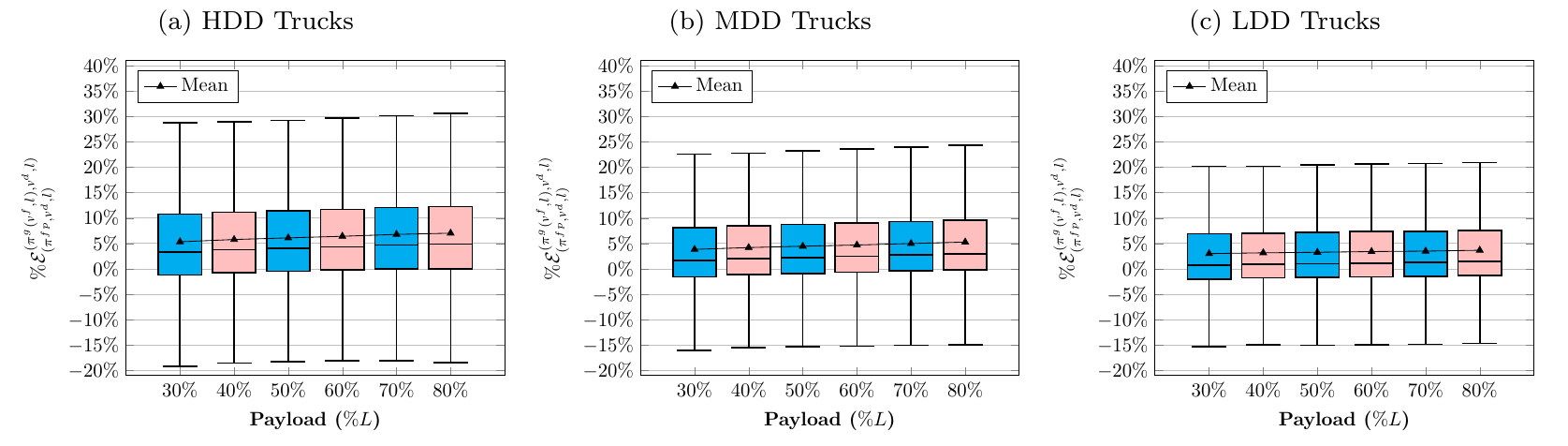}\\
	\caption{\textsf{$\%\mathcal{E}^{(\pi^g(\varv^f,l), v^d , l)}_{(\pi^{fp}, v^d , l)}$ across truck types and payloads in traffic condition.}}
  \label{fig:traffic_fpvd_gvfvd}
\end{figure}
Although $\pi^g(\varv^f,l)$ is not the optimal path for minimizing \cotwo emissions when $\varv^d$ is decided, a comparison of Figures \ref{fig:traffic_fpvd_gvd} and \ref{fig:traffic_fpvd_gvfvd} reveals that choosing $\pi^g(\varv^f,l)$ instead of $\pi^{fp}$ can achieve more than half of the potential \cotwo emissions reduction in most instances (the same argument holds under $\varv^s$).
It is worth noting that \cotwo reduction potential of $\pi^g(\varv^d,l)$ (or $\pi^g(\varv^s,l)$) over $\pi^g(\varv^f,l)$ is slightly higher for lower payloads (see Figure \ref{fig:traffic_gvfvd_gvd}).
This phenomenon is due to the convergence of the greenest paths to the asymptotic greenest paths.
\begin{figure}[htb!]
	\centering
  \includegraphics[width=1\textwidth]{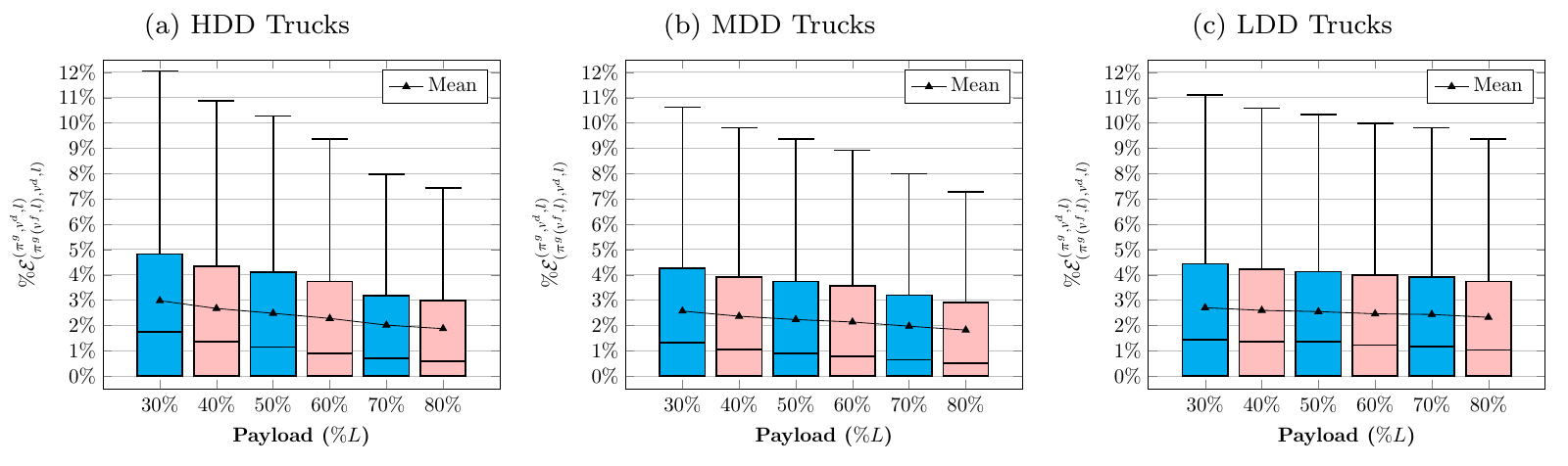}\\
	\caption{\textsf{$\%\mathcal{E}^{(\pi^g, v^d , l)}_{(\pi^g(\varv^f,l), v^d , l)}$ across truck types and payloads in traffic condition.}}
  \label{fig:traffic_gvfvd_gvd}
\end{figure}

\subsection{Results: Increased Travel Duration}
\label{subsec:traveltime_traffic}
Figures \ref{fig:time_fpvf_gvd} presents statistics on the increased travel duration when trucks travel on the greenest paths with $\varv^d$ instead of the fastest path with $\varv^f$.
\begin{figure}[htb!]
	\centering
  \includegraphics[width=1\textwidth]{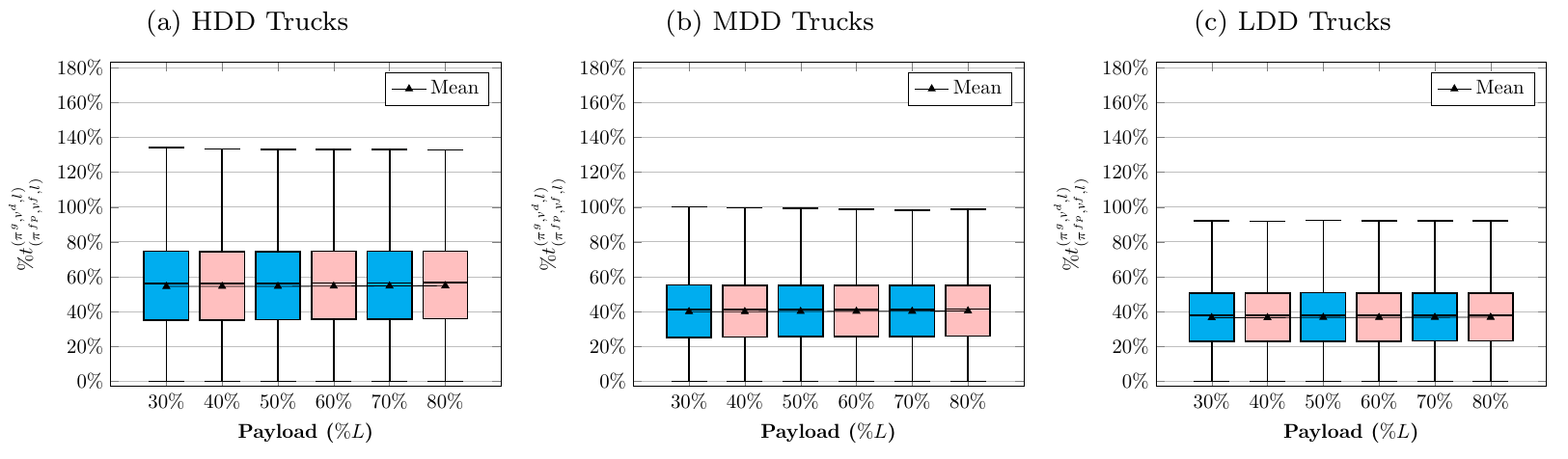}\\
	\caption{\textsf{$\%t^{(\pi^g, v^d , l)}_{(\pi^{fp}, \varv^f, l)}$ across truck types and payloads in traffic condition.}}
  \label{fig:time_fpvf_gvd}
\end{figure}
The figure indicates that, on average, the travel duration increases relative to the fastest path (incorporating traffic speed) from 36.53\% for LDD trucks with 30\% payload to 54.48\% for HDD trucks with 80\% payload.
However, if the speed policy for both greenest path and fastest path is $\varv^d$, the average advantage of selecting the fastest path in terms of travel duration is less than 2.28\% (see Figures \ref{fig:time_fpvd_gvd}).
\begin{figure}[htb!]
	\centering
  \includegraphics[width=1\textwidth]{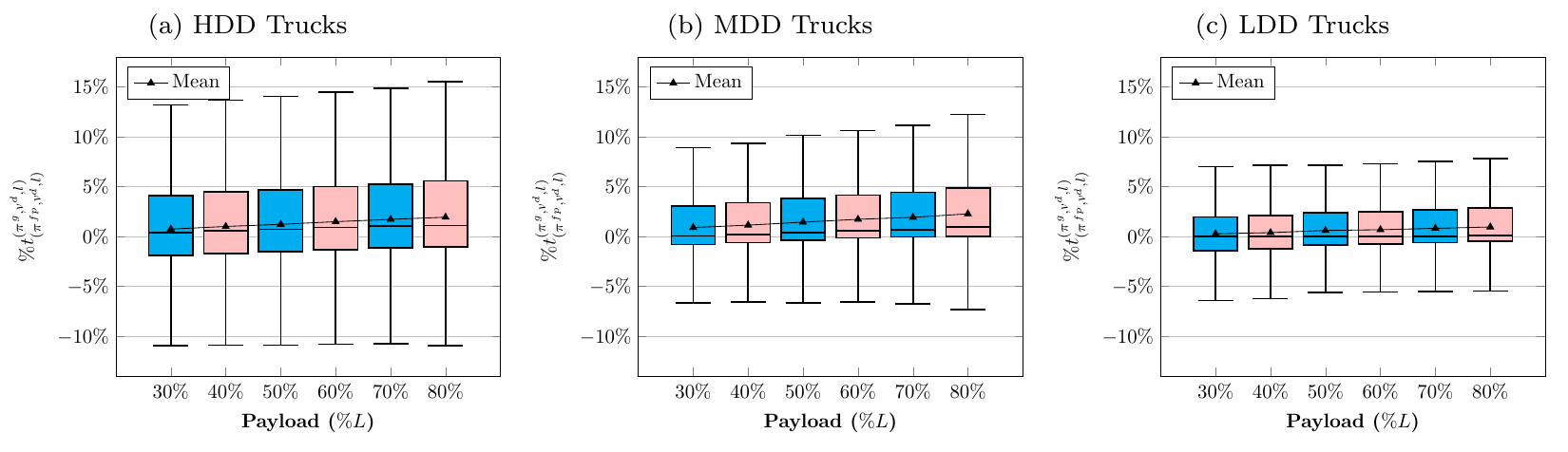}\\
	\caption{\textsf{$\%t^{(\pi^g, v^d , l)}_{(\pi^{fp}, \varv^d, l)}$ across truck types and payloads in traffic condition.}}
  \label{fig:time_fpvd_gvd}
\end{figure}
In several instances, a truck traverses the greenest path even faster than the fastest path when the selected speed policy is $\varv^d$ because the fastest path is found under the assumption of the traffic speed policy.
The statistics presented in Sections \ref{subsec:emissionreduction_traffic} and \ref{subsec:traveltime_traffic} clearly show that when the dynamic speed policy is selected, the \cotwo emissions reduction is larger than the increased travel duration.
This argument holds for the static speed policy, $\varv^s$. The statistics of the static speed policy are presented in Appendix \ref{section:NumericalExperimentTraffic}.

\subsection{Results: Paths of the $\pi^g(\varv^d, l)$, $\pi^g(\varv^s, l)$, $\pi^g(\varv^f, l)$, and $\pi^{fp}$}
\label{subsec:trajectories_traffic}
The dissimilarities between the fastest path and the greenest paths is an important factor in the \cotwo reduction potential of the greenest paths, as highlighted in Section \ref{subsec:emissionreduction_traffic}.
Figures \ref{fig:path_fp_gvd} to \ref{fig:path_fp_gvs} present the statistics for these dissimilarities.
\begin{figure}[htb!]
	\centering
  \includegraphics[width=1\textwidth]{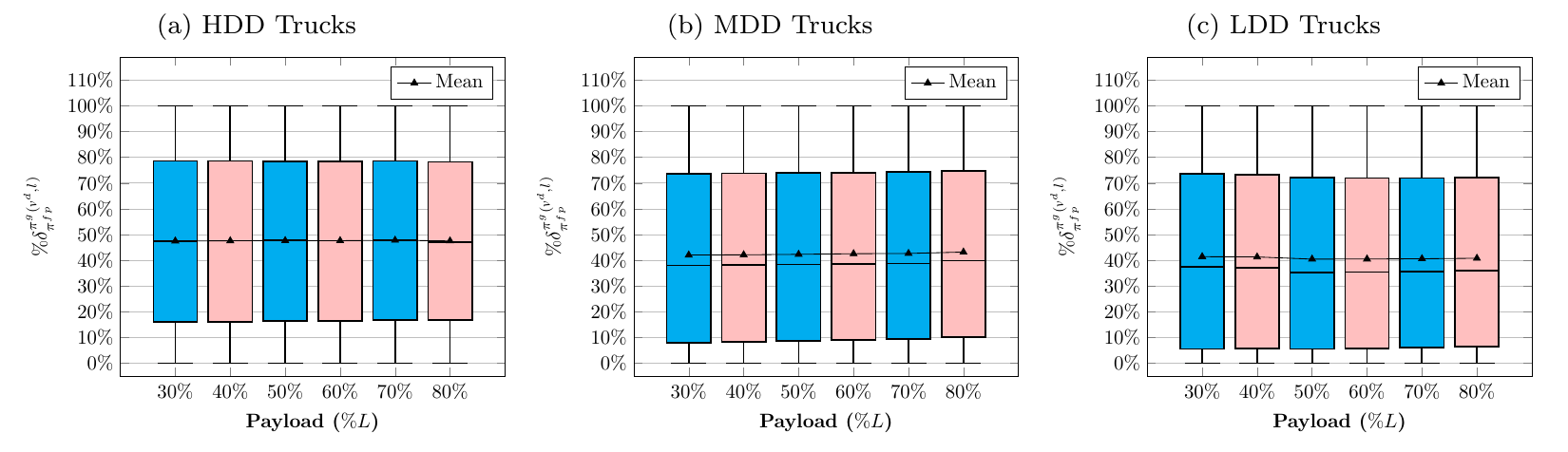}\\
	\caption{\textsf{$\%\delta^{\pi^g(\varv^d,l)}_{\pi^{fp}}$ across truck types and payloads in traffic condition.}}
  \label{fig:path_fp_gvd}
\end{figure}
\begin{figure}[htb!]
	\centering
  \includegraphics[width=1\textwidth]{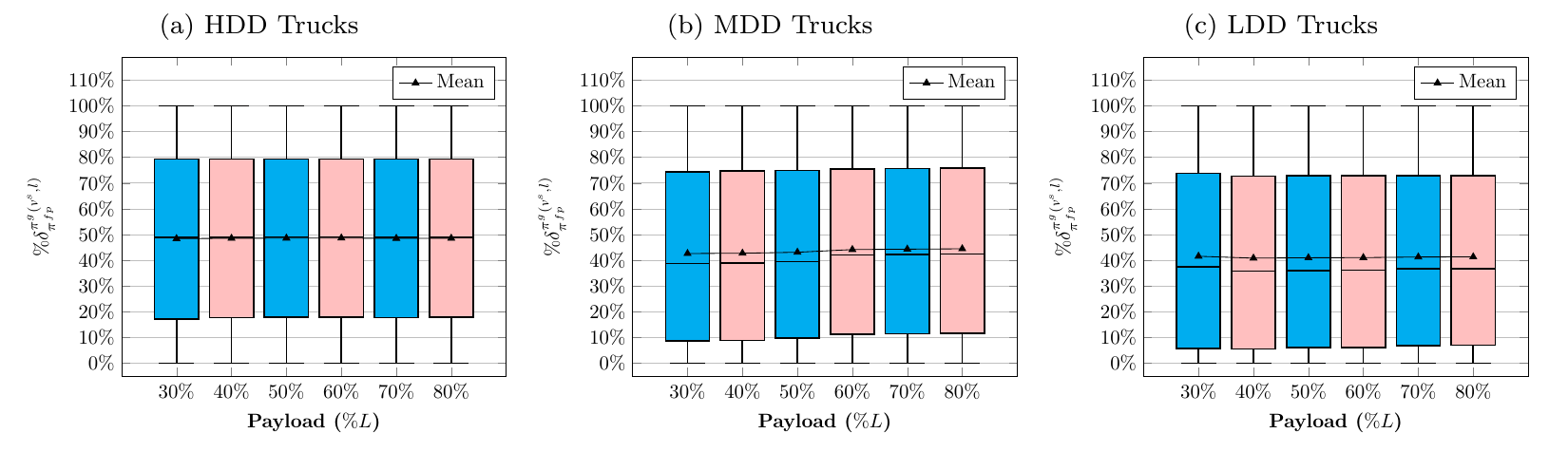}\\
	\caption{\textsf{$\%\delta^{\pi^g(\varv^s,l)}_{\pi^{fp}}$ across truck types and payloads in traffic condition.}}
  \label{fig:path_fp_gvs}
\end{figure}
Figure \ref{fig:path_fp_gvd} shows that, on average, $\pi^{fp}$ does not share 40.47\% to 47.81\% of its paths with $\pi^g(\varv^d, l)$.
A similar statistic for $\pi^g(\varv^s, l)$, i.e., $\overline{\%\delta}^{\pi^g(\varv^s, l)}_{\pi^{fp}}$, ranges from 41.01\% to 48.62\%, which is slightly higher than that of $\overline{\%\delta}^{\pi^g(\varv^d, l)}_{\pi^{fp}}$ (see Figure \ref{fig:path_fp_gvs}).
We observed a similar pattern in Section \ref{subsec:trajectories} (see Figures \ref{fig:baseLength_new} and \ref{fig:baseLength_old}), where we compared the greenest paths and the shortest path in the free flow conditions.
Regarding the dissimilarity between $\pi^g(\varv^d, l)$ and $\pi^g(\varv^s, l)$, our experiments show that $\%\delta^{\pi^g(\varv^s, l)}_{\pi^g(\varv^d, l)}$ is zero or negligible for the majority of instances.
Figure \ref{fig:path_gvd_gvs} indicates that the first, second, and third quartiles, as well as the upper whisker of $\%\delta^{\pi^g(\varv^s, l)}_{\pi^g(\varv^d, l)}$, are zero, and the maximum $\overline{\%\delta}^{\pi^g(\varv^s, l)}_{\pi^g(\varv^d, l)}$ is 2.62\% (cf. Figure \ref{fig:baseLength_GPP}).
\begin{figure}[htb!]
	\centering
  \includegraphics[width=1\textwidth]{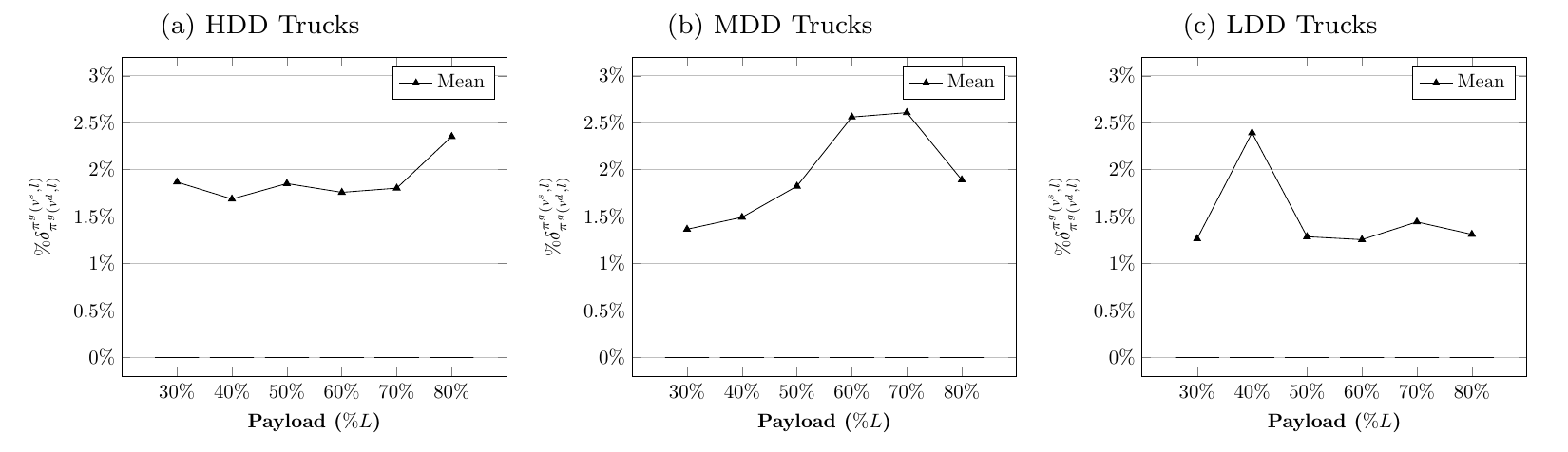}\\
	\caption{\textsf{$\%\delta^{\pi^g(\varv^s,l)}_{\pi^g(\varv^d,l)}$ across truck types and payloads in traffic condition.}}
  \label{fig:path_gvd_gvs}
\end{figure}
As mentioned in Section \ref{subsec:emissionreduction_traffic} this result stems from the limitations on maximum speed on the downhill arcs due to traffic conditions.

\subsection{Results: Asymptotic Greenest Path under Traffic}
\label{subsec:AP_num_traffic}
In Section \ref{subsec:AP_num}, we explained that the convergence of greenest paths to the asymptotic greenest paths is observable for all truck types as the payload increases.
This section examines the asymptotic greenest paths under traffic conditions.
We focus exclusively on $\pi^{\infty}(\varv^g)$, since similar trends can be expected for $\pi^{\infty}(\varv^s)$, as discussed in Section \ref{subsec:trajectories_traffic}.
Figure \ref{fig:traffic_fpvd_avd} illustrates the average \cotwo reduction achieved by the asymptotic greenest path compared to the fastest path with the same dynamic speed policy.
\begin{figure}[htb!]
	\centering
  \includegraphics[width=1\textwidth]{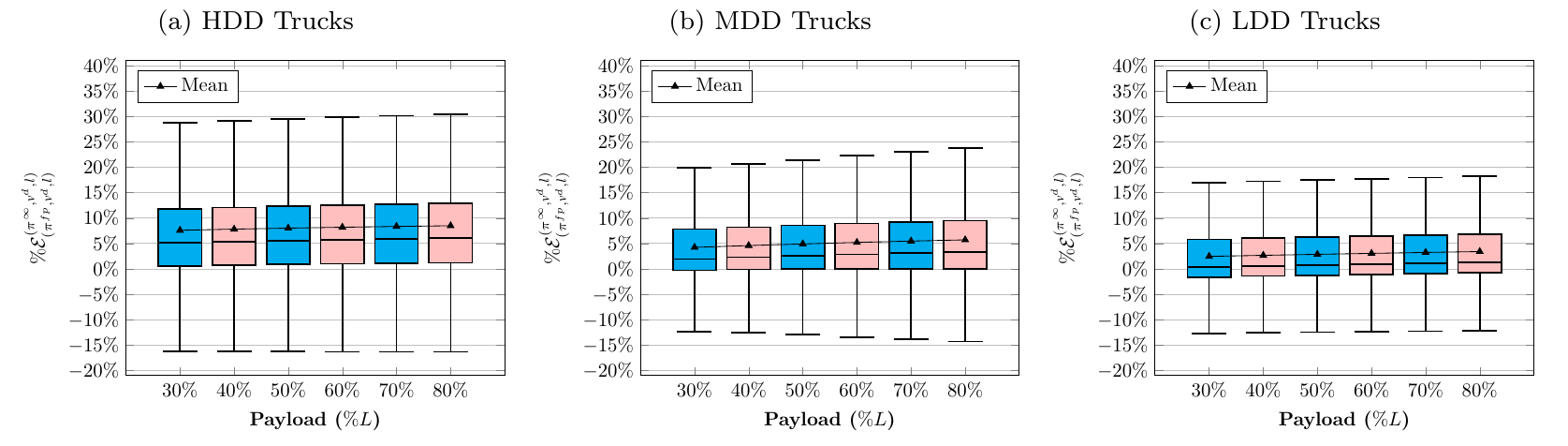}\\
	\caption{\textsf{$\%\mathcal{E}^{(\pi^{\infty}, v^d , l)}_{(\pi^{fp}, \varv^d, l)}$ across truck types and payloads in traffic condition.}}
  \label{fig:traffic_fpvd_avd}
\end{figure}
Specifically, $\overline{\%\mathcal{E}}^{(\pi^{\infty}, \varv^d, l)}_{(\pi^{fp}, \varv^d, l)}$ ranges from 2.47\% for LLD trucks with 30\% payload to 8.50\% for HDD trucks with 80\% payload.
This result indicates that the \cotwo reduction potential of $\pi^{\infty}(\varv^g)$ is similar to that of $\pi^g(\varv^f, l)$, even slightly higher for MDD and HDD trucks (cf. Figure \ref{fig:traffic_fpvd_gvfvd}).
Figure \ref{fig:traffic_gvd_avd} shows the \cotwo emissions reduction of the asymptotic greenest path relative to the greenest path, $\%\mathcal{E}^{\pi^{\infty}, \varv^d, l}_{\pi^g, \varv^d, l}$.
\begin{figure}[htb!]
	\centering
  \includegraphics[width=1\textwidth]{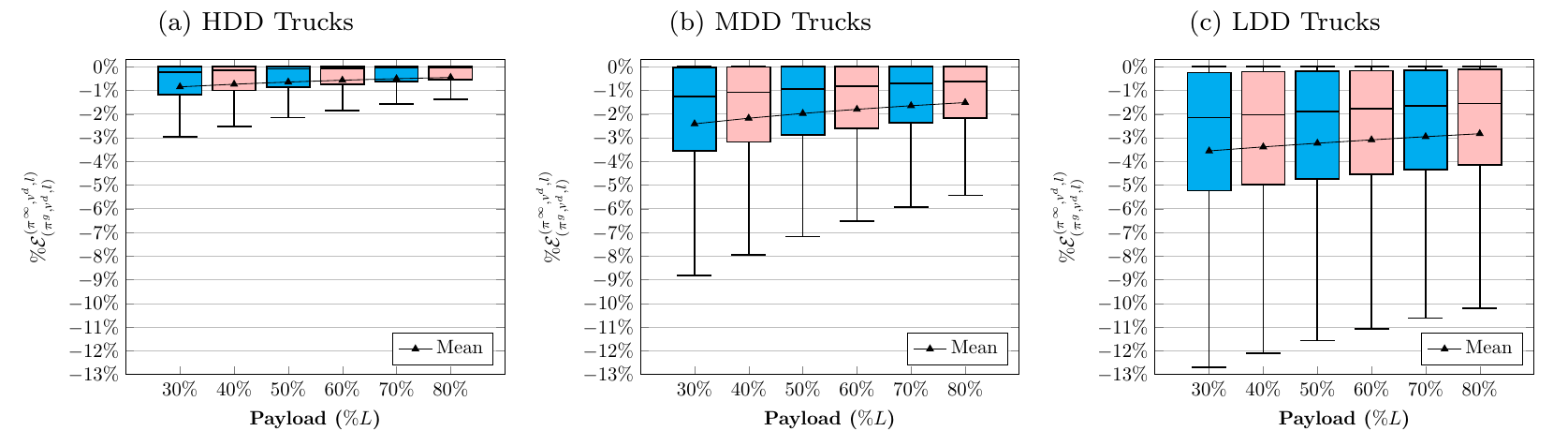}\\
	\caption{\textsf{$\%\mathcal{E}^{(\pi^{\infty}, v^d , l)}_{(\pi^g, \varv^d, l)}$ across truck types and payloads in traffic condition.}}
  \label{fig:traffic_gvd_avd}
\end{figure}
It is straightforward to see that $\overline{\%\mathcal{E}}^{\pi^{\infty}, \varv^d, l}_{\pi^g, \varv^d, l}$ increases with truck weight, rising from -3.56\% for LLD trucks with 30\% payload to -0.45\% for HDD trucks with 80\% payload.
Figure \ref{fig:traffic_path_gvd_avd} highlights the convergence of the greenest paths to the asymptotic greenest paths, similar to the tendency observed in Section \ref{subsec:AP_num} for free flow conditions.
\begin{figure}[htb!]
	\centering
  \includegraphics[width=1\textwidth]{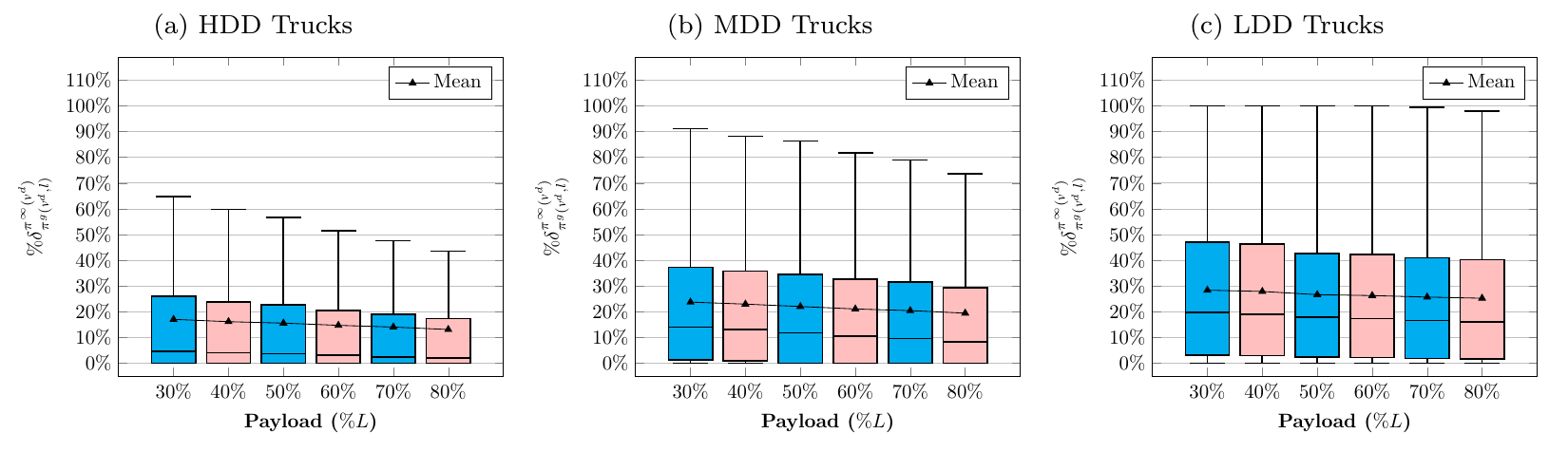}\\
	\caption{\textsf{$\%\delta^{\pi^{\infty}(\varv^d)}_{\pi^g(\varv^d,l)}$ across truck types and payloads in traffic condition.}}
  \label{fig:traffic_path_gvd_avd}
\end{figure}
Comparing the convergence results under free flow and traffic conditions in New York city, we can infer that, on average, convergence under traffic occurs more rapidly than free-flow condition.
This phenomenon is due to the limitations on speed choices in traffic. 
In such a case, the selection of paths with less ascent, i.e. $h^\prime(a)$, $a \in A$, plays a crucial role in the reduction of emissions \cotwo.

\subsection{Conclusions on the role of traffic information}

The incorparation of slope information to find the greenest path is even more important in heavy traffic than it is in free flow traffic. 
Comparing the results of dynamic and static speed optimizations in traffic, we find that optimizing speed on uphill arcs can significantly reduce \cotwo emissions.
However, using gravity to accelerate on downhill arcs is limited due to traffic congestion. From a public policy perspective, this finding reinforces arguments for scheduling truck deliveries during ``Off-Hours" when more environmentally-friendly options for path and speed selection are available \citep[see e.g.][]{offhour}.
Truck delivery during ``Off-Hours'' can, additionally, reduce traffic congestion when trucks traverse uphill roads with optimized speeds that may be lower than the traffic speed.
In Section \ref{subsec:AP_num_traffic}, we demonstrate that the greenest path converges to the asymptotic greenest path even faster under traffic conditions than in free flow conditions. 

\end{rev}

\begin{rev2}
    \section{Summary of Key Findings from Numerical Experiments}
    
    The first outcome of our experiments in Sections \ref{section:NumericalExperiments} and \ref{section:NumericalExperimentTraffic} is that high-resolution topographical data should be incorporated into urban truck transportation decisions when minimizing \cotwo emissions is the objective. Specifically, pre-computation of the greenest paths is not feasible due to the non-linear effects of speed decisions, road gradients, and payload. A similar argument has previously been made regarding the need to integrate high-resolution traffic speed data into emissions-minimizing transportation decisions \citep[see e.g.][]{Ehmke2016b}.
    
    Secondly, our results show that optimal speed decisions are dynamic, with dynamic speed choices reducing \cotwo emissions by 2\% to 4\% in free-flow conditions compared to static speed choices. While the difference between dynamic and static speed decisions is less significant in traffic, we found that optimized speeds still achieve significantly lower emissions than traffic speeds, even when acceleration is restricted by traffic congestion.

    Thirdly, we observed that the greenest path is relatively insensitive to whether speed decisions are static or dynamic, even in free-flow conditions. Additionally, the greenest path begins to converge to the asymptotic greenest path at low payload ratios under both free-flow and traffic conditions. Therefore, a pre-computed greenest path for a given speed decision (e.g., static) and payload level (e.g., 50\% or 100\%) can be a good approximation for the greenest paths across different speed decisions and payloads. This approximation can help reduce the computational complexity of green transportation problems like PRP.
\end{rev2}

\section{Conclusions}
\label{section:conclusion}
In this paper, we studied the greenest path selection problem for a logistics service provider that operates a fleet of heavy-, medium-, and light-duty trucks in an urban environment. 
We established that the policies for the speed and path that minimize \cotwo emissions are slope-dependent (dynamic). We also showed that the greenest path converges to a fixed path as the payload increases and provided an efficient algorithm to compute the asymptotic greenest path.
We conducted extensive numerical experiments using elevation data of 25 cities around the world to investigate the potential \cotwo reduction by such dynamic policies \begin{rev}under free flow traffic conditions.
\end{rev}
The results in section \ref{subsec:emissionreduction} showed that, on average, the combined dynamic path and speed selection can reduce \cotwo emissions by 1.19\% to 10.15\% based on the truck type and city. 
Our analysis also showed that in most cities, the average emissions reduction potential of dynamic speed optimization lies between 2\% to 4\% regardless of the truck type. Nonetheless, the effect of slope-dependent path selection (the greenest path) depends on the payload and truck type. \begin{rev}\begin{rev2}In section \ref{section:NumericalExperimentTraffic},\end{rev2} for the city of New York we also studied the effect of effective speed limits due to traffic congestion and found that choosing the greenest path can lead to even larger \cotwo reduction than in free flow traffic conditions.
\end{rev}
In section \ref{subsec:trajectories}, we explained that the greenest path significantly differs from the shortest path.
While the greenest path depends on the speed policy, the experiments show that this dependence is weak and that the greenest path under the static speed policy is usually optimal or near optimal \begin{rev}
 especially when speed is determined by traffic.
\end{rev}
Moreover, we demonstrated, in Sections \ref{subsect:asymcal} and \ref{subsec:AP_num}, that the greenest path diverges from the shortest path as the payload increases and converges to the asymptotic greenest path, i.e. the greenest path for the arbitrary large payloads. \begin{rev}
    Convergence to the asymptotic greenest path is faster under heavy traffic conditions.
\end{rev}
These results could be used for the approximation of the greenest path to simplify complex transportation problems.
The analysis of variance (ANOVA) indicated that the potential \cotwo emissions reduction by the greenest path and the dynamic speed policy is associated positively with the variability of arc gradients along the shortest path, and negatively with the relative elevation of the source and target. 

\ACKNOWLEDGMENT{%
We thank Dennis Davydov for sharing initial explorations on this topic and Tiffany Nguyen for extensive feedback on early drafts.
}

\section*{Declarations}
The authors did not receive support from any organization for the submitted work.
The authors have no relevant financial or non-financial interests to disclose.

\bibliographystyle{informs2014} 
\bibliography{bib} 

%
%
%

\newpage
\begin{APPENDICES}

\section{ANOVA Results}
\label{section:ANOVA}
\vspace{1cm}
\begin{table*}[!htbp]
  \fontsize{8pt}{9pt}\selectfont
  \scriptsize
  \parbox[t]{.45\linewidth}{
  \centering
  \caption{\textsf{Response: $\%\mathcal{E}_{(\pi^{sp}, \varv^s, l)}^{(\pi^g, \varv^d, l)}$}}
  \label{tab:ENO}
  \begin{tabular}{l r r r r r}
        \toprule
        $Feature$ & $df$ & SS & $MSS$ & $F-value$ & $p-value$ \\ \midrule
        $\sigma^{sp} (\theta)$ & 1 & 11771 & 11771 & 9061762 & $<10^{-15}$ \\ 
        $\Delta h$ & 1 & 22347 & 22347 & 17203995 & $<10^{-15}$ \\ 
        $l$ & 1 & 293 & 293 & 225597 & $<10^{-15}$ \\ 
        $\delta^{sp}$ & 1 & 47 & 47 & 36175 & $<10^{-15}$ \\ 
        City & 24 & 12477 & 520 & 400209 & $<10^{-15}$ \\ 
        Truck & 2 & 926 & 463 & 356459 & $<10^{-15}$ \\ 
        (Intercept) & 1 & 7265 & 7265 & 5592723 & $<10^{-15}$ \\ 
        Residuals & 55415417 & 71983 & ~ & ~ & ~ \\ \bottomrule
    \end{tabular}
  }
  \hfill
  \parbox[t]{.45\linewidth}{
  \centering
  \caption{\textsf{Response: $\%\mathcal{E}_{(\pi^{sp}, \varv^s, l)}^{(\pi^g, \varv^s, l)}$}}
  \label{tab:EOO}
  \begin{tabular}{l r r r r r}
        \toprule
        $Feature$ & $df$ & SS & $MSS$ & $F-value$ & $p-value$ \\ \midrule
        $\sigma^{sp} (\theta)$ & 1 & 10020 & 10020 & 5973305 & $<10^{-15}$ \\ 
        $\Delta h$ & 1 & 6261 & 6261 & 3732651 & $<10^{-15}$ \\ 
        $l$ & 1 & 908 & 908 & 541007 & $<10^{-15}$ \\ 
        $\delta^{sp}$ & 1 & 225 & 225 & 134066 & $<10^{-15}$ \\ 
        City & 24 & 16178 & 674 & 401852 & $<10^{-15}$ \\ 
        Truck & 2 & 1156 & 578 & 344674 & $<10^{-15}$ \\ 
        (Intercept) & 1 & 420 & 420 & 250138 & $<10^{-15}$ \\ 
        Residuals & 55415417 & 92958 & ~ & ~ & ~ \\ \bottomrule
    \end{tabular}
  }
\end{table*}

\begin{table*}[!htbp]
  \fontsize{8pt}{9pt}\selectfont
  \scriptsize
  \parbox[t]{.45\linewidth}{
  \centering
  \caption{\textsf{Response: $\%\mathcal{E}_{(\pi^{sp}, \varv^d, l)}^{(\pi^g, \varv^d, l)}$}}
  \label{tab:ENN}
  \begin{tabular}{l r r r r r}
        \toprule
        $Feature$ & $df$ & SS & $MSS$ & $F-value$ & $p-value$ \\ \midrule
        $\sigma^{sp} (\theta)$ & 1 & 7092 & 7092 & 6099148 & $<10^{-15}$ \\ 
        $\Delta h$ & 1 & 4822 & 4822 & 4147274 & $<10^{-15}$ \\ 
        $l$ & 1 & 816 & 816 & 701716 & $<10^{-15}$ \\ 
        $\delta^{sp}$ & 1 & 82 & 82 & 70510 & $<10^{-15}$ \\ 
        City & 24 & 10056 & 419 & 360351 & $<10^{-15}$ \\ 
        Truck & 2 & 918 & 459 & 394610 & $<10^{-15}$ \\ 
        (Intercept) & 1 & 347 & 347 & 298743 & $<10^{-15}$ \\ 
        Residuals & 55415417 & 64434 & ~ & ~ & ~ \\ \bottomrule
    \end{tabular}
  }
  \hfill
  \parbox[t]{.45\linewidth}{
  \centering
  \caption{\textsf{Response: $\%\mathcal{E}_{(\pi^g, \varv^s, l)}^{(\pi^g, \varv^d, l)}$}}
  \label{tab:EGP}
  \begin{tabular}{l r r r r r}
        \toprule
        $Feature$ & $df$ & SS & $MSS$ & $F-value$ & $p-value$ \\ \midrule
        $\sigma^{sp} (\theta)$ & 1 & 1431 & 1431 & 3270532 & $<10^{-15}$ \\ 
        $\Delta h$ & 1 & 16335 & 16335 & 37344012 & $<10^{-15}$ \\ 
        $l$ & 1 & 147 & 147 & 335105 & $<10^{-15}$ \\ 
        $\delta^{sp}$ & 1 & 870 & 870 & 1989427 & $<10^{-15}$ \\ 
        City & 24 & 3102 & 129 & 295466 & $<10^{-15}$ \\ 
        Truck & 2 & 81 & 40 & 91981 & $<10^{-15}$ \\ 
        (Intercept) & 1 & 10731 & 10731 & 24531295 & $<10^{-15}$ \\ 
        Residuals & 55415417 & 24240 & ~ & ~ & ~ \\ \bottomrule
    \end{tabular}
  }
\end{table*}

\vspace{2em}

\begin{table*}[!htbp]
  \fontsize{8pt}{9pt}\selectfont
  \scriptsize
  \parbox[t]{.45\linewidth}{
  \centering
  \caption{\textsf{Response: $\%\delta_{\pi^{sp}}^{\pi^g (\varv^s, l)}$}}
  \label{tab:RO}
  \begin{tabular}{l r r r r r}
        \toprule
        $Feature$ & $df$ & SS & $MSS$ & $F-value$ & $p-value$ \\ \midrule
        $\sigma^{sp} (\theta)$ & 1 & 204303 & 204303 & 3344542 & $<10^{-15}$ \\ 
        $\Delta h$ & 1 & 15691 & 15691 & 256873 & $<10^{-15}$ \\ 
        $l$ & 1 & 5270 & 5270 & 86270 & $<10^{-15}$ \\ 
        $\delta^{sp}$ & 1 & 44990 & 44990 & 736513 & $<10^{-15}$ \\ 
        City & 24 & 640819 & 26701 & 437105 & $<10^{-15}$ \\ 
        Truck & 2 & 13408 & 6704 & 109750 & $<10^{-15}$ \\ 
        (Intercept) & 1 & 2732 & 2732 & 44716 & $<10^{-15}$ \\ 
        Residuals & 55415417 & 3385080 & ~ & ~ & ~ \\ \bottomrule
    \end{tabular}
  }
  \hfill
  \parbox[t]{.45\linewidth}{
  \centering
  \caption{\textsf{Response: $\%\delta_{\pi^{sp}}^{\pi^g (\varv^d, l)}$}}
  \label{tab:RN}
  \begin{tabular}{l r r r r r}
        \toprule
        $Feature$ & $df$ & SS & $MSS$ & $F-value$ & $p-value$ \\ \midrule
        $\sigma^{sp} (\theta)$ & 1 & 208923 & 208923 & 3601979 & $<10^{-15}$ \\ 
        $\Delta h$ & 1 & 17790 & 17790 & 306719 & $<10^{-15}$ \\ 
        $l$ & 1 & 7602 & 7602 & 131065 & $<10^{-15}$ \\ 
        $\delta^{sp}$ & 1 & 37788 & 37788 & 651499 & $<10^{-15}$ \\ 
        City & 24 & 577881 & 24078 & 415128 & $<10^{-15}$ \\ 
        Truck & 2 & 18821 & 9411 & 162244 & $<10^{-15}$ \\ 
        (Intercept) & 1 & 3005 & 3005 & 51805 & $<10^{-15}$ \\ 
        Residuals & 55415417 & 3214225 & ~ & ~ & ~ \\ \bottomrule
    \end{tabular}
  }
\end{table*}

\vspace{2em}

\begin{table*}[!htbp]
  \fontsize{8pt}{9pt}\selectfont
  \scriptsize
  \parbox[t]{.48\linewidth}{
  \centering
  \caption{\textsf{Response: $\%\delta_{(\pi^g, \varv^d, l)}^{(\pi^g, \varv^s, l)}$}}
  \label{tab:RGP}
  \begin{tabular}{l r r r r r}
        \toprule
        $Feature$ & $df$ & SS & $MSS$ & $F-value$ & $p-value$ \\ \midrule
        $\sigma^{sp} (\theta)$ & 1 & 1190 & 1190 & 74122 & $<10^{-15}$ \\ 
        $\Delta h$ & 1 & 596 & 596 & 37147 & $<10^{-15}$ \\ 
        $l$ & 1 & 390 & 390 & 24271 & $<10^{-15}$ \\ 
        $\delta^{sp}$ & 1 & 1096 & 1096 & 68270 & $<10^{-15}$ \\ 
        City & 24 & 22738 & 947 & 59004 & $<10^{-15}$ \\ 
        Truck & 2 & 310 & 155 & 9667 & $<10^{-15}$ \\ 
        (Intercept) & 1 & 203 & 203 & 12647 & $<10^{-15}$ \\ 
        Residuals & 55415417 & 889779 & ~ & ~ & ~ \\ \bottomrule
    \end{tabular}
  }
\end{table*}

\newpage
\section{Results: Performance of the Asymptotic Paths with Static Speed policies}
\label{section:APC_Performance}

\begin{figure}[htb!]
	\centering
  \includegraphics[width=1\textwidth]{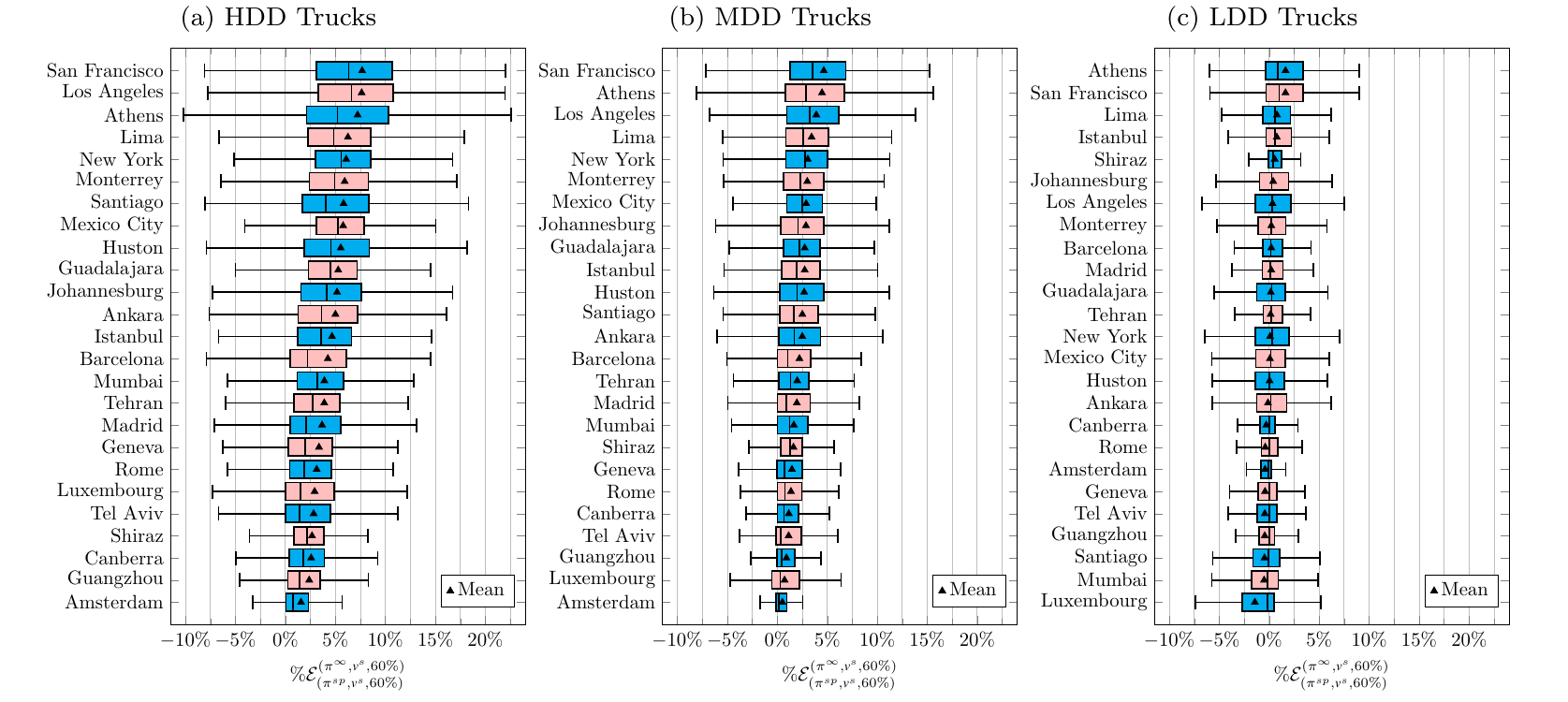}\\
	\caption{\textsf{Relative \cotwo emissions reduction by selecting $(\pi^{\infty}, \varv^s, 60\%)$ rather than $(\pi^{sp}, \varv^s, 60\%)$.}}
  \label{fig:relFCSPC}
\end{figure}
\begin{figure}[htb!]
	\centering
  \includegraphics[width=1\textwidth]{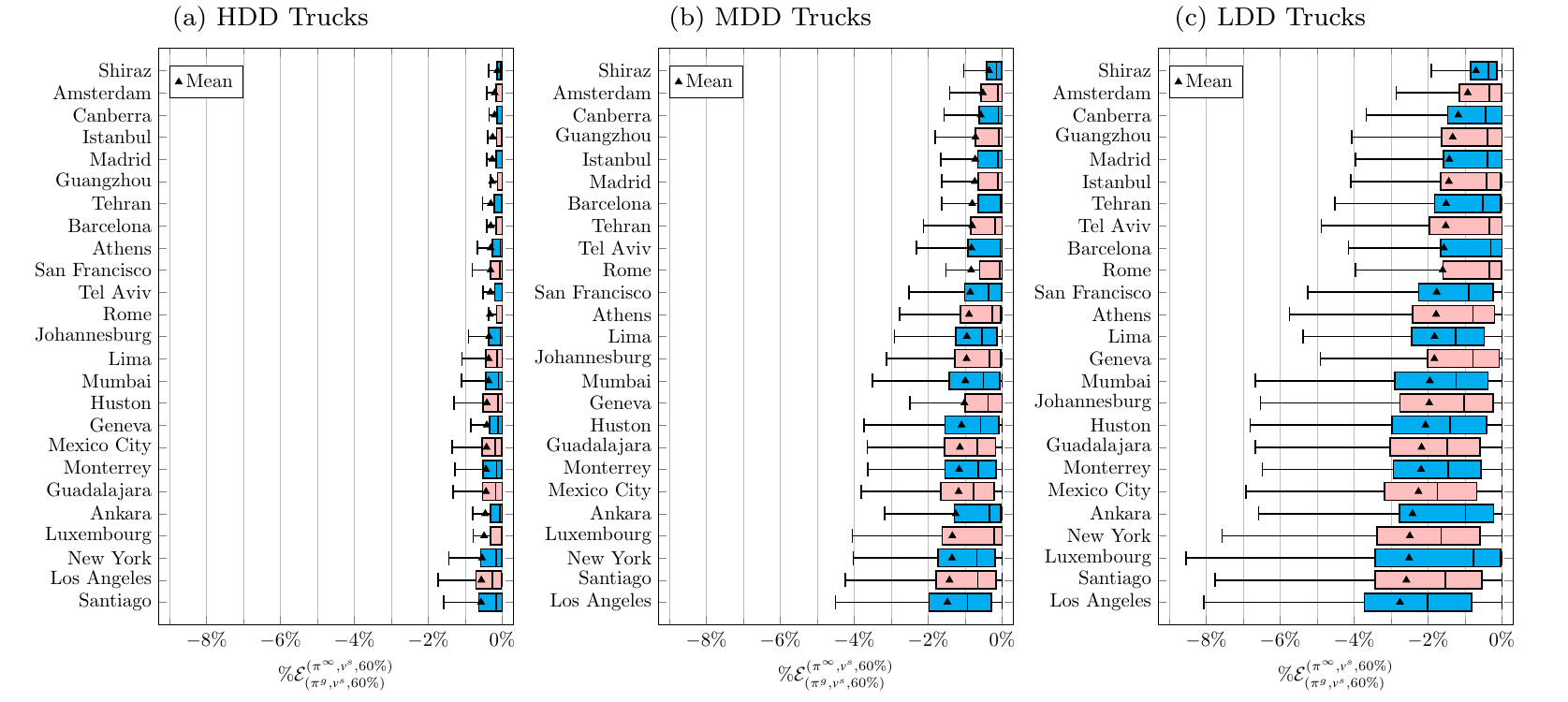}\\
	\caption{\textsf{Relative \cotwo emissions reduction by selecting $(\pi^{\infty}, \varv^s, 60\%)$ rather than $(\pi^g, \varv^s, 60\%)$.}}
  \label{fig:relFCGPC}
\end{figure}
\begin{figure}[htb!]
	\centering
  \includegraphics[width=1\textwidth]{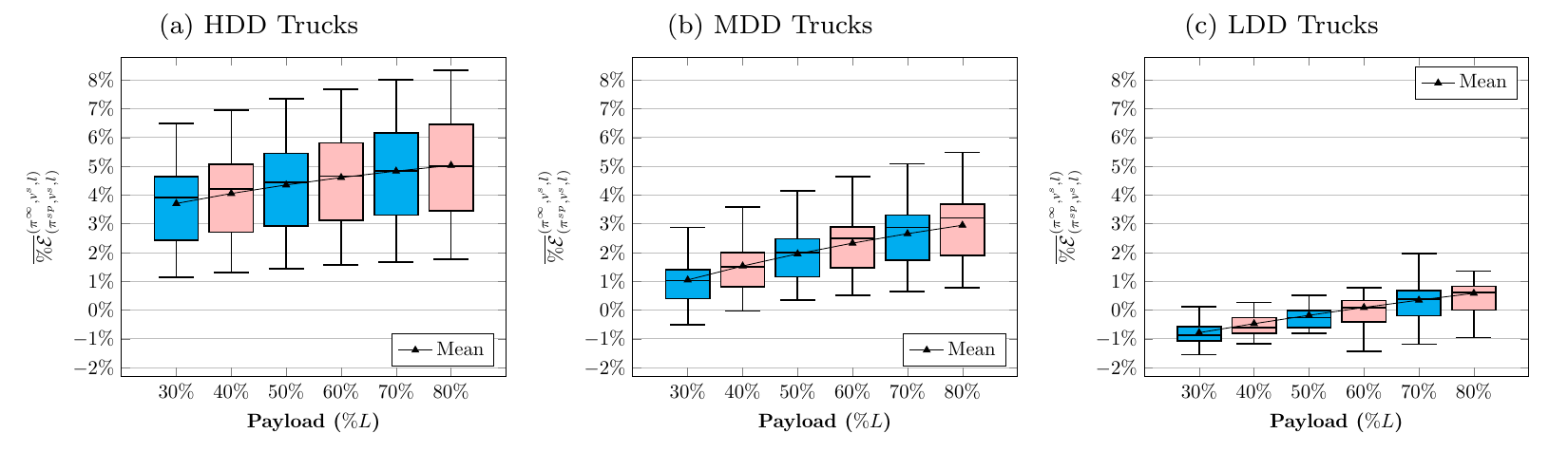}\\
	\caption{\textsf{Effect of payload on $\overline{\%\mathcal{E}}_{(\pi^{sp}, \varv^s, l)}^{(\pi^{\infty}, \varv^s, l)}$ across 25 cities.}}
  \label{fig:varLrelFCSPC}
\end{figure}
\begin{figure}[htb!]
	\centering
  \includegraphics[width=1\textwidth]{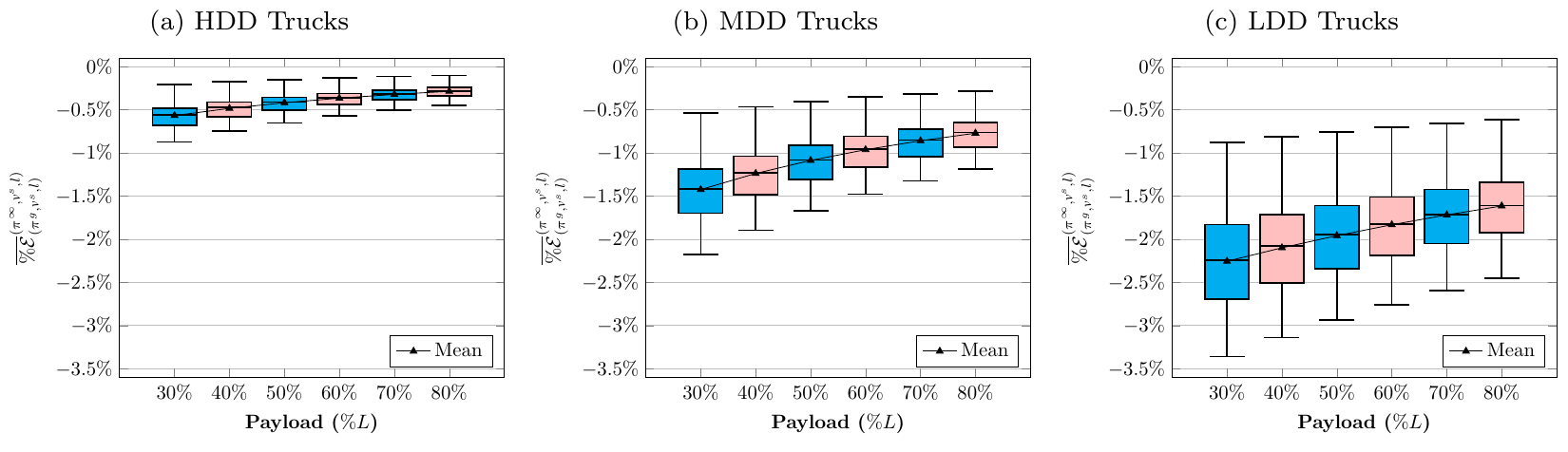}\\
	\caption{\textsf{Effect of payload on $\overline{\%\mathcal{E}}_{(\pi^g, \varv^s, l)}^{(\pi^{\infty}, \varv^s, l)}$ across 25 cities.}}
  \label{fig:varLrelFCGPC}
\end{figure}
\begin{figure}[htb!]
	\centering
  \includegraphics[width=1\textwidth]{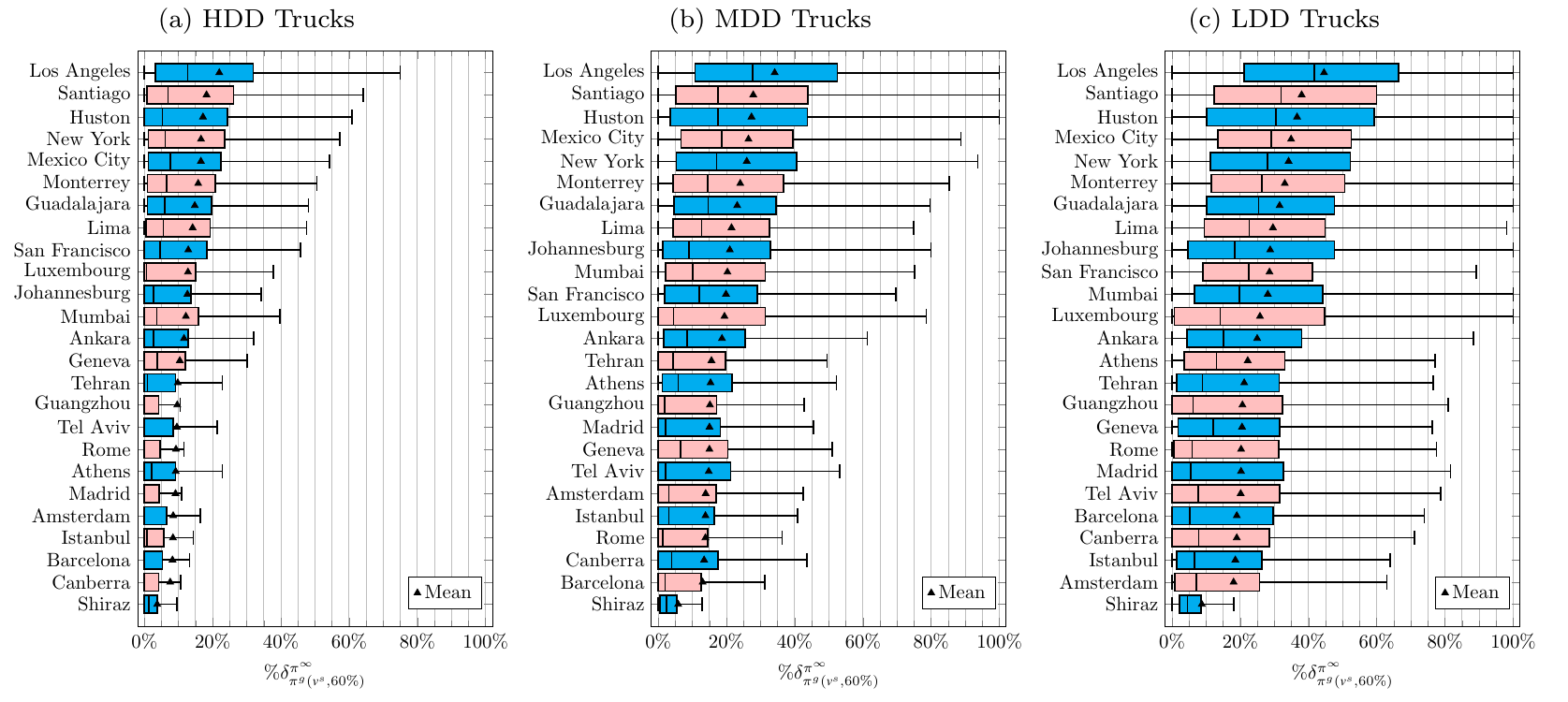}\\
	\caption{\textsf{Ratio of the length of $\pi^g (\varv^s, 60\%)$ that is not shared with $\pi^{\infty}$.}}
  \label{fig:relGPC}
\end{figure}
\begin{figure}[htb!]
	\centering
  \includegraphics[width=1\textwidth]{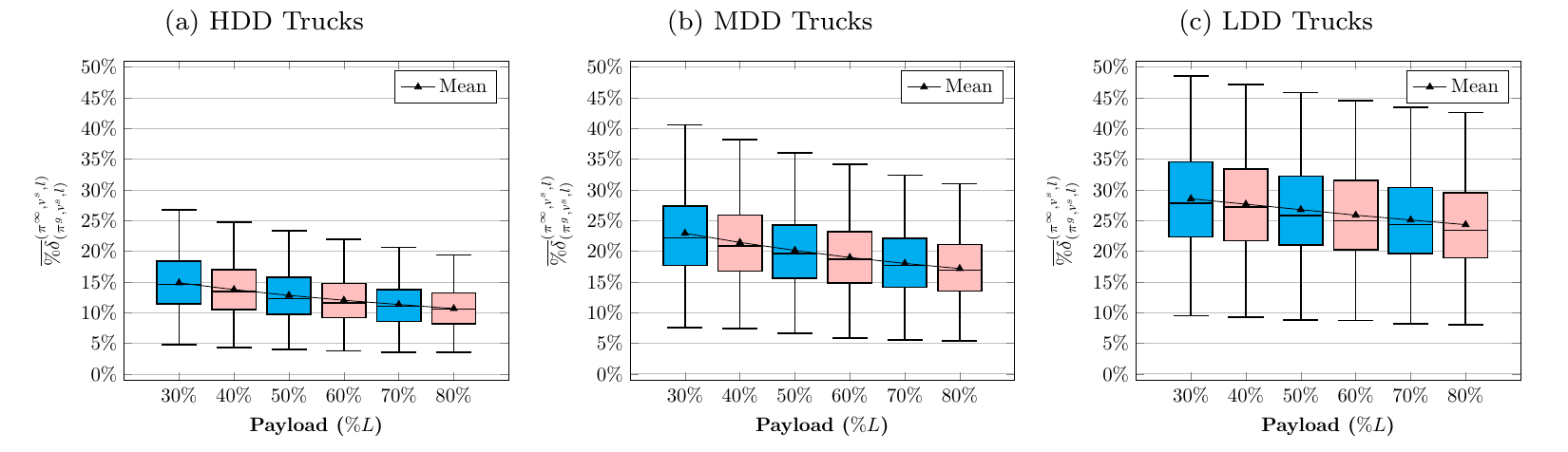}\\
	\caption{\textsf{Effect of payload on $\overline{\%\delta}_{(\pi^g, \varv^s, l)}^{(\pi^{\infty}, \varv^s, l)}$ across 25 cities.}}
  \label{fig:varLrelGPC}
\end{figure}

\newpage
\begin{rev}
\section{Results: Increased Travel Duration under $\varv^s$}
\label{section:traveltime_traffic_appendix}
\begin{figure}[H]
	\centering
  \includegraphics[width=1\textwidth]{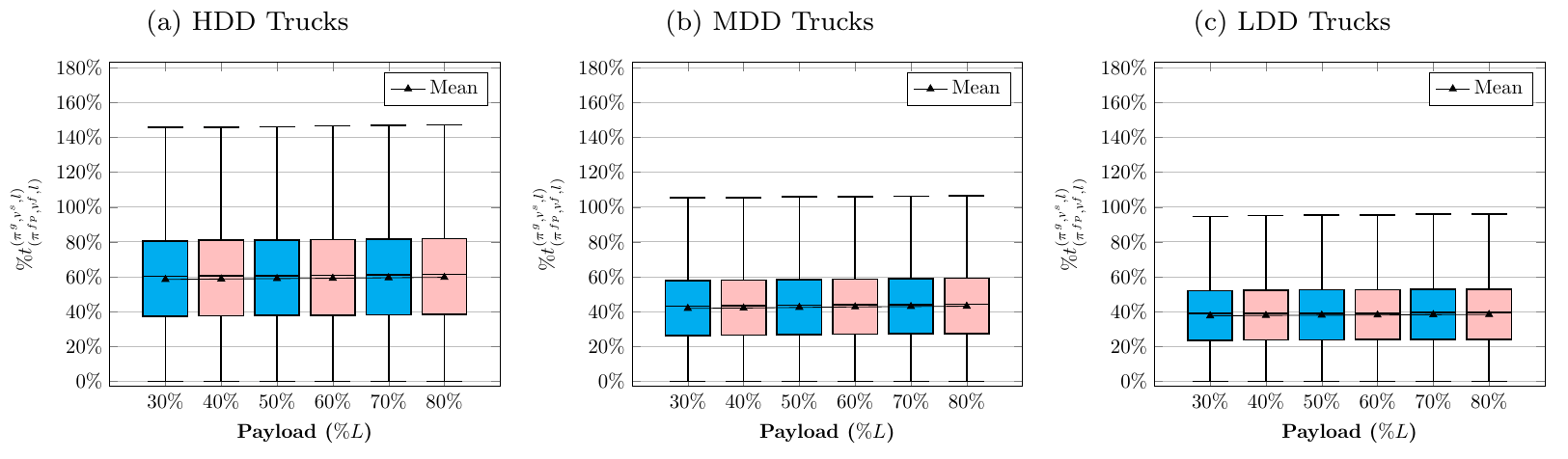}\\
	\caption{\textsf{$\%t^{(\pi^g, v^s , l)}_{(\pi^{fp}, \varv^f, l)}$ across truck types and payload in traffic condition.}}
  \label{fig:time_fpvf_gvs}
\end{figure}
\begin{figure}[H]
	\centering
  \includegraphics[width=1\textwidth]{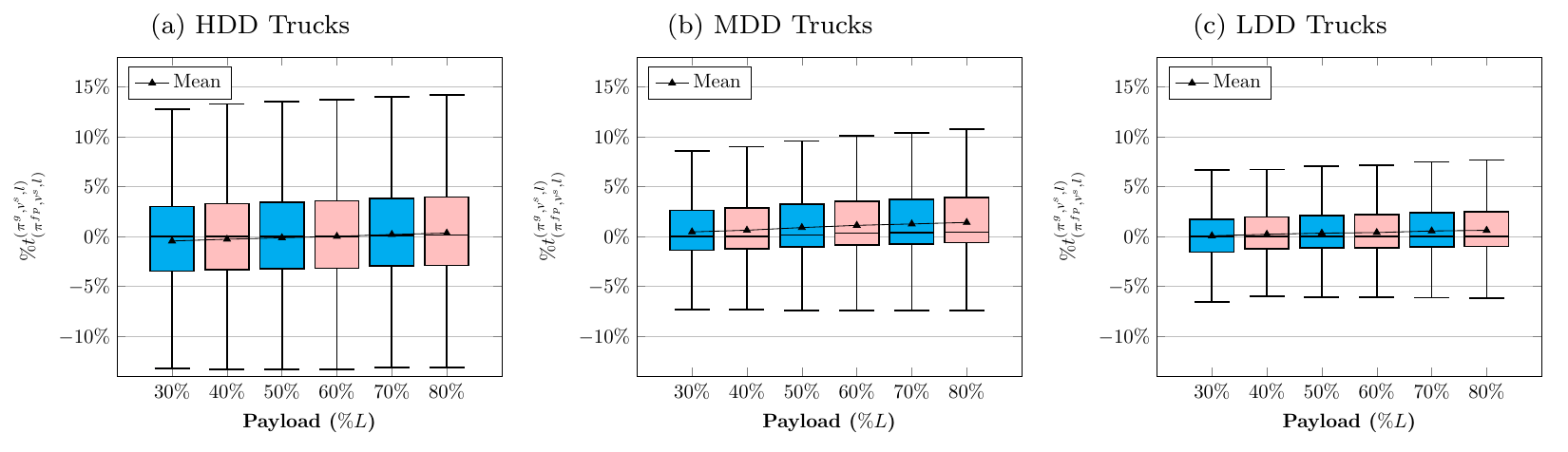}\\
	\caption{\textsf{$\%t^{(\pi^g, v^s , l)}_{(\pi^{fp}, \varv^s, l)}$ across truck types and payloads in traffic condition.}}
  \label{fig:time_fpvs_gvs}
\end{figure}
\end{rev}

\end{APPENDICES}

\end{document}